\documentclass[leqno,11pt]{amsart}
\usepackage{amsmath,amscd,amsthm,amsxtra}
\usepackage{epsfig,graphics,color,colortbl}
\usepackage{amssymb,latexsym}
\usepackage{mathabx}
\usepackage{mathrsfs,eucal,upgreek}
\usepackage[poly,all]{xy}
\usepackage{hyperref}
\usepackage{tikz-cd}
\usepackage{arydshln}
\usepackage{ytableau}
\usepackage{adjustbox}
\usepackage{hyperref}
\usepackage{tikz-cd}
\usepackage{adjustbox}
\usepackage{centernot}
\usepackage{mathtools}
\usepackage{stmaryrd}
\usepackage{caption}
\usepackage{upgreek}

\newtheorem{thm}{\bf Theorem}[section]
\newtheorem{df}[thm]{\bf Definition}
\newtheorem{prop}[thm]{\bf Proposition}
\newtheorem{cor}[thm]{\bf Corollary}
\newtheorem{lem}[thm]{\bf Lemma}
\newtheorem{rem}[thm]{\bf Remark}
\newtheorem{ex}[thm]{\bf Example}

\numberwithin{equation}{section}

\newcommand{\bs}{\boldsymbol}
\newcommand{\B}{\mathbf{B}}

\newcommand{\cP}{\mathscr{P}}
\newcommand{\pf}{\noindent{\bfseries Proof. }}
\newcommand{\ov}{\overline}

\newcommand{\ba}{\bs{\rm a}}
\newcommand{\bb}{\bs{\rm b}}

\newcommand{\bi}{\bs{\rm i}}

\newcommand{\V}{{\bf V}}

\newcommand{\N}{\mathbb{N}}

\newcommand{\Z}{\mathbb{Z}}

\newcommand{\te}{\widetilde{e}}
\newcommand{\tf}{\widetilde{f}}
\newcommand{\g}{\mathfrak{g}}
\newcommand{\td}{\widetilde}

\newcommand{\mc}{\mathcal}
\newcommand{\mf}{\mathfrak}

\newcommand{\la}{\lambda}
\newcommand{\ep}{\epsilon}

\newcommand{\tl}[1]{\substack{\scalebox{0.75}{#1}}}
\newcommand{\Tl}[1]{\substack{\scalebox{0.95}{#1}}}

\newcommand{\lt}{\triangleleft}

\begin{document}

\title[Lusztig data of KN tableaux in type D]
{Lusztig data of Kashiwara-Nakashima tableaux in type D}

\author{IL-SEUNG JANG}

\address{Department of Mathematical Sciences, Seoul National University, Seoul 08826, Korea}
\email{is\_jang@snu.ac.kr}

\author{JAE-HOON KWON}

\address{Department of Mathematical Sciences and RIM, Seoul National University, Seoul 08826, Korea}
\email{jaehoonkw@snu.ac.kr}


\thanks{This work was supported by the National Research Foundation of Korea(NRF) grant funded by the Korea government(MSIT) (No. 2019R1A2C1084833).}

\begin{abstract}
We describe the embedding from the crystal of Kashiwara-Nakashima tableaux in type $D$ of an arbitrary shape into that of $\mathbf{i}$-Lusztig data associated to a family of reduced expressions $\mathbf{i}$ which are compatible with the maximal Levi subalgebra of type $A$. The embedding is described explicitly in terms of well-known combinatorics of type $A$ including the Sch\"{u}tzenberger's jeu de taquin and an analog of RSK algorithm.
\end{abstract}

\maketitle
\setcounter{tocdepth}{2}

\section{Introduction}

A {\em Kashiwara-Nakashima tableau} (KN tableau for short) is a combinatorial model for the crystal of a finite-dimensional irreducible representation of the classical Lie algebra, which is equal to a Young tableau in case of type $A$ \cite{KN}. 
Let $\bi$ be a reduced expression of the longest element $w_0$ in the Weyl group of a semisimple Lie algebra. 
An {\em $\bi$-Lusztig datum} is a parametrization of a PBW basis of the negative part of the associated quantized enveloping algebra, and hence parametrizes its crystal \cite{Lu90, Lu93}. Let us denote by ${\bf KN}_\la$ the crystal of KN tableaux of shape $\la$, and denote by ${\bf B}_{\bi}$ the set of $\bi$-Lusztig data.

Suppose that $\mf g$ is a classical Lie algebra. Let $\mf{l}$ be a proper maximal Levi subalgebra  of type $A$, and let $\mf u$ be the nilradical of the parabolic subalgebra $\mf p=\mf l + \mf b$, where $\mf b$ is a Borel subalgebra of $\mf g$. 
We take a reduced expression $\bi$ of $w_0$ such that
the positive roots of ${\mf u}$ precede those of ${\mf l}$ with respect to the corresponding convex ordering on the set of positive roots of $\mf{g}$.
Then an explicit combinatorial description of the crystal embedding 
\begin{equation}\label{eq:embedding KN to Bi}
\xymatrixcolsep{4pc}\xymatrixrowsep{3pc}\xymatrix{
{\bf KN}_\la \ \ar@{^{(}->}[r] &\ \B_{\bi}},
\end{equation}
(up to a shift of weight)
is given in \cite{K18-1} when $\mf g$ is of type $A$ and in \cite{K18-2} when $\mf g$ is of type $B$ and $C$. The embedding \eqref{eq:embedding KN to Bi} has several interesting features. For example, it can be described {\em only} in terms of the well-known combinatorics of type $A$ including the Sch\"{u}tzenberger's jeu de taquin sliding and RSK algorithm. The crystal structure $\B_{\bi}$ also has a very simple description in this case.

In this paper, we provide a combinatorial description of \eqref{eq:embedding KN to Bi} when $\mf g$ is of type $D$ (Theorem \ref{thm:main}). 
As in \cite{K18-1,K18-2}, the map is given by a composition of two embeddings, where one is the embedding of $\bf{KN}_\la$ into the crystal of a parabolic Verma module with respect to $\mf l$, say ${\bf V}_\la$, and the other is the embedding of ${\bf V}_\la$ into $\B_{\bi}$ (up to shift of weights). 
The embedding of ${\bf KN}_\la$ into ${\bf V}_\la$ is given by a non-trivial sequence of the Sch\"{u}tzenberger's jeu de taquin sliding called {\em separation}, which is the main result in this paper (Theorem \ref{thm:separation}). 
Here ${\bf KN}_\la$ is replaced without difficulty by another combinatorial realization ${\bf T}_\la$ of the crystal of integrable highest weight module, called {\em spinor model} \cite{K16}. 
The embedding of ${\bf V}_\la$ into $\B_{\bi}$ is obtained by applying an analogue of RSK due to Burge, which is a morphism of crystals \cite{JK19-1}, and the embedding \eqref{eq:embedding KN to Bi} in case of type $A$ \cite{K18-1}.

We remark that the algorithm of separation is a generalization of the one introduced in \cite{JK19-2}, which yields a new combinatorial formula for the branching rule associated to $(GL_n,O_n)$, and the embedding of ${\bf V}_\la$ into $\B_{\bi}$ also has an application to a combinatorial description of Kirillov-Reshetikhin crystals of type $D_n^{(1)}$ \cite{JK19-1}.

The paper is organized as follows. 
In Section \ref{sec:notations}, we introduce some necessary notations.
In Section \ref{sec:KN and Spinor}, we recall the notions of ${\bf KN}_\la$ and ${\bf T}_\la$, and describe the isomorphism between them.
Then we describe the embedding of ${\bf T}_\la$ into ${\bf V}_\la$ in Section \ref{sec:spinor into parabolic verma}, and embedding of ${\bf V}_\la$ into ${\bf B}_{\bi}$ in Section \ref{sec:parabolic into PBW} (up to shift of weights).

\section{Notations}\label{sec:notations}
\subsection{Semistandard tableaux}
Let $\Z_+$ denote the set of non-negative integers.
Let $\cP$ be the set of partitions or Young diagrams. 
We let $\cP_{n}=\{\,\la\in\cP\,|\,\ell(\la)\leq n\,\}$ for $n\geq 1$, 
where $\ell(\la)$ is the length of $\la$.
Let $\cP^{(1,1)}=\{\,\la \in \cP \,|\, \la' = (\la_i')_{i \ge 1}\,, \,\la_i'  \in 2\mathbb{Z}_+\,\}$, 
where $\la'$ is the conjugate of $\la$.
Let 
$\cP^{(1,1)}_n = \cP^{(1,1)} \cap \cP_n$.
Let $\lambda^\pi$ be the skew Young diagram obtained by $180^{\circ}$-rotation of $\lambda$.

Let $\N$ be the set of positive integers with the usual linear ordering and let $\ov{\mathbb{N}}$ be the set consisting of $\ov{i}$ $(i \in \mathbb{N})$ with the linear ordering $\ov{i} > \ov{j}$ for $i < j \in \mathbb{N}$.
For $n\in \N$,
we put $[n]=\{\,1,\dots,n\,\}$ and $[\ov{n}]=\{\,\ov{1},\dots,\ov{n}\,\}$, and assume that $[n] \cup [\ov{n}]$ has an ordering given by
\begin{equation*}
1 < 2 < \dots < n-1 <
\begin{array}{c}
n \\
\ov{n}
\end{array} < \ov{n-1} < \dots < \ov{2} < \ov{1}.
\end{equation*}

For a skew Young diagram $\lambda/\mu$, we denote by ${SST}_{\mc{A}}(\lambda/\mu)$ the set of semistandard tableaux of shape $\lambda/\mu$ with entries in a subset $\mc{A}$ of $\N \cup\ov{\N}$.
We put ${SST}(\lambda/\mu)={SST}_{\N}(\lambda/\mu)$ for short.
For $T\in {SST}_{\mc{A}}(\lambda/\mu)$,
let $w(T)$ be the word given by reading the entries of $T$ column by column from right to left and from top to bottom in each column, and let ${\rm sh}(T)$ denote the shape of $T$.
For $a\in [\ov{n}]$ and $U\in {SST}(\lambda^{\pi})$ with $\la\in\cP_n$, let $V\leftarrow a$ be the tableau obtained by applying the Schensted's column insertion of $a$ into $V$ in a reverse way starting from the rightmost column (cf. \cite{Ful})

\subsection{Tableaux with two columns}
We also use the following notations for description of our main result (cf.\,\cite{JK19-2}).

For $a,b,c\in \Z_+$, let $\lambda(a, b, c)$ be a skew Young diagram with at most two columns given by $(2^{b+c}, 1^a)/(1^b)$.
Let $T$ be a tableau of shape $\lambda(a, b, c)$. We denote the left and right columns of $T$ by $T^{\texttt{L}}$ and $T^{\texttt{R}}$, respectively. 
%

Let $T$ be a tableau. If necessary, we assume that it is placed on the plane with a horizontal line $L$, say $\mathbb{P}_L$, such that any box in $T$ is either below or above $L$, and at least one edge of a box in $T$ meets $L$.
We denote by $T^{\texttt{body}}$ and $T^{\texttt{tail}}$ the subtableaux of $T$ above and below $L$, respectively.
For example, 
\begin{equation*}
\begin{split}
T=\raisebox{5ex}{
\ytableausetup {mathmode, boxsize=0.9em} 
\begin{ytableau}
 \none & \none & \none & \none \\
 \none & \none & \tl{1} & \none \\
 \none & \none & \tl{2} & \none \\
 \none & \tl{1} & \tl{3} & \none \\
 \none[\!\!\!\!\mathrel{\raisebox{-0.5ex}{$\scalebox{0.45}{\dots\dots\dots\dots}$ }}] 
& \tl{2} & \tl{4} & \none[\quad\quad \mathrel{\raisebox{-0.5ex}{$\scalebox{0.45}{\dots\dots\dots\dots}$\ ${}_{\scalebox{0.75}{$L$}}$}}] \\
 \none & \tl{3} & \none & \none \\
 \none & \tl{5} & \none & \none \\
  \none & \none & \none & \none \\
\end{ytableau}}\quad\quad\quad
T^{\texttt{body}}=\raisebox{5ex}{
\ytableausetup {mathmode, boxsize=0.9em} 
\begin{ytableau}
 \none & \none & \none \\
 \none & \tl{1} & \none \\
 \none & \tl{2} & \none \\
 \tl{1} & \tl{3} & \none \\
 \tl{2} & \tl{4} & \none \\
\end{ytableau}}\quad
T^{\texttt{tail}}=\raisebox{5ex}{
\ytableausetup {mathmode, boxsize=0.9em} 
\begin{ytableau}
 \none & \none & \none \\
 \none & \none & \none \\
 \none & \none & \none \\
 \none & \none & \none \\
 \none & \none & \none \\
 \tl{3} & \none & \none \\
 \tl{5} & \none & \none \\
\end{ytableau}}\quad\quad
\end{split}
\end{equation*}
where the dotted line denotes $L$.

For a tableau $U$ with the shape of a single column, let {$\textrm{ht}(U)$} denote the height of $U$ and we put $U(i)$ (resp. $U[i]$) to be $i$-th entry of $U$ from bottom (resp. top). 
We also write
\begin{equation*} 
	U=\left(U(\ell),\dots , U(1)\right)=\left(U[1], \dots, U[\ell] \right),
\end{equation*} 
where $\ell={\rm ht}(U)$.
Suppose that $U$ is a tableau in $\mathbb{P}_L$. To emphasize gluing and cutting tableaux with respect to $L$,
we also write
\begin{equation*} 
\begin{split}
& U^{\texttt{body}}\boxplus U^{\texttt{tail}} = U,\quad U \boxminus U^{\texttt{tail}} = U^{\texttt{body}}.
\end{split}
\end{equation*}

For a sequence of tableaux $U_1,U_2,\dots,U_m$ in $\mathbb{P}_L$, whose shapes are single columns, let us say that $(U_1,U_2,\dots,U_m)$ is {\em semistandard along $L$} if they form a semistandard tableau $T$ of a skew shape with $U_i$ the $i$-th column of $T$ from the left.

\section{Kashiwara-Nakashima tableaux and spinor tableaux}\label{sec:KN and Spinor}

In this section, we recall two combinatorial models for the crystals of type $D_n$, say KN tableaux model and spinor model, and then give an explicit combinatorial description of the isomorphism between them.

\subsection{Lie algebra of type $D$} \label{subsec:basics}
We assume that ${\mf g}$ is the simple Lie algebra of type $D_n$ $(n \ge 4)$.
The weight lattice is $P=\bigoplus_{i=1}^n\Z\epsilon_i$ with a symmetric bilinear form $(\, \, ,\,)$ such that $(\epsilon_i|\epsilon_j)=\delta_{ij}$ for $i,j$.
Put $I=\{\,1,\dots,n\,\}$.
The set of simple roots is $\{\,\alpha_i\,|\,i\in I\,\}$, 
where $\alpha_i=\ep_i-\ep_{i+1}$ for $1\leq i\leq n-1$, and $\alpha_n=\ep_{n-1}+\ep_n$,
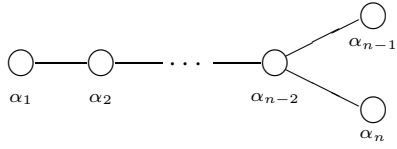
\begin{center} 
\setlength{\unitlength}{0.19in}
\begin{picture}(15,4.5)
\put(3.4,2){\makebox(0,0)[c]{$\bigcirc$}}
\put(5.6,2){\makebox(0,0)[c]{$\bigcirc$}}
\put(10.4,2){\makebox(0,0)[c]{$\bigcirc$}}
\put(13.1,3.3){\makebox(0,0)[c]{$\bigcirc$}}
\put(13.1,0.7){\makebox(0,0)[c]{$\bigcirc$}}
\put(3.8,2){\line(1,0){1.4}}
\put(6,2){\line(1,0){1.3}} 
\put(8.7,2){\line(1,0){1.3}} 
%
\put(10.7,2.2){\line(2,1){2}}
\put(10.7,1.8){\line(2,-1){2}}

\put(8,1.95){\makebox(0,0)[c]{$\cdots$}}
\put(3.4,1){\makebox(0,0)[c]{\tiny ${\alpha}_1$}}
\put(5.6,1){\makebox(0,0)[c]{\tiny ${\alpha}_2$}}
\put(10.4,1){\makebox(0,0)[c]{\tiny ${\alpha}_{n-2}$}}
\put(13.1,2.5){\makebox(0,0)[c]{\tiny ${\alpha}_{n-1}$}}
\put(13.1,0.0){\makebox(0,0)[c]{\tiny ${\alpha}_{n}$}}

\end{picture}
\end{center} 
\vskip 3mm
and the set of positive roots is $\Phi^+=\{\,\ep_i{\pm}\ep_j\,|\,1\leq i<j\leq n\,\}$. 
The fundamental weights $\Lambda_i$ ($i\in I$) are given by $\Lambda_i=\sum_{k=1}^i\epsilon_k$ for $i=1,\ldots, n-2$, $\Lambda_{n-1}=(\epsilon_1+\cdots+\epsilon_{n-1}-\epsilon_n)/2$ and $\Lambda_{n}=(\epsilon_1+\cdots+\epsilon_{n-1}+\epsilon_n)/2$. 

Let 
\begin{equation*}
\begin{split}
\mc{P}_n & = \left\{\, (\lambda_1, \dots, \lambda_n) \, \Big\vert \, \lambda_i \in \tfrac{1}{2}\mathbb{Z}, \, \lambda_i - \lambda_{i+1} \in \mathbb{Z}_+, \, \lambda_{n-1} \ge |\lambda_n| \, \right\}.\\
\end{split}
\end{equation*}
For $\la\in \mc{P}_n$, we put 
$$\omega_{\lambda} = \sum_{i=1}^n \lambda_i \epsilon_i\,.$$
Then $P^+ = \{ \omega_{\lambda} \, | \, \lambda \in \mc{P}_n \}$ is the set of dominant integral weights.
We put ${\texttt{sp}^+}={\left( \left( \frac{1}{2} \right)^n \right)}$ and ${\texttt{sp}^-}={( ( \frac{1}{2} )^{n-1} , -\frac{1}{2} )})$ for simplicity. We also identify $\lambda \in \mc P_n$ with a (generalized) Young diagram, which may have a half-width box on the leftmost column \cite[Section 6.7]{KN},

Put $J = I \setminus \{ n \}$.
Let ${\mf l}$ be the Levi subalgebra of ${\mf g}$ associated to $\{\, \alpha_i \, | \, i \in J \,\}$, which is of type $A_{n-1}$. 
Let $W$ be the Weyl group generated by the simple reflections $s_i$ ($i\in I$) with the longest element $w_0$, and 
let $R(w_0)$ be the set of reduced expressions of $w_0$.
Recall that $W$ acts faithfully on $P$ by $s_i(\ep_i)=\ep_{i+1}$, $s_i(\ep_k)=\ep_k$ for $1\leq i\leq n-1$ and $k\neq i, i+1$, and $s_n(\ep_{n-1})=-\ep_n$ and $s_n(\ep_k)=\ep_k$ for $k\neq n-1, n$.


Let $U_q({\mf g})$ be the quantized universal enveloping algebra associated to $\g$.
For $\Lambda\in P^+$, we denote by $B(\Lambda)$ the crystal associated to an irreducible highest weight $U_q(\g)$-module with highest weight $\Lambda$. For $\la\in P$, let $T_\la=\{\,t_\la\,\}$ be the crystal, where ${\rm wt}(t)=\la$, $\te_i t_\la=\tf_i t_\la=\bf{0}$, and $\varepsilon_i(t_\la)=\varphi_i(t_\la)=-\infty$ for $i\in I$.
%
We refer the reader to \cite{Kas91, Kas95} for more details of crystals.

\subsection{Kashiwara-Nakashima tableaux} 
Suppose that $\la = (\lambda_1, \dots, \lambda_n)  \in \mc{P}_n$ is given.
The notion of Kashiwara-Nakashima tableaux (KN tableaux, for short) of type $D$ \cite{KN} is a combinatorial model of $B(\omega_\la)$.
In this paper, we need an analogue, which is obtained from the one in \cite{KN} by applying $180^\circ$ rotation and replacing $i$ and $\ov{i}$ (resp. $\ov{i}$ with $i$). For the reader's convenience, let us give its definition and crystal structure.

\begin{df} \label{df:d-KN tableau}
{\rm 
For $\lambda = (\lambda_1, \dots, \lambda_n) \in \mc{P}_n$, let $T$ be a tableau of shape $\la^{\pi}$ with entries in $[n] \cup [\ov{n}]$ such that
\begin{itemize}
	\item[(1)] $T(i,j)\not\geq T(i+1,j)$ and $T(i,j)\leq T(i,j+1)$ for each $i$ and $j$,
	\item[(2)] $n$ and $\ov{n}$ can appear successively in $T$ other than half-width boxes,
	\item[(3)] $i$ and $\ov{i}$ do not appear simultaneously in the half-width boxes,
\end{itemize}
where $T(i, \, j)$ denotes the entry in $T$ located in the $i$-th row from the bottom and the $j$-th column from the right. Then $T$ is called a {\em KN tableau of type $D_n$} if it satisfies the following conditions:
\begin{itemize}
	\item[(${\mf d}$-1)] If $T(p, j) = \ov{i}$ and $T(q, j) = i$ for some $i \in [n]$ with $p < q$, then $(q-p)+i > \lambda_j'$.
	\item[(${\mf d}$-2)] Suppose $\lambda_n \ge 0$ and $\lambda_j' = n$. If $T(k, j) = n$ (resp. $\ov{n}$), then $k$ is odd (resp. even).
	\item[(${\mf d}$-3)] Suppose $\lambda_n < 0$ and $\lambda_j' = n$. If $T(k, j) = n$ (resp. $\ov{n}$), then $k$ is even (resp. odd).
	\item[(${\mf d}$-4)] If either $T(p, j)=\ov{a}$, $T(q, j)=\ov{b}$, $T(r, j)=b$ and $T(s, j+1)=a$ or $T(p, j) = \ov{a}$, $T(q, j+1)=\ov{b}$, $T(r, j+1)=b$ and $T(s, j+1) = a$ with $p \le q < r \le s$ and $a \le b < n$, then $(q-p)+(s-r) < b-a$.
	\item[(${\mf d}$-5)] Suppose $T(p, j) = \ov{a}$, $T(s, j+1) = a$ with $p < s$. If there exists $p \le q < s$ such that either $T(q, j),\, T(q+1, j) \in \{ \, n, \ov{n} \, \}$ with $T(q, j) \neq T(q+1, j)$ or $T(q, j+1),\, T(q+1, j+1) \in \{ \, n, \ov{n} \, \}$ with $T(q, j+1) \neq T(q+1, j+1)$, then $s-p \le n-a$.
	\item[(${\mf d}$-6)] It is not possible that $T(p, j) \in \{ \, n, \ov{n} \, \}$ and $T(s, j+1) \in \{ \, n, \ov{n} \, \}$ with $p < s$.
	\item[(${\mf d}$-7)] Suppose $T(p, j) = \ov{a}$, $T(s, j+1) = a$ with $p < s$. If $T(q, j+1) \in \{ n, \ov{n} \}$, $T(r, j) \in \{ n, \ov{n} \}$ and $s-q+1$ is either odd or even with $p \le q < r \le s$ and $a < n$, then $s-p < n-a$.
\end{itemize}	
  
We denote by ${\bf KN}_{\lambda}$ the set of KN tableaux of shape $\la^\pi$.
}
\end{df}


%
Recall that ${\bf KN}_{(1)}$ has the following crystal structure isomorphic to that of $B(\Lambda_1)$ 
\vskip 1mm
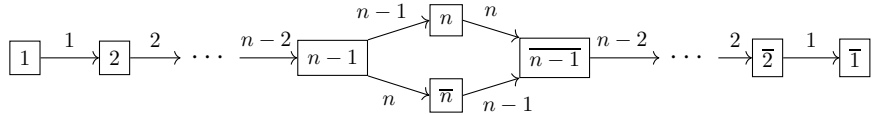
\begin{equation*}
\begin{tikzpicture}[
squarednode/.style={rectangle, draw=black, fill=white, minimum size=0.95mm},]
	\node[squarednode](1) {$\tl{1}$};
	\node[squarednode](2) [right=28] {$\tl{2}$};
	\node(3) [right=58.9] {$\cdots$};
	\node[squarednode](4) [right=103] {$\tl{$n-1$}$};
	\node[squarednode](5) at (5.6,0.5) {$\tl{$n$}$};
	\node[squarednode](6) at (5.6,-0.5) {$\tl{$\ov{n}$}$};
	\node[squarednode](7) [right=187] {$\tl{$\ov{n-1}$}$};
	\node(8) [right=240] {$\cdots$};
	\node[squarednode](9) [right=275] {$\tl{$\ov{2}$}$};
	\node[squarednode](10) [right=308] {$\tl{$\ov{1}$}$};
	
	\draw[->] (1.east) -- node[above]{$\tl{1}$} (2.west);
	\draw[->] (2.east) -- node[above]{$\tl{2}$} (3.west);
	\draw[->] (3.east) -- node[above]{$\tl{$n-2$}$\,} (4.west);
	\draw[->] (4.north east) -- node[above]{$\tl{$n-1$}$\, \,\,\,\,} (5.west);
	\draw[->] (4.south east) -- node[below]{$\tl{$n$}$\, \,} (6.west);
	\draw[->] (5.east) -- node[above]{$\tl{$n$}$} (7.north west);
	\draw[->] (6.east) -- node[below]{\,\, \, $\tl{$n-1$}$} (7.south west);
	\draw[->] (7.east) -- node[above]{\, $\tl{$n-2$}$\, \,} (8.west);
	\draw[->] (8.east) -- node[above]{$\tl{2}$} (9.west);
	\draw[->] (9.east) -- node[above]{$\tl{1}$} (10.west);
\end{tikzpicture}
\end{equation*}
where 
$\raisebox{-.6ex}{{\tiny ${\def\lr#1{\multicolumn{1}{|@{\hspace{.6ex}}c@{\hspace{.6ex}}|}{\raisebox{0ex}{$#1$}}}\raisebox{0ex}
{$\begin{array}[b]{c}
\cline{1-1}
\lr{ a }\\
\cline{1-1}
\end{array}$}}$}}
\overset{i}{\longrightarrow}
\raisebox{-.6ex}{{\tiny ${\def\lr#1{\multicolumn{1}{|@{\hspace{.6ex}}c@{\hspace{.6ex}}|}{\raisebox{0ex}{$#1$}}}\raisebox{0ex}
{$\begin{array}[b]{c}
\cline{1-1}
\lr{ b }\\
\cline{1-1}
\end{array}$}}$}}$ 
means 
$\widetilde{f}_i
\raisebox{-.6ex}{{\tiny ${\def\lr#1{\multicolumn{1}{|@{\hspace{.6ex}}c@{\hspace{.6ex}}|}{\raisebox{0ex}{$#1$}}}\raisebox{0ex}
{$\begin{array}[b]{c}
\cline{1-1}
\lr{ a }\\
\cline{1-1}
\end{array}$}}$}} =
\raisebox{-.6ex}{{\tiny ${\def\lr#1{\multicolumn{1}{|@{\hspace{.6ex}}c@{\hspace{.6ex}}|}{\raisebox{0ex}{$#1$}}}\raisebox{0ex}
{$\begin{array}[b]{c}
\cline{1-1}
\lr{ b }\\
\cline{1-1}
\end{array}$}}$}}
$ with $\tf_i$ the Kashiwara operator for $i\in I$, and
$
{\rm wt}
(\,\raisebox{-.6ex}{{\tiny ${\def\lr#1{\multicolumn{1}{|@{\hspace{.9ex}}c@{\hspace{.9ex}}|}{\raisebox{0ex}{$#1$}}}\raisebox{0.2ex}
{$\begin{array}[b]{c}
\cline{1-1}
\lr{ i }\\
\cline{1-1}
\end{array}$}}$}}\,) = \epsilon_i$, 
${\rm wt}
(\,\raisebox{-.6ex}{{\tiny ${\def\lr#1{\multicolumn{1}{|@{\hspace{.9ex}}c@{\hspace{.9ex}}|}{\raisebox{-0.2ex}{$#1$}}}\raisebox{0.2ex}
{$\begin{array}[b]{c}
\cline{1-1}
\lr{ \ov{i} }\\
\cline{1-1}
\end{array}$}}$}}\,) = -\epsilon_i.
$ 
for $i = 1, \dots, n$.
On the other hand, ${\bf KN}_{\texttt{sp}+}$ and
${\bf KN}_{\texttt{sp}-}$ have crystal structures isomorphic to those of $B(\Lambda_n)$ and $B(\Lambda_{n-1})$ which are the crystals of spin representations with highest weights $\Lambda_n$ and $\Lambda_{n-1}$, respectively.
For $i\in I$, $\widetilde{f}_i$ on ${\bf KN}_{\texttt{sp}\pm}$ is given by
\begin{equation} \label{eq:Kashiwara operator on spin crystal}
\begin{split}
\raisebox{-.6ex}{{\tiny ${\def\lr#1{\multicolumn{1}{|@{\hspace{.6ex}}c@{\hspace{.6ex}}|}{\raisebox{-.3ex}{$#1$}}}\raisebox{-.6ex}
{$\begin{array}[b]{c}
\lr{\tl{$\vdots$\,}\!}\\
\cline{1-1}
\lr{ i }\\
\cline{1-1}
\lr{\tl{$\vdots$\,}\!}\\
\cline{1-1}
\lr{ \!\ov{i+1}\!}\\
\cline{1-1}
\lr{\tl{$\vdots$\,}\!}\\
\end{array}$}}$}}
\end{split} \, \overset{\widetilde{f}_i\,\, (i \neq n)}{\xrightarrow{\hspace*{1.2cm}}} \,
\begin{split}
\raisebox{-.0ex}{{\tiny ${\def\lr#1{\multicolumn{1}{|@{\hspace{.6ex}}c@{\hspace{.6ex}}|}{\raisebox{-.3ex}{$#1$}}}\raisebox{-.6ex}
{$\begin{array}[b]{c}
\lr{\tl{$\vdots$\,}\!}\\
\cline{1-1}
\lr{ i+1 }\\
\cline{1-1}
\lr{\tl{$\vdots$\,}\!}\\
\cline{1-1}
\lr{ \!\ov{i}\!}\\
\cline{1-1}
\lr{\tl{$\vdots$\,}\!}\\
\end{array}$}}$}}
\end{split} \quad \quad \quad \quad \quad
\begin{split}
\raisebox{-.0ex}{{\tiny ${\def\lr#1{\multicolumn{1}{|@{\hspace{.6ex}}c@{\hspace{.6ex}}|}{\raisebox{-.3ex}{$#1$}}}\raisebox{-.6ex}
{$\begin{array}[b]{c}
\lr{\tl{$\vdots$\,}\!}\\
\cline{1-1}
\lr{ n-1 }\\
\cline{1-1}
\lr{ \! n\! }\\
\cline{1-1}
\lr{\tl{$\vdots$\,}\!}\\
\end{array}$}}$}}
\end{split} \, \overset{\widetilde{f}_n}{\xrightarrow{\hspace*{1cm}}} \,
\begin{split}
\raisebox{-.6ex}{{\tiny ${\def\lr#1{\multicolumn{1}{|@{\hspace{.6ex}}c@{\hspace{.6ex}}|}{\raisebox{-.3ex}{$#1$}}}\raisebox{-.6ex}
{$\begin{array}[b]{c}
\lr{\tl{$\vdots$\,}\!}\\
\cline{1-1}
\lr{ \!\ov{n}\! }\\
\cline{1-1}
\lr{ \!\ov{n-1}\!}\\
\cline{1-1}
\lr{\tl{$\vdots$\,}\!}\\
\end{array}$}}$}}
\end{split}\,.
\end{equation}
\vskip 1.5mm

%
Let $\lambda \in \mc{P}_n$ be given. Let us identify $T \in {\bf KN}_{\lambda}$ with its word $w(T)$ so that we may regard
\begin{equation*}
	{\bf KN}_{\lambda} \subset
	\left\{
\begin{array}{ll}
	\left( {\bf KN}_{(1)} \right)^{\otimes \rm N}, & \,\,\textrm{if $\la_n \in \Z$}, \\
	{\bf KN}_{\texttt{sp}^{\pm}} \otimes \left( {\bf KN}_{(1)} \right)^{\otimes \rm N},   & \,\,\textrm{if $\la_n\not\in \Z$}, 
\end{array}
\right.
\end{equation*}
where ${\rm N}$ is the number of letters in $w(T)$ except for the one in half-width boxes.
Then ${\bf KN}_{\lambda}$ is invariant under $\widetilde{e}_i$ and $\widetilde{f}_i$ for $i\in I$, and 
\begin{equation*} \label{eq:KN regular}
	{\bf KN}_{\lambda} \cong B(\omega_{\lambda}),
\end{equation*}
(\cite[Theorem 6.7.1]{KN}).




\subsection{Spinor model} \label{subsec:spinor}
Let us briefly recall another combinatorial model for $\B(\omega_\la)$ (see \cite{JK19-2,K16} for more details and examples). We keep the notations used in \cite{JK19-2}.

For $T \in {SST}_{[\ov{n}]}(\lambda(a, b, c))$ and $0 \le k \le \min\{ a, b \}$, we slide down $T^{\texttt{R}}$ by $k$ positions to have a tableau $T'$ of shape $\lambda(a-k, b-k, c+k)$. We define $\mathfrak{r}_T$ to to be the maximal $k$ such that $T'$ is semistandard.

For $T\in {SST}_{[\ov{n}]}(\la(a,b,c))$ with ${\mf r}_T=0$, we define $\mc E T$ and $\mc F T$ as follows:
\begin{itemize}
\item[(1)] $\mc E T$ is tableau in ${SST}_{[\ov{n}]}(\la(a-1,b+1,c))$ obtained from $T$ by applying Sch\"{u}tenberger's jeu de taquin sliding to the position below the bottom of $T^{\texttt{R}}$, when $a>0$,

\item[(2)] $\mc F T$ is tableau in ${SST}_{[\ov{n}]}(\la(a+1,b-1,c))$ obtained from $T$ by applying jeu de taquin sliding to the position above the top of $T^{\texttt{L}}$, when $b>0$.

\end{itemize}
Here we assume that $\mathcal{E}T = \mathbf{0}$ and $\mathcal{F}T = \mathbf{0}$ when $a=0$ and $b=0$, respectively, where ${\bf 0}$ is a formal symbol. In general, if $\mathfrak{r}_T=k$, then we define $\mc{E}T = \mc{E}T'$ and $\mc{F}T = \mc{F}T'$, where $T'$ is obtained from $T$ by sliding down $T^{\texttt{R}}$ by $k$ positions and hence $\mf{r}_{T'}=0$.

Let
\begin{equation} \label{eq:def_T}
\begin{split}
	\mathbf{T}(a)  & = \left\{\, T \,|\, T \in {SST}_{[\ov{n}]}(\lambda(a,b,c)),\ b,c \in 2\Z_+,\ \mathfrak{r}_T \le 1 \,\right\}\quad (0\leq a\leq n-1), \\
	\overline{\mathbf{T}}(0) &   = \bigsqcup_{b, c\, \in 2\mathbb{Z}_+} {SST}_{[\ov{n}]}(\lambda(0, b, c+1)),  \quad \ \ \mathbf{T}^{\textrm{sp}}  = \bigsqcup_{a \in \mathbb{Z}_+}{SST}_{[\ov{n}]}((1^a)), \\
	\mathbf{T}^{\textrm{sp}+} &  = \{\, T \, |\,  T \in \mathbf{T}^{\textrm{sp}}, \, \mathfrak{r}_T = 0\, \}, \ \ \mathbf{T}^{\textrm{sp}-} = \{ \,T \, |\, T \in \mathbf{T}^{\textrm{sp}}, \, \mathfrak{r}_T = 1\, \},
\end{split}
\end{equation}
where $\mathfrak{r}_T$ of $T \in \mathbf{T}^{\textrm{sp}}$ is defined to be the residue of $\textrm{ht}(T)$ modulo $2$.
%
%
For $T \in \mathbf{T}(a)$, we define $(T^{{\texttt{L}}*}, T^{{\texttt{R}}*})$ when $\mathfrak{r}_T = 1$, and $({}^{\texttt{L}} T, {}^{\texttt{R}} T)$ by
\begin{equation*} 
\begin{split}
	(T^{{\texttt{L}}*}, T^{{\texttt{R}}*})  & = \ ((\mc{F}T)^{\texttt{L}},(\mc{F}T)^{\texttt{R}}), \\
 	({}^{\texttt{L}} T, {}^{\texttt{R}} T) &  = ((\mc{E}^{a^\ast} T)^{\texttt{L}},(\mc{E}^{a^\ast} T)^{\texttt{R}}) \quad(a^\ast=a-\mf{r}_T)
\end{split}	
\end{equation*}

\begin{df}\label{def:admissibility}
{\rm
Let $a, a'$ be given with $0\leq a'\leq a \leq n-1$. We say a pair $(T,S)$ is {\em admissible}, and write $T \prec S$ if it is one of the following cases:
\begin{itemize}
	\item[(1)] $(T,S) \in \mathbf{T}(a) \times \mathbf{T}(a')$ or $\mathbf{T}(a) \times \mathbf{T}^{\textrm{sp}}$ with 
		\begin{equation*}
		\begin{split}
			& (\text{i}) \ \ \ \  \textrm{ht}(T^{\texttt{R}}) \le \textrm{ht}(S^{\texttt{L}}) - a' + 2\mathfrak{r}_T \mathfrak{r}_S, \\
			& (\text{ii}) \ \ \left\{ \begin{array}{ll} 
				T^{\texttt{R}}(i) \le {}^{\texttt{L}} S(i), & \textrm{if} \ \mf{r}_T \mf{r}_S  = 0, \\
				T^{{\texttt{R}}*}(i) \le {}^{\texttt{L}} S(i), & \textrm{if} \ \mathfrak{r}_T   \mathfrak{r}_S = 1, 
			\end{array} \right. \\
			& (\text{iii}) \ \ \left\{ \begin{array}{ll} 
				{}^{\texttt{R}} T(i+a-a') \le S^{\texttt{L}}(i), & \textrm{if} \ \mathfrak{r}_T  \mathfrak{r}_S = 0, \\
				{}^{\texttt{R}} T(i+a-a'+\varepsilon) \le S^{{\texttt{L}}*}(i), & \textrm{if} \ \mathfrak{r}_T   \mathfrak{r}_S = 1, 
			\end{array} \right. \\
		\end{split}
		\end{equation*}
		for $i \ge 1$. Here $\varepsilon = 1$ if $S \in \mathbf{T}^{\textrm{sp}-}$ and $0$ otherwise, and we assume that $a' = \mathfrak{r}_S$, $S = S^{\texttt{L}} = {}^{\texttt{L}} S = S^{{\texttt{L}}*}$ when $S \in \mathbf{T}^{\textrm{sp}}$.
		
	\item[(2)] $(T,S) \in \mathbf{T}(a) \times \overline{\mathbf{T}}(0)$ with $T \prec S^{\texttt{L}}$ in the sense of (1), regarding $S^{\texttt{L}} \in \mathbf{T}^{\textrm{sp}-}$.
	
	\item[(3)] $(T,S) \in \overline{\mathbf{T}}(0) \times \overline{\mathbf{T}}(0)$ or $\overline{\mathbf{T}}(0) \times \mathbf{T}^{\textrm{sp}-}$ with $(T^{\texttt{R}}, S^{\texttt{L}}) \in \overline{\mathbf{T}}(0)$.
\end{itemize}}
\end{df}



\begin{rem}\label{rem:admissibilty for spin column}
{\rm 
(1) For $T\in {\bf T}(a)$, we assume that $T\in \mathbb{P}_L$ such that the subtableau of single column with height $a$ is below $L$ and hence equal to $T^{\texttt{tail}}$.
\begin{equation*}
\begin{split}
T=\ \raisebox{6ex}{
\ytableausetup {mathmode, boxsize=0.9em} 
\begin{ytableau}
 \none & \none & \none & \none \\
 \none & \none & \tl{1} & \none \\
 \none & \none & \tl{2} & \none \\
 \none & \tl{1} & \tl{3} & \none \\
 \none[\ \mathrel{\!\!\!\!\raisebox{-0.5ex}{$\scalebox{0.45}{\dots\dots\dots\dots}$ }}] 
& \tl{2} & \tl{4} & \none[\quad\quad \mathrel{\raisebox{-0.5ex}{$\scalebox{0.45}{\dots\dots\dots\dots}$\ ${}_{\scalebox{0.7}{$L$}}$}}] \\
 \none & \tl{3} & \none & \none \\
 \none & \tl{5} & \none & \none \\
\end{ytableau}}\quad\quad  \in {\bf T}(2)
\end{split}
\end{equation*}

(2)
Let $S \in \mathbf{T}^{\textrm{sp}}$ with $\varepsilon=\mf{r}_S$. 
We may assume that $S= U^{\texttt{L}}$ for some $U\in {\bf T}(\varepsilon)$, where $U^{\texttt{R}}(i)$ ($i\ge 1$) are sufficiently large
so that $S = U^{\texttt{L}} = {}^{\texttt{L}} U$. Then we may understand the condition Definition \ref{def:admissibility}(1) for $(T,S) \in \mathbf{T}(a) \times \mathbf{T}^{\textrm{sp}}$ as induced from the one for $(T,U) \in \mathbf{T}(a) \times \mathbf{T}(\varepsilon)$.

}
\end{rem}
\vskip 2mm

Let $\B$ be one of $\mathbf{T}(a)$ $(0\leq a \leq n-1)$, $\mathbf{T}^{\textrm{sp}}$, and $\overline{\mathbf{T}}(0)$. The $\mf g$-crystal structure on $\B$ \cite{K16} is given as follows.
Let $T \in \B$ given. 
For $i\in J$, we define $\te_i, \tf_i$ by
regarding $\B$ as an $\mf{l}$-subcrystal of $\bigsqcup_{\lambda \in \cP_n} {SST}_{[\ov{n}]}(\lambda)$ \cite{KN}, where we consider the set $[\ov{n}]$ as the dual crystal of $[n]$ that is the crystal of vector representation of ${\mf l}$.
%
For $i=n$ and $T\in \B$, we define $\te_nT$ and $\tf_nT$ as follows:

\begin{itemize}
	\item[(1)] if $\B = \mathbf{T}^{\textrm{sp}}$, then $\te_nT$ is the tableau obtained by removing a domino 
\raisebox{-.6ex}{{\tiny ${\def\lr#1{\multicolumn{1}{|@{\hspace{.6ex}}c@{\hspace{.6ex}}|}{\raisebox{-.3ex}{$#1$}}}\raisebox{-.6ex}
{$\begin{array}[b]{c}
\cline{1-1}
\lr{ \ov{n} }\\
\cline{1-1}
\lr{ \!\ov{n-1}\!}\\
\cline{1-1}
\end{array}$}}$}}
from $T$ if it is possible, and \textbf{0} otherwise, and $\tf_n T$ is given in a similar way by adding
\raisebox{-.6ex}{{\tiny ${\def\lr#1{\multicolumn{1}{|@{\hspace{.6ex}}c@{\hspace{.6ex}}|}{\raisebox{-.3ex}{$#1$}}}\raisebox{-.6ex}
{$\begin{array}[b]{c}
\cline{1-1}
\lr{ \ov{n} }\\
\cline{1-1}
\lr{ \!\ov{n-1}\!}\\
\cline{1-1}
\end{array}$}}$}},

	\item[(2)] if $\B = \mathbf{T}(a)$ or $\overline{\mathbf{T}}(0)$, then
	$\te_n T=\te_n \left( T^{\texttt{R}} \otimes T^{\texttt{L}} \right)$ and $\tf_n T=\tf_n \left( T^{\texttt{R}} \otimes T^{\texttt{L}} \right)$ regarding $\B \subset \left( \mathbf{T}^{\textrm{sp}} \right)^{\otimes 2}$.
\end{itemize}
The weight of $T\in \B$ is given by
\begin{equation*}
\textrm{wt}(T) =
\begin{cases} 
2\omega_n + \sum_{i\geq 1}m_i\epsilon_i , & \textrm{if $T \in \mathbf{T}(a)$ or $\ov{\bf T}(0)$}, \\
\omega_n + \sum_{i\geq 1}m_i\epsilon_i, & \textrm{if $T\in {\bf T}^{\rm sp}$}.
\end{cases},\\
\end{equation*}
where $m_i$ is the number of occurrences of $\ov{i}$ in $T$.
Then $\B$ is a regular $\mathfrak{g}$-crystal with respect to $\te_i$ and $\tf_i$ for $i \in I$, and   
\begin{equation*} 
\begin{split}
	& \mathbf{T}(a) \cong \B(\omega_{n-a}) \ \  (2 \le a \le n-1), \\ 
	& \mathbf{T}(0) \cong \B(2\omega_{n}), \quad \overline{\mathbf{T}}(0) \cong \B(2\omega_{n-1}), \quad \mathbf{T}(1) \cong \B(\omega_{n-1}+\omega_n), \\
	& \mathbf{T}^{\textrm{sp}-} \cong \B(\omega_{n-1}), \quad \mathbf{T}^{\textrm{sp}+} \cong \B(\omega_n).
\end{split}
\end{equation*}
(\cite[Proposition 4.2]{K16}).
Note that the highest weight element $H$ of ${\bf B}$ is of the following form:
\begin{equation} \label{eq:hwvs in spinor}
	H=\left\{
		\begin{array}{ll}
			\emptyset \, \boxplus \, H_{(1^a)}, & \textrm{if ${\bf B} = {\bf T}(a)$ with $2 \le a \le n-1$,} \\
			\emptyset \, \boxplus \, \raisebox{-.0ex}{{\tiny ${\def\lr#1{\multicolumn{1}{|@{\hspace{.6ex}}c@{\hspace{.6ex}}|}{\raisebox{-.3ex}{$#1$}}}\raisebox{-.6ex}
{$\begin{array}[b]{c}
\cline{1-1}
\lr{ \!\ov{n}\!}\\
\cline{1-1}
\end{array}$}}$}}\,, & \textrm{if ${\bf B}={\bf T}(1)$,} \\
			\emptyset \,, & \textrm{if ${\bf B} = {\bf T}(0)$ or ${\bf B} = {\bf T}^{\rm sp+}$,} \\
			\, \raisebox{-.0ex}{{\tiny ${\def\lr#1{\multicolumn{1}{|@{\hspace{.6ex}}c@{\hspace{.6ex}}|}{\raisebox{-.3ex}{$#1$}}}\raisebox{-.6ex}
{$\begin{array}[b]{c}
\cline{1-1}
\lr{ \!\ov{n}\!}\\
\cline{1-1}
\end{array}$}}$}}\,, & \textrm{if ${\bf B} = {\bf T}^{\rm sp-}$,} 		
\end{array}
	\right.
\end{equation}
where $\emptyset$ is the empty tableau and $H_{(1)}\in SST_{[n]}((1^a))$ $(2 \le a \le n-1)$ such that $H_{(1^a)}[k]=\ov{n-k+1}$ ($1\leq k\leq a$), that is,
\begin{equation}\label{eq:H_(1^a)}
\hspace{5cm} H_{(1^a)} = \hspace{-6.2cm}
\begin{split}
\raisebox{-.6ex}{{\tiny ${\def\lr#1{\multicolumn{1}{|@{\hspace{.6ex}}c@{\hspace{.6ex}}|}{\raisebox{-.3ex}{$#1$}}}\raisebox{-.6ex}
{$\begin{array}[b]{c}
\cline{1-1}
\lr{ \!\ov{n}\! }\\
\cline{1-1}
\lr{ \!\ov{n-1}\!}\\
\cline{1-1}
\lr{ \!\vdots\!}\\
\cline{1-1}
\lr{ \!\ov{n-a+1}\!}\\
\cline{1-1}
\end{array}$}}$}}
\end{split}\,.
\end{equation}
Note that the empty tableau $\emptyset$ is an element of ${SST}_{[\ov{n}]}((0))$.

Let $\lambda = (\lambda_1, \dots, \lambda_n) \in \mc{P}_n$ be given. 
We have 
\begin{equation} \label{eq:weight expression}
	\omega_{\lambda} = 
	\left\{
\begin{array}{ll}
	\sum_{i=1}^{\ell}\omega_{n-a_i}+p\omega_n
	+q(2\omega_n)+r\omega_n , & \textrm{if $\lambda_n \ge 0$}, \\
	\sum_{i=1}^{\ell}\omega_{n-a_i}+p\omega_n
	+\ov{q}(2\omega_{n-1})+\ov{r}\omega_{n-1},	& \textrm{if $\lambda_n < 0$},
\end{array}
\right.
\end{equation}
where 
$a_{\ell} \ge \dots \ge a_1 \ge 1$, $p$ is the number of $i$ such that $a_i = 1$ and
($q$, $r$) (resp. ($\ov{q}$, $\ov{r}$)) is given by $2\lambda_n = 2q+r$ with $r \in \{ \, 0, 1 \, \}$ (resp. $-2\lambda_n = 2\ov{q}+\ov{r}$ with  $\ov{r} \in \{ \, 0, 1 \, \}$).

Let
\begin{equation}\label{eq:hat{T}lambda}
	\widehat{\bf T}_{\lambda} =
	\left\{
\begin{array}{ll}
	{\bf T}(a_{\ell}) \times \dots \times {\bf T}(a_1) \times {\bf T}(0)^{\times q} \times ( {\bf T}^{\rm sp+} )^{r} , & \textrm{if $\lambda_n \ge 0$}, \\
	{\bf T}(a_{\ell}) \times \dots \times {\bf T}(a_1) \times \ov{{\bf T}}(0)^{\times \ov{q}} \times ( {\bf T}^{\rm sp-} )^{\ov{r}} , & \textrm{if $\lambda_n < 0$},
\end{array}
\right.
\end{equation}
and regard it as a crystal by identifying ${\bf T}=(\dots,T_2,T_1)\in \widehat{\bf T}_{\lambda}$ with $T_1\otimes T_2\otimes \dots$.
Define
\begin{equation*}
	{\bf T}_{\lambda} = \{ \, {\bf T}=(\dots, T_2, T_1) \in \widehat{\bf T}_{\lambda} \, | \, T_{i+1} \prec T_i \  \textrm{for all $i$}\, \},
\end{equation*}
where $\prec$ is given in Definition \ref{def:admissibility}.
%
%
Then ${\bf T}_{\lambda}\subset \widehat{\bf T}_{\lambda}$ is invariant under $\te_i$ and $\tf_i$ for $i\in I$, and 
\begin{equation}\label{thm:fundamental_theorem}
{\bf T}_{\lambda} \cong \mathbf{B}(\omega_{\lambda}),
\end{equation}
(\cite[Theorem 4.3--4.4]{K16}).
The highest weight element ${\bf H}_{\lambda}$ of ${\bf T}_{\lambda}$ is of the form:
\begin{equation} \label{eq:hwe of spinor general}
	{\bf H}_{\lambda} =
	\left\{
\begin{array}{ll}
	 H_{\ell} \otimes \dots \otimes H_1 \otimes H_0^{\, \otimes \,q} \otimes H_{+}^{\, \otimes \,r}, & \textrm{if $\lambda_n \ge 0$}, \\
	 H_{\ell} \otimes \dots \otimes H_1 \otimes H_0^{\, \otimes \,\ov{q}} \otimes H_{-}^{\, \otimes \,\ov{r}}, & \textrm{if $\lambda_n < 0$},
\end{array}
\right.
\end{equation}
where $H_i$ and $H_{\pm}$ are the highest weight element of ${\bf T}(a_i)$ and ${\bf T}^{\rm sp\pm}$ given in \eqref{eq:hwvs in spinor}, respectively.
We call $\mathbf{T}_{\lambda}$ the {\em spinor model for $\mathbf{B}(\omega_{\lambda})$} since ${\bf T}_\la$ is a subcrystal of $({\bf T}^{\rm sp})^{\otimes N}$ for some $N\geq 1$.

\begin{ex} \label{ex:spinor}
{\rm
Let $n=8$ and $\lambda = (4, 4, 4, 4, 4, 2) \in \mc{P}_8$ be given. 
By \eqref{eq:weight expression}, we have
\begin{equation*}
	\omega_{\lambda} = 2\omega_{8-3} + 2\omega_{8-2}=2\omega_{5} + 2\omega_{6},
\end{equation*}
with $\ell=4$ and $(a_4, a_3, a_2, a_1) = (3, 3, 2, 2)$.
Let ${\bf T} = (T_4, T_3, T_2, T_1)$
given by
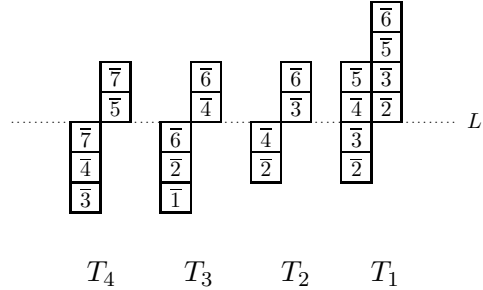
\begin{equation*}
\begin{split}
\ytableausetup {mathmode, boxsize=1.0em} 
& \begin{ytableau}
\none &\none & \none & \none & \none & \none & \none & \none & \none & \none & \none & \tl{$\ov{6}$} & \none  \\
\none &\none & \none & \none & \none & \none & \none & \none & \none & \none & \none & \tl{$\ov{5}$} & \none  \\
\none &\none & \tl{$\ov{7}$} & \none & \none & \tl{$\ov{6}$} & \none & \none & \tl{$\ov{6}$} & \none & \tl{$\ov{5}$} & \tl{$\ov{3}$} & \none \\
\none[\!\!\!\!\mathrel{\raisebox{-0.7ex}{$\scalebox{0.45}{\dots\dots\dots\dots}$}}] &\none & \tl{$\ov{5}$} & \none[\mathrel{\raisebox{-0.7ex}{$\scalebox{0.45}{\dots\dots}$}}] & \none & \tl{$\ov{4}$} & \none[\mathrel{\raisebox{-0.7ex}{$\scalebox{0.45}{\dots\dots}$}}] & \none & \tl{$\ov{3}$} & \none[\mathrel{\raisebox{-0.7ex}{$\scalebox{0.45}{\dots\dots}$}}]& \tl{$\ov{4}$} & \tl{$\ov{2}$} & \none[\ \mathrel{\raisebox{-0.7ex}{\quad $\scalebox{0.45}{\dots\dots\dots\dots}$\ ${}_{\scalebox{0.75}{$L$}}$}}] \\
\none &\tl{$\ov{7}$} & \none & \none & \tl{$\ov{6}$} & \none & \none & \tl{$\ov{4}$} & \none & \none & \tl{$\ov{3}$} & \none & \none  \\
\none &\tl{$\ov{4}$} & \none & \none & \tl{$\ov{2}$} & \none & \none & \tl{$\ov{2}$} & \none & \none & \tl{$\ov{2}$} & \none & \none  \\
\none &\tl{$\ov{3}$} & \none & \none & \tl{$\ov{1}$} & \none & \none & \none & \none & \none & \none & \none & \none  \\
\none &\none & \none & \none & \none & \none & \none & \none & \none & \none & \none & \none & \none  \\
\end{ytableau} \quad\quad \\
&\ \hskip 5mm T_4 \hskip 9mm T_3 \hskip 9mm T_2 \hskip 8mm T_1 \hskip 13mm 
\end{split}
\end{equation*}
where the dotted line is the common horizontal line $L$.
Then 
$T_4 \prec T_3 \prec T_2 \prec T_1$, and hence ${\bf T} \in {\bf T}_{\lambda}$.
}	
\end{ex}

We also need the following in Section \ref{subsec:generalization}.

\begin{df} \label{df:order}
{\em Let $\B$ be one of $\mathbf{T}(b)$ $(0\leq b <n)$, $\mathbf{T}^{\textrm{sp}}$, and $\overline{\mathbf{T}}(0)$.
For $(T, S) \in {\bf T}(a) \times {\bf B}$ with $a \in \mathbb{Z}_+$,
we write  $T \lt S$ if the pair $({}^{\texttt{R}} T, S^{\texttt{L}})$ forms a semistandard tableau of a skew shape, where we assume that ${}^{\texttt{R}} T$ and $S^{\texttt{L}}$ are arranged along $L$ as follows:
	\begin{equation*}
	\begin{split}
		& {}^{\texttt{R}} T  = (\dots, {}^{\texttt{R}} T(a+1)) \boxplus ({}^{\texttt{R}} T(a), \dots, {}^{\texttt{R}} T(1)), \\
		& S^{\texttt{L}} = (\dots, S^{\texttt{L}}(b+1)) \boxplus (S^{\texttt{L}}(b), \dots, S^{\texttt{L}}(1)).
	\end{split}
	\end{equation*}
Here we understand $S$ in the sense of Remark \ref{rem:admissibilty for spin column} and put $b = {\rm ht}(S^{\texttt{tail}})$ when $S\in {\bf T}^{\rm sp-}$ or $\ov{\bf T}(0)$.
}
\end{df}
%
\vskip 1.5mm


\subsection{Isomorphisms}
Let us give an explicit description of the isomorphisms between ${\bf KN}_\la$ and ${\bf T}_\la$ for $\la\in \mc{P}_n$ (cf. \cite[Section 3.3]{K18-2} for type $B_n$ and $C_n$).

Let $\B$ be one of $\mathbf{T}(a)$ $(0\leq a \leq n-1)$, $\mathbf{T}^{\textrm{sp}\pm}$, and $\overline{\mathbf{T}}(0)$. For $T\in \B$, we define a tableau $\td{T}$ as follows:

\begin{itemize}
\item[{\em Case 1}.] Suppose that $\B = {\bf T}^{{\rm sp}\pm}$. 
Let $\td{T}$ be the unique tableau in $SST_{[\ov{n}]}((1^n))$ such that $\ov{i}$ appears in $T$ if and only if $\ov{i}$ appears in $\widetilde{T}$ for $1\leq i\leq n$.

\item[{\em Case 2}.] Suppose that ${\bf B} = {\bf T}(a)$\, $(0 \le a < n-1)$ or $\ov{\bf T}(0)$. \\
(i) First, let $\widetilde{{}^{\texttt{R}} T}$ be the unique tableau in ${SST}_{[n]}(1^m)$ with $m = n-{\rm ht}({}^{\texttt{R}} T)$ such that $i$ appears in $\widetilde{{}^{\texttt{R}} T}$ if and only if $\ov{i}$ does not appear in ${}^{\texttt{R}}T$ for $i\in [n]$.\\
(ii) Then define $\widetilde{T}$ to be the tableau of single column obtained by putting the single-column tableau consisting of\,
\raisebox{-.6ex}{{\tiny ${\def\lr#1{\multicolumn{1}{|@{\hspace{.6ex}}c@{\hspace{.6ex}}|}{\raisebox{-.3ex}{$#1$}}}\raisebox{-.6ex}
{$\begin{array}[b]{c}
\cline{1-1}
\lr{ \ov{n} }\\
\cline{1-1}
\lr{ \!n\!}\\
\cline{1-1}
\end{array}$}}$}}
with height $b-2\mf{r}_T$ between ${}^{\texttt{L}}T$ and $\widetilde{{}^{\texttt{R}}T}$, 
where ${}^{\texttt{L}}T$ is located below $\widetilde{{}^{\texttt{R}}T}$.
\end{itemize}

\begin{ex}
{\em 
Let $n=8$ and let $T \in {\bf T}(2)$ be $T_1$ in Example \ref{ex:spinor} 
with ${\rm sh}(T) = \lambda(2, 2, 2)$ and ${\mf r}_T = 1$, where we have
\begin{equation*}
\hspace{2.5cm} T^{\texttt{L}} = \hspace{-3cm} \begin{split}
\ytableausetup {mathmode, boxsize=0.9em} 
\begin{ytableau}
	\none & \none \\
	\none & \none \\
	\scalebox{0.7}{$\ov{5}$} & \none \\
	\scalebox{0.7}{$\ov{4}$} & \none \\
	\scalebox{0.7}{$\ov{3}$} & \none \\
	\scalebox{0.7}{$\ov{2}$} & \none \\
\end{ytableau}
\end{split}\,,\quad
T^{\texttt{R}} = \begin{split}
\ytableausetup {mathmode, boxsize=0.9em} 
\begin{ytableau}
	\none & \scalebox{0.7}{$\ov{6}$} \\
	\none & \scalebox{0.7}{$\ov{5}$} \\
	\none & \scalebox{0.7}{$\ov{3}$} \\
	\none & \scalebox{0.7}{$\ov{2}$} \\
	\none & \none \\
	\none & \none \\
\end{ytableau}
\end{split}\,,\quad
{}^{\texttt{L}}T = \begin{split}
\ytableausetup {mathmode, boxsize=0.9em} 
\begin{ytableau}
	\none & \none \\
	\none & \none \\
	\none & \none \\
	\scalebox{0.7}{$\ov{5}$} & \none \\
	\scalebox{0.7}{$\ov{3}$} & \none \\
	\scalebox{0.7}{$\ov{2}$} & \none \\
\end{ytableau}
\end{split}\,,\quad
{}^{\texttt{R}}T = \begin{split}
\ytableausetup {mathmode, boxsize=0.9em} 
\begin{ytableau}
	\none & \none \\
	\none & \scalebox{0.7}{$\ov{6}$} \\
	\none & \scalebox{0.7}{$\ov{5}$} \\
	\none & \scalebox{0.7}{$\ov{4}$} \\
	\none & \scalebox{0.7}{$\ov{3}$} \\
	\none & \scalebox{0.7}{$\ov{2}$} \\
\end{ytableau}
\end{split}\,.
\end{equation*}
Then $\widetilde{{}^{\texttt{R}} T}$ is given by 
\begin{equation*}
\hspace{3.9cm} {}^{\texttt{R}}T = \hspace{-4.8cm}
\begin{split}
\ytableausetup {mathmode, boxsize=0.9em} 
\begin{ytableau}
	\none & \scalebox{0.7}{$\ov{6}$} \\
	\none & \scalebox{0.7}{$\ov{5}$} \\
	\none & \scalebox{0.7}{$\ov{4}$} \\
	\none & \scalebox{0.7}{$\ov{3}$} \\
	\none & \scalebox{0.7}{$\ov{2}$} \\
\end{ytableau}
\end{split}
\quad {\xrightarrow{\hspace*{1.2cm}}} \quad
\begin{split}
\ytableausetup {mathmode, boxsize=0.9em} 
\begin{ytableau}
	\none & \none \\
	\scalebox{0.7}{$1$} & \none \\
	\scalebox{0.7}{$7$} & \none \\
	\scalebox{0.7}{$8$} & \none \\
	\none & \none \\
\end{ytableau}
\end{split} = \, \widetilde{{}^{\texttt{R}} T}.
\end{equation*}
and hence 
\begin{equation*}
\hspace{3.9cm} T = \hspace{-4.5cm}
\begin{split}
\ytableausetup {mathmode, boxsize=0.9em} 
\begin{ytableau}
	\none & \scalebox{0.7}{$\ov{6}$} \\
	\none & \scalebox{0.7}{$\ov{5}$} \\
	\scalebox{0.7}{$\ov{5}$} & \scalebox{0.7}{$\ov{3}$} \\
	\scalebox{0.7}{$\ov{4}$} & \scalebox{0.7}{$\ov{2}$} \\
	\scalebox{0.7}{$\ov{3}$} & \none \\
	\scalebox{0.7}{$\ov{2}$} & \none \\
\end{ytableau}
\end{split}
\quad {\xrightarrow{\hspace*{1.2cm}}} \quad
\begin{split}
\ytableausetup {mathmode, boxsize=0.9em} 
\begin{ytableau}
	\scalebox{0.7}{$1$} & \none \\
	\scalebox{0.7}{$7$} & \none \\
	\scalebox{0.7}{$8$} & \none \\
	\scalebox{0.7}{$\ov{5}$} & \none \\
	\scalebox{0.7}{$\ov{3}$} & \none \\
	\scalebox{0.7}{$\ov{2}$} & \none \\
\end{ytableau}
\end{split} = \,\widetilde{T}.
\end{equation*}
Note that since $b = 2$ and ${\mf r}_T=1$, there is no domino 
\raisebox{-.6ex}{{\tiny ${\def\lr#1{\multicolumn{1}{|@{\hspace{.6ex}}c@{\hspace{.6ex}}|}{\raisebox{-.3ex}{$#1$}}}\raisebox{-.6ex}
{$\begin{array}[b]{c}
\cline{1-1}
\lr{ \ov{n} }\\
\cline{1-1}
\lr{ \!n\!}\\
\cline{1-1}
\end{array}$}}$}}
in $\widetilde{T}$.
On the other hand,
if $T \in {\bf T}^{{\rm sp}-}$ is given as follows, then we have 
\begin{equation*}
\hspace{4.2cm} T  = \hspace{-5cm}
\begin{split}
	\ytableausetup {mathmode, boxsize=0.9em} 
\begin{ytableau}
	\none  \\
	\scalebox{0.7}{$\ov{5}$}  \\
	\scalebox{0.7}{$\ov{4}$}  \\
	\scalebox{0.7}{$\ov{1}$}  \\
	\none \\
\end{ytableau}
\end{split}
\quad {\xrightarrow{\hspace*{1.2cm}}} \quad
\begin{split}
	\ytableausetup {mathmode, boxsize=0.9em} 
\begin{ytableau}
	\scalebox{0.7}{$2$}  \\
	\scalebox{0.7}{$3$}  \\
	\scalebox{0.7}{$6$}  \\
	\scalebox{0.7}{$7$}  \\
	\scalebox{0.7}{$8$}  \\
	\scalebox{0.7}{$\ov{5}$}  \\
	\scalebox{0.7}{$\ov{4}$}  \\
	\scalebox{0.7}{$\ov{1}$}  \\
\end{ytableau}
\end{split} \,=\, \widetilde{T}.
\end{equation*}
}
\end{ex}

\begin{lem} \label{lem: T(a) 2 <= a <= n-1}
The map sending $T$ to $\td{T}$ gives an isomorphisms of crystals
\begin{equation*} 
	\xymatrixcolsep{3pc}\xymatrixrowsep{0pc}
	\xymatrix{ \Phi:
	{\bf B}  \ \ar@{->}[r] &  
	{\begin{cases} 
	 {\bf KN}_{{\texttt {\em sp}}\pm}, & \ \text{if $\B= {\bf T}^{{\rm sp}\pm}$}, \\
	 {\bf KN}_{(1^{n-a})}, & \ \text{if ${\bf B} = {\bf T}(a)$ or $\ov{\bf T}(0)$}. \\
	\end{cases}}
}
\end{equation*} 
\end{lem}

\pf
{\em Case 1.} $\B={\bf T}^{\rm sp\pm}$.
It is straightforward to see that $\Phi$ is a weight preserving bijection and by \eqref{eq:Kashiwara operator on spin crystal}, it is a morphism of crystals. 
Hence $\Phi$ is an isomorphism.
\vskip 2.5mm

{\em Case 2.} Suppose that ${\bf B} = {\bf T}(a)$\, $(0 \le a < n-1)$ or $\ov{\bf T}(0)$. Let $T \in {\bf B}$ given.
First we claim $\widetilde{T} \in {\bf KN}_{(1^{n-a})}$. 
Suppose that $\widetilde{T} \notin {\bf KN}_{(1^{n-a})}$. Then by the condition (${\mf d}$-1)
there exists $i \in [n]$ such that 
\begin{equation} \label{eq1:T to KN}
(q-p)+i \le n-a \quad (p < q).
\end{equation}
Put $x = n-a-q$ and $y = p$. We note that $x$ is the number of entries in $[n]$ smaller than $i$ in $\widetilde{T}$, and $y$ is the number of entries in $[\ov{n}]$ equal to or larger than $\ov{i}$ in $\widetilde{T}$.
Take $k$ such that ${}^{\texttt{L}}T(k) = \ov{i}$. Then we have
${}^{\texttt{L}}T(k) > {}^{\texttt{R}}T(k)$  by \eqref{eq1:T to KN}. This contradicts to the fact that the pair $({}^{\texttt{L}} T, {}^{\texttt{R}} T)$ forms a semistandard tableau when the two columns are placed on the common bottom line.
Hence $\widetilde{T} \in {\bf KN}_{(1^{n-a})}$.

Second we show that $T \equiv_{\mf l} \widetilde{T}$, where $\equiv_{\mf l}$ denotes the crystal equivalence as elements of $\mf l$-crystals. By the construction of $\widetilde{{}^{\texttt{R}}T}$, it is not difficult to check that
$\widetilde{{}^{\texttt{R}}T} \equiv_{\mf l} {}^{\texttt{R}}T$ (more precisely as elements of $\mf{sl}_n$-crystals).
Put 
$\texttt{D}_0$
to be the single-column tableau consisting of the domino 
\raisebox{-.6ex}{{\tiny ${\def\lr#1{\multicolumn{1}{|@{\hspace{.6ex}}c@{\hspace{.6ex}}|}{\raisebox{-.3ex}{$#1$}}}\raisebox{-.6ex}
{$\begin{array}[b]{c}
\cline{1-1}
\lr{ \ov{n} }\\
\cline{1-1}
\lr{ \!n\!}\\
\cline{1-1}
\end{array}$}}$}}
with height $b-2\mf{r}_T$. 
By the tensor product rule of crystals, we see that $\{\texttt{D}_0\}$ is the crystal of the trivial representation of $\mf l$. 
This implies that 
$\widetilde{{}^{\texttt{R}} T}  \equiv_{\mf l}\widetilde{{}^{\texttt{R}} T} \otimes \texttt{D}_0 $ and thus 
$$
	T \equiv_{\mf l} {}^{\texttt{R}} T \otimes {}^{\texttt{L}} T \equiv_{\mf l} \widetilde{{}^{\texttt{R}} T} \otimes {}^{\texttt{L}} T \equiv_{\mf l}
	\widetilde{{}^{\texttt{R}} T} \otimes \texttt{D}_0 \otimes {}^{\texttt{L}} T \equiv_{\mf l} \widetilde{T}
$$

Next we claim that
$\widetilde{T'} = \widetilde{f}_n \widetilde{T}$, 
where $T' = \widetilde{f}_n T$. 
Let $T \in {SST}_{[\ov{n}]}(\lambda(a, b, c))$ and $T' \in {SST}_{[\ov{n}]}(\lambda(a, b', c'))$. 
Let us consider the case when
$\widetilde{f}_n (T^{\texttt{R}} \otimes T^{\texttt{L}}) \neq {\bf 0}$ and 
$\widetilde{f}_n (T^{\texttt{R}} \otimes T^{\texttt{L}}) =  T^{\texttt{R}} \otimes (\widetilde{f}_n T^{\texttt{L}})$. The proof of the other cases is similar.
In this case, we have $b' = b-2$, $c' = c+2$ by definition of $\widetilde{f}_n$, and that $T^{\texttt{L}}[1]$ and $T^{\texttt{R}}[1]$ must satisfy that
	\begin{equation*}
		\ov{\Tl{$n-2$}}\, \le T^{\texttt{L}}[1], \, \quad \, T^{\texttt{R}}[1] \le \ov{\Tl{$n-1$}}.
	\end{equation*}
	%
	Then we have
	\begin{equation*}
		\widetilde{f}_n \widetilde{T} = 
		\left\{ \begin{array}{ll}
	\widetilde{{}^{\texttt{R}}T} \otimes \widetilde{f}_n(\texttt{D}_0) \otimes {}^{\texttt{L}} T, & \textrm{if $b > 2$}, \\
	\widetilde{f}_n(\widetilde{{}^{\texttt{R}}T})\otimes  {}^{\texttt{L}} T, & \textrm{if $b = 2$},
	\end{array}
		\right.
	\end{equation*}
	where $\widetilde{f}_n(\texttt{D}_0)$ is obtained from $\texttt{D}_0$ by replacing 
	\raisebox{-.6ex}{{\tiny ${\def\lr#1{\multicolumn{1}{|@{\hspace{.6ex}}c@{\hspace{.6ex}}|}{\raisebox{-.3ex}{$#1$}}}\raisebox{-.6ex}
{$\begin{array}[b]{c}
\cline{1-1}
\lr{ \!\ov{n}\! }\\
\cline{1-1}
\lr{ \!n\!}\\
\cline{1-1}
\end{array}$}}$}}
	by 
	\raisebox{-.6ex}{{\tiny ${\def\lr#1{\multicolumn{1}{|@{\hspace{.6ex}}c@{\hspace{.6ex}}|}{\raisebox{-.3ex}{$#1$}}}\raisebox{-.6ex}
{$\begin{array}[b]{c}
\cline{1-1}
\lr{ \!\ov{n}\!}\\
\cline{1-1}
\lr{ \!\ov{n-1}\!}\\
\cline{1-1}
\end{array}$}}$}}
	at the bottom of $\texttt{D}_0$.
	Note that the bottom entry of $\widetilde{f}_n(\widetilde{{}^{\texttt{R}}T})$ is given by
	\begin{equation*} \label{eq:top of f R T}
		\left\{ 
		\begin{array}{ll}
	\raisebox{-.6ex}{{\tiny ${\def\lr#1{\multicolumn{1}{|@{\hspace{.6ex}}c@{\hspace{.6ex}}|}{\raisebox{-.1ex}{$#1$}}}\raisebox{0.6ex}
{$\begin{array}[b]{c}
\cline{1-1}
\lr{ \!\ov{n}\! }\\
\cline{1-1}
\end{array}$}}$}},
	& \textrm{if $T^{\texttt{R}}[1] = \ov{n}$}, \\
	\raisebox{-.6ex}{{\tiny ${\def\lr#1{\multicolumn{1}{|@{\hspace{.6ex}}c@{\hspace{.6ex}}|}{\raisebox{-.5ex}{$#1$}}}\raisebox{0.6ex}
{$\begin{array}[b]{c}
\cline{1-1}
\lr{ \!\ov{n-1}\! }\\
\cline{1-1}
\end{array}$}}$}},
	& \textrm{if $T^{\texttt{R}}[1] = \ov{\Tl{$n-1$}}$}.
		\end{array}
		\right.
	\end{equation*}
	On the other hand, we can check that ${}^{\texttt{L}} T' $ is obtained from ${}^{\texttt{L}}T$ by putting the domino 
	\raisebox{-.6ex}{{\tiny ${\def\lr#1{\multicolumn{1}{|@{\hspace{.6ex}}c@{\hspace{.6ex}}|}{\raisebox{-.3ex}{$#1$}}}\raisebox{-.6ex}
{$\begin{array}[b]{c}
\cline{1-1}
\lr{ \!\ov{n}\!}\\
\cline{1-1}
\lr{ \!\ov{n-1}\!}\\
\cline{1-1}
\end{array}$}}$}}
	(resp.
	\raisebox{-.6ex}{{\tiny ${\def\lr#1{\multicolumn{1}{|@{\hspace{.6ex}}c@{\hspace{.6ex}}|}{\raisebox{-.1ex}{$#1$}}}\raisebox{0.6ex}
{$\begin{array}[b]{c}
\cline{1-1}
\lr{ \!\ov{n}\! }\\
\cline{1-1}
\end{array}$}}$}}
	if $T^{\texttt{R}}[1]=\ov{n}$, and
	\raisebox{-.6ex}{{\tiny ${\def\lr#1{\multicolumn{1}{|@{\hspace{.6ex}}c@{\hspace{.6ex}}|}{\raisebox{-.5ex}{$#1$}}}\raisebox{0.6ex}
{$\begin{array}[b]{c}
\cline{1-1}
\lr{ \!\ov{n-1}\! }\\
\cline{1-1}
\end{array}$}}$}}
	if $T^{\texttt{R}}[1]=\ov{\Tl{$n-1$}}$
	)
	on the top of ${}^{\texttt{L}}T$ when $b > 2$ (resp. $b=2$).
	Now it is easy to see that $\widetilde{f}_n \widetilde{T}$ is equal to $\widetilde{T'}$.
	%
	%

Consequently, $\Phi$
is a morphism of crystals. Since $\Phi$ is injective and sends the highest weight elements of ${\bf T}(a)$ to that of ${\bf KN}_{(1^{n-a})}$ , $\Phi$ is an isomorphism.
\qed
\vskip 2mm

Next let us describe the inverse map of $\Phi$. Let $T \in {\bf KN}_{(1^a)}\cup {\bf KN}_{\texttt{sp}^{\pm}}$ ($0<a\leq n$) be given. Then we define $\Psi_a(T)$ and $\Psi_{{\rm sp}\pm}(T)$ if $T\in {\bf KN}_{(1^a)}$ and $T\in {\bf KN}_{\texttt{sp}^{\pm}}$ respectively as follows:
\begin{itemize}
	\item[(1)] Let $T_+$ (resp. $T_-$) be the subtableau in $T$ with entries in $[n]$ (resp. $[\ov{n}]$) except for dominos
\raisebox{-.6ex}{{\tiny ${\def\lr#1{\multicolumn{1}{|@{\hspace{.6ex}}c@{\hspace{.6ex}}|}{\raisebox{-.3ex}{$#1$}}}\raisebox{-.6ex}
{$\begin{array}[b]{c}
\cline{1-1}
\lr{ \ov{n} }\\
\cline{1-1}
\lr{ \!n\!}\\
\cline{1-1}
\end{array}$}}$}}

	\item[(2)] Let $\widetilde{T_+}$ be the single-column tableau with height $n-{\rm ht}(T_+)$ such that $i$ appears in $T_+$ if and only if $\ov{i}$ does not appear in $\widetilde{T_+}$.
	
	\item[(3)] 
%
We define $\Psi_a(T)$ and $\Psi_{{\rm sp}\pm}(T)$ by
\begin{equation}\label{eq:map from KN to spinor for fundamentals}
\begin{split}
\Psi_{{\rm sp}\pm}(T)  = T_- ,\quad \Psi_a(T) & = \mathcal{F}^{n-a-\epsilon}(T_-, \widetilde{T_+}), \\
\end{split}
\end{equation}
where 
\begin{equation*}
	\epsilon = 
	\left\{
\begin{array}{ll}
	0  , & \textrm{if ${\rm ht}(\widetilde{T_+})-a$ is even}, \\
	1 , & \textrm{if ${\rm ht}(\widetilde{T_+})-a$ is odd} .
\end{array}
\right.
\end{equation*}
We note that $\Psi_a(T)$ has residue $1$ if ${\rm ht}(\widetilde{T_+})-a$ is odd, otherwise $0$.
\end{itemize}
It is not difficult to check that the map $\Psi_a$ (resp. $\Psi_{\rm sp{\pm}}$) is the inverse of $\Phi$. Hence by Lemma \ref{lem: T(a) 2 <= a <= n-1}, we have the following.
%

\begin{lem} \label{lem:KN to spin}
The maps $\Psi_a$ and $\Psi_{\rm sp{{\pm}}}$ are isomorphisms of crystals
\begin{equation*}
\begin{split}
\Psi_{{\rm sp}\pm} : \xymatrixcolsep{3pc}\xymatrixrowsep{0pc}
\xymatrix{
{\bf KN}_{\texttt{\em sp}\pm}  \ar@{->}[r] & {\bf T}^{{\rm sp}\pm}},\quad
\Psi_a : \xymatrixcolsep{3pc}\xymatrixrowsep{0pc}
\xymatrix{
{\bf KN}_{(1^a)}  \ar@{->}[r] & {\bf T}(n-a)} \quad (0<a\leq n),
\end{split}	
\end{equation*}
\end{lem}
\qed
\vskip 1.5mm


Now we consider the isomorphism for any $\lambda \in \mc{P}_n$.
Let 
$\mu'=(a_{\ell}, \dots, a_1)$, where $a_1, \dots, a_{\ell}$ are given in \eqref{eq:weight expression}, with $\ell = \mu_1$.
For $T \in {\bf KN}_{\lambda}$, let
\begin{equation*}
	\left\{
\begin{array}{ll}
	(T_{\ell}, \dots, T_1) , & \textrm{if $\lambda_n\in\Z$}, \\
	(T_{\ell}, \dots, T_1, T_0),	& \textrm{if $\lambda_n\not\in\Z$},
\end{array}
\right.	
\end{equation*}
denote the sequence of columns of $T$, where $T_0$ is the column of $T$ with half-width boxes, and $T_1,T_2,\dots$ are the other columns enumerated from right to left.
\begin{thm} \label{thm:KN to spinor}
For $\lambda \in \mc{P}_n$, the map
\begin{equation} \label{eq:Psi}
\xymatrixcolsep{3pc}\xymatrixrowsep{0pc}
\xymatrix{
\Psi_{\lambda} : {\bf KN}_{\lambda} \ \ar@{->}[r] & {\bf T}_{\lambda}
}	
\end{equation}
defined by
\begin{equation*}
	\Psi_{\lambda}(T) =
	\left\{
\begin{array}{ll}
	(\Psi_{\mu_{\ell}'}(T_{\ell}), \dots, \Psi_{\mu_1'}(T_1)  ) , & \textrm{if $\lambda_n\in\Z$}, \\
	(\Psi_{\mu_{\ell}'}(T_{\ell}), \dots, \Psi_{\mu_1'}(T_1), \Psi_{\rm sp{\pm}}(T_0)),	& \textrm{if $\lambda_n\not\in\Z$},
\end{array}
\right.
\end{equation*}
is an isomorphism of crystals from ${\bf KN}_{\lambda}$ to ${\bf T}_{\lambda}$,
where  we take $\Psi_{\rm sp+}$ and $\Psi_{\rm sp-}$ if $\lambda_n \ge 0$ and $\lambda_n < 0$, respectively.

\end{thm}
\pf
By Lemma \ref{lem:KN to spin}, the map $\Psi_{\lambda}$ is an embedding of crystals into $\widehat{\bf T}_{\lambda}$.
Also the map $\Psi_{\lambda}$
sends the highest weight element of ${\bf KN}_{\lambda}$ to the one of ${\bf T}_{\lambda}$ (cf. \eqref{eq:hwe of spinor general}).
Then by \eqref{thm:fundamental_theorem}, the image of $\Psi_{\lambda}$ is isomorphic to ${\bf T}_{\lambda}$.
\qed

\begin{ex} \label{ex:KN to spinor}
{\rm
Let $\lambda = (4^5, 2) \in \mc{P}_8$ be given. Consider
\begin{equation*}
\hspace{5.5cm} T = \hspace{-5.3cm} 
\begin{split}
\ytableausetup {mathmode, boxsize=0.8em} 
\begin{ytableau}
 \none & \none & \none & \none \\
 \none & \none & \tl{$1$} & \tl{$1$} \\
 \tl{$1$} & \tl{$3$} & \tl{$4$} & \tl{$7$} \\
 \tl{$2$} & \tl{$5$} & \tl{$5$} & \tl{$8$} \\
 \tl{$6$} & \tl{$7$} & \tl{$7$} & \tl{$\ov{5}$} \\
 \tl{$8$} & \tl{$8$} & \tl{$8$} & \tl{$\ov{3}$} \\
 \tl{$\ov{7}$} & \tl{$\ov{6}$} & \tl{$\ov{4}$} & \tl{$\ov{2}$}
\end{ytableau}
\end{split} \,\, \in \,\, {\bf KN}_{\lambda}\,.
\end{equation*}
Put
\begin{equation*}
\hspace{3.5cm} \scalebox{0.8}{$T_4$} = \hspace{-3.9cm} \begin{split}
\ytableausetup {mathmode, boxsize=0.8em} 
\begin{ytableau}
	\tl{$1$} \\
	\tl{$2$} \\
	\tl{$6$} \\
	\tl{$8$} \\
	\tl{$\ov{7}$}
\end{ytableau}
\end{split}\,, \quad
\scalebox{0.8}{$T_3$} = \begin{split}
\ytableausetup {mathmode, boxsize=0.8em} 
\begin{ytableau}
	\tl{$3$} \\
	\tl{$5$} \\
	\tl{$7$} \\
	\tl{$8$} \\
	\tl{$\ov{6}$}
\end{ytableau}
\end{split}\,, \quad
\scalebox{0.8}{$T_2$} = \begin{split}
\ytableausetup {mathmode, boxsize=0.8em} 
\begin{ytableau}
	\tl{$1$} \\
	\tl{$4$} \\
	\tl{$5$} \\
	\tl{$7$} \\
	\tl{$8$} \\
	\tl{$\ov{4}$}
\end{ytableau}
\end{split}\,, \quad
\scalebox{0.8}{$T_1$} = \begin{split}
\ytableausetup {mathmode, boxsize=0.8em} 
\begin{ytableau}
	\tl{$1$} \\
	\tl{$7$} \\
	\tl{$8$} \\
	\tl{$\ov{5}$} \\
	\tl{$\ov{3}$} \\
	\tl{$\ov{2}$}
\end{ytableau}
\end{split}\,.
\end{equation*}
By definition of $((T_i)_- \,, \widetilde{(T_i)_+})$,
\begin{equation*}
\begin{split}
& (\scalebox{0.8}{$(T_4)_-$}\,, \scalebox{0.8}{$\widetilde{(T_4)_+}$}) =\, 
\Bigg(\, \raisebox{-.8ex}{{\tiny ${\def\lr#1{\multicolumn{1}{|@{\hspace{1ex}}c@{\hspace{1ex}}|}{\raisebox{-.3ex}{$#1$}}}\raisebox{-3.5ex}
{$\begin{array}[b]{c}
\cline{1-1}
\lr{ \!\ov{7}\! }\\
\cline{1-1}
\end{array}$}}$}}
\raisebox{-2.5ex}{\,\scalebox{0.9}{,}} \,\, 
\raisebox{-.8ex}{{\tiny ${\def\lr#1{\multicolumn{1}{|@{\hspace{1ex}}c@{\hspace{1ex}}|}{\raisebox{-.3ex}{$#1$}}}\raisebox{-3.5ex}
{$\begin{array}[b]{c}
\cline{1-1}
\lr{ \!\ov{7}\! }\\
\cline{1-1}
\lr{ \!\ov{5}\! }\\
\cline{1-1}
\lr{ \!\ov{4}\! }\\
\cline{1-1}
\lr{ \!\ov{3}\! }\\
\cline{1-1}
\end{array}$}}$}} \,\, \Bigg) \,, \quad
(\scalebox{0.8}{$(T_3)_-$}, \scalebox{0.8}{$\widetilde{(T_3)_+}$}) =\, 
\Bigg( \, \raisebox{-.8ex}{{\tiny ${\def\lr#1{\multicolumn{1}{|@{\hspace{1ex}}c@{\hspace{1ex}}|}{\raisebox{-.3ex}{$#1$}}}\raisebox{-3.5ex}
{$\begin{array}[b]{c}
\cline{1-1}
\lr{ \!\ov{6}\! }\\
\cline{1-1}
\end{array}$}}$}}
\raisebox{-2.5ex}{\,\scalebox{0.9}{,}} \,\,  
\raisebox{-.8ex}{{\tiny ${\def\lr#1{\multicolumn{1}{|@{\hspace{1ex}}c@{\hspace{1ex}}|}{\raisebox{-.3ex}{$#1$}}}\raisebox{-3.5ex}
{$\begin{array}[b]{c}
\cline{1-1}
\lr{ \!\ov{6}\! }\\
\cline{1-1}
\lr{ \!\ov{4}\! }\\
\cline{1-1}
\lr{ \!\ov{2}\! }\\
\cline{1-1}
\lr{ \!\ov{1}\! }\\
\cline{1-1}
\end{array}$}}$}} \, \Bigg)\,, \\
\\
& (\scalebox{0.8}{$(T_2)_-$}\,, \scalebox{0.8}{$\widetilde{(T_2)_+}$}) =\, 
\Bigg(\, \raisebox{-.8ex}{{\tiny ${\def\lr#1{\multicolumn{1}{|@{\hspace{1ex}}c@{\hspace{1ex}}|}{\raisebox{-.3ex}{$#1$}}}\raisebox{-3.5ex}
{$\begin{array}[b]{c}
\cline{1-1}
\lr{ \!\ov{4}\! }\\
\cline{1-1}
\end{array}$}}$}}
\raisebox{-2.5ex}{\,\scalebox{0.9}{,}} \,\, 
\raisebox{-.8ex}{{\tiny ${\def\lr#1{\multicolumn{1}{|@{\hspace{1ex}}c@{\hspace{1ex}}|}{\raisebox{-.3ex}{$#1$}}}\raisebox{-3.5ex}
{$\begin{array}[b]{c}
\cline{1-1}
\lr{ \!\ov{6}\! }\\
\cline{1-1}
\lr{ \!\ov{3}\! }\\
\cline{1-1}
\lr{ \!\ov{2}\! }\\
\cline{1-1}
\end{array}$}}$}} \,\, \Bigg) \,, \quad
(\scalebox{0.8}{$(T_1)_-$}\,, \scalebox{0.8}{$\widetilde{(T_1)_+}$}) =\, 
\Bigg(\, \raisebox{-.8ex}{{\tiny ${\def\lr#1{\multicolumn{1}{|@{\hspace{1ex}}c@{\hspace{1ex}}|}{\raisebox{-.3ex}{$#1$}}}\raisebox{-3.5ex}
{$\begin{array}[b]{c}
\cline{1-1}
\lr{ \!\ov{5}\! }\\
\cline{1-1}
\lr{ \!\ov{3}\! }\\
\cline{1-1}
\lr{ \!\ov{2}\! }\\
\cline{1-1}
\end{array}$}}$}}
\raisebox{-2.5ex}{\,\scalebox{0.9}{,}} \,\, 
\raisebox{-.8ex}{{\tiny ${\def\lr#1{\multicolumn{1}{|@{\hspace{1ex}}c@{\hspace{1ex}}|}{\raisebox{-.3ex}{$#1$}}}\raisebox{-3.5ex}
{$\begin{array}[b]{c}
\cline{1-1}
\lr{ \!\ov{6}\! }\\
\cline{1-1}
\lr{ \!\ov{5}\! }\\
\cline{1-1}
\lr{ \!\ov{4}\! }\\
\cline{1-1}
\lr{ \!\ov{3}\! }\\
\cline{1-1}
\lr{ \!\ov{2}\! }\\
\cline{1-1}
\end{array}$}}$}} \,\, \Bigg).
\end{split}
\end{equation*}
Since ${\rm ht}(\widetilde{(T_i)_+})-{\rm ht}(T_i)$ is odd for $i=1,2,3,4$, we have by \eqref{eq:map from KN to spinor for fundamentals}

\begin{equation*}
\hspace{2.3cm} \scalebox{0.85}{$\Psi_5(T_4)$} = \hspace{-1.7cm}
\begin{split}
\ytableausetup {mathmode, boxsize=0.8em} 
\begin{ytableau}
	\none & \none \\
	\none & \none \\
	\none & \scalebox{0.65}{$\ov{7}$} \\
	\none & \scalebox{0.65}{$\ov{5}$} \\
	\scalebox{0.65}{$\ov{7}$} & \none \\
	\scalebox{0.65}{$\ov{4}$} & \none \\
	\scalebox{0.65}{$\ov{3}$} & \none \\
\end{ytableau}
\end{split}\,, \quad
\scalebox{0.85}{$\Psi_5(T_3)$} = \,
\begin{split}
\ytableausetup {mathmode, boxsize=0.8em} 
\begin{ytableau}
	\none & \none \\
	\none & \none \\
	\none & \scalebox{0.65}{$\ov{6}$} \\
	\none & \scalebox{0.65}{$\ov{4}$} \\
	\scalebox{0.65}{$\ov{6}$} & \none \\
	\scalebox{0.65}{$\ov{2}$} & \none \\
	\scalebox{0.65}{$\ov{1}$} & \none \\
\end{ytableau}
\end{split}\,, \quad
\scalebox{0.85}{$\Psi_6(T_2)$} = \,
\begin{split}
\ytableausetup {mathmode, boxsize=0.8em} 
\begin{ytableau}
	\none & \none \\
	\none & \none \\
	\none & \scalebox{0.65}{$\ov{6}$} \\
	\none & \scalebox{0.65}{$\ov{3}$} \\
	\scalebox{0.65}{$\ov{4}$} & \none \\
	\scalebox{0.65}{$\ov{2}$} & \none \\
	\none & \none \\
\end{ytableau}
\end{split}\,, \quad
\scalebox{0.85}{$\Psi_6(T_1)$} = \,
\begin{split}
\ytableausetup {mathmode, boxsize=0.8em} 
\begin{ytableau}
	\none & \scalebox{0.65}{$\ov{6}$} \\
	\none & \scalebox{0.65}{$\ov{5}$} \\
	\scalebox{0.65}{$\ov{5}$} & \scalebox{0.65}{$\ov{3}$} \\
	\scalebox{0.65}{$\ov{4}$} & \scalebox{0.65}{$\ov{2}$} \\
	\scalebox{0.65}{$\ov{3}$} & \none \\
	\scalebox{0.65}{$\ov{2}$} & \none \\
	\none & \none \\
\end{ytableau}
\end{split}\,. \quad \quad \quad
\end{equation*}
Hence, we obtain
\begin{equation*}
\hspace{3.5cm} \Psi_{\lambda}(T) =\quad  \hspace{-3.5cm} 
\begin{split}
\ytableausetup {mathmode, boxsize=0.9em} 
& \begin{ytableau}
\none &\none & \none & \none & \none & \none & \none & \none & \none & \none & \none & \tl{$\ov{6}$} & \none  \\
\none &\none & \none & \none & \none & \none & \none & \none & \none & \none & \none & \tl{$\ov{5}$} & \none  \\
\none &\none & \tl{$\ov{7}$} & \none & \none & \tl{$\ov{6}$} & \none & \none & \tl{$\ov{6}$} & \none & \tl{$\ov{5}$} & \tl{$\ov{3}$} & \none \\
\none[\!\!\!\!\!\!\mathrel{\raisebox{-0.6ex}{$\scalebox{0.45}{\dots\dots\dots}$}}] &\none & \tl{$\ov{5}$} & \none[\mathrel{\raisebox{-0.6ex}{$\scalebox{0.45}{\dots\dots}$}}] & \none & \tl{$\ov{4}$} & \none[\mathrel{\raisebox{-0.6ex}{$\scalebox{0.45}{\dots\dots}$}}] & \none & \tl{$\ov{3}$} & \none[\mathrel{\raisebox{-0.6ex}{$\scalebox{0.45}{\dots\dots}$}}]& \tl{$\ov{4}$} & \tl{$\ov{2}$} & \none[\ \ \mathrel{\raisebox{-0.6ex}{\quad $\scalebox{0.45}{\dots\dots\dots\dots}$\ }}] \\
\none &\tl{$\ov{7}$} & \none & \none & \tl{$\ov{6}$} & \none & \none & \tl{$\ov{4}$} & \none & \none & \tl{$\ov{3}$} & \none & \none  \\
\none &\tl{$\ov{4}$} & \none & \none & \tl{$\ov{2}$} & \none & \none & \tl{$\ov{2}$} & \none & \none & \tl{$\ov{2}$} & \none & \none  \\
\none &\tl{$\ov{3}$} & \none & \none & \tl{$\ov{1}$} & \none & \none & \none & \none & \none & \none & \none & \none  \\
\none &\none & \none & \none & \none & \none & \none & \none & \none & \none & \none & \none & \none  \\
\end{ytableau} \quad\quad  \\
\end{split}\in {\bf T}_\la.
\end{equation*}
}	
\end{ex}


\section{Embedding into the crystal of parabolic Verma module}\label{sec:spinor into parabolic verma}

In this section, we describe an embedding of ${\bf T}_\la$ into the crystal of a parabolic Verma module corresponding to $\la\in \mc{P}_n$ with respect to $\mf{l}$. The embedding maps ${\bf T}\in {\bf T}_\la$ to a pair of semistandard tableaux, which can be described in terms of a combinatorial algorithm called {\em separation}.

\subsection{Sliding} \label{subsec:generalization}
Suppose that $\la\in \mc{P}_n$ is given.
The algorithm of separation consists of certain elementary steps denoted by $\mc{S}_j$ ($j=2,4,\dots$) as operators on ${\bf T}\in {\bf T}_\la$, each of which moves a tail in ${\bf T}$ by one position left based on the Sch\"{u}tzeberger's jeu de taquin. 
%

Note that ${\bf T}_\la$ is a subcrystal of $\left({\bf T}^{\rm sp}\right)^{\otimes N}$ for some $N$. We may identify $\left({\bf T}^{\rm sp}\right)^{\otimes N}$ with
\begin{equation*}
	\mathbf{E}^N := \underset{{(u_N, \dots, u_1) \in \mathbb{Z}_+^n}}{\bigsqcup} 
	{SST}_{[\ov{n}]}(1^{u_N}) \times \dots \times {SST}_{[\ov{n}]}(1^{u_1}).
\end{equation*}
In order to give a precise description of the operator $\mc{S}_j$, we use an $(\mf{l},\mf{sl}_N)$-bicrystal structure on ${\bf E}^N$ as in \cite[Lemma 5.1]{K18-2}:
The ${\mf l}$-crystal structure on ${\bf E}^N$ with respect to $\widetilde{e}_i$ and $\widetilde{f}_i$ for $i \in I$ is naturally induced from that of $\left({\bf T}^{\rm sp}\right)^{\otimes N}$. 
On the other hand, the $\mf{sl}_N$-crystal structure is defined as follows.
Let $(U_N, \dots, U_1) \in \mathbf{E}^N$ given.
For $1 \le j \le N-1$ and $\mathcal{X} = \mathcal{E}, \mathcal{F}$, we define
\begin{equation*}
	\mathcal{X}_j (U_N, \dots, U_1) = 
	\begin{cases}
	(U_r, \dots, \mathcal{X}(U_{j+1}, U_j), \dots, U_1), 
	& \textrm{if $\mc{X}(U_{j+1}, U_j) \neq \mathbf{0}$,} \\
	\mathbf{0}, &  \textrm{if $\mc{X}(U_{j+1}, U_j) = \mathbf{0}$}.
	\end{cases}
\end{equation*} 
where $\mathcal{X}(U_{j+1}, U_j)$ is understood to be $\mc{X}U$ for some $U\in SST_{[\ov{n}]}(\la(a,b,c))$ with ${\mf r}_U=0$, $U^L = U_{j+1}$ and $U^R = U_j$ (see Section \ref{subsec:spinor}).

Let $l$ be the number of components in $\widehat{{\bf T}}_\la$ in \eqref{eq:hat{T}lambda} except ${\bf T}^{\rm sp \pm}$. Consider an embedding of sets
\begin{equation}\label{eq:identification}
\xymatrixcolsep{3pc}\xymatrixrowsep{0pc}\xymatrix{
\quad \quad {\bf T}_{\lambda}  \ \ar@{->}[r] & \ \mathbf{E}^{2l+1} \quad  \\
\mathbf{T} = (T_{l}, \dots, T_1, T_0) \ar@{|->}[r] & (T_{l}^{\texttt{L}}, T_{l}^{\texttt{R}}, \dots, T_1^{\texttt{L}}, T_1^{\texttt{R}}, T_0) },
\end{equation}
where $T_0$ is regarded as 
\begin{equation*}
	T_0 \in
\left\{
\begin{array}{ll}
	\{ \, \emptyset \, \} , & \textrm{if $\la_n\in\Z$}, \\
	{\bf T}^{\rm sp^{\pm}},	& \textrm{if $\la_n\not\in\Z$}.
\end{array}
\right.
\end{equation*}
Here $\{\emptyset\}$ is the crystal of trivial module.
We identify $\mathbf{T} = (T_{l}, \dots, T_1, T_0)\in {\bf T}_{\lambda}$ with its image ${\bf U}=(U_{2l},\dots,U_1, U_0)$ under \eqref{eq:identification} so that $T_0 = U_0$ and $(T_{i+1}, T_i)$ is given by
\begin{equation}\label{eq:identification-1}
(T_{i+1}, T_i)
=(U_{j+2}, U_{j+1}, U_j, U_{j-1})
=(T^{\texttt{L}}_{i+1},T^{\texttt{R}}_{i+1},T^{\texttt{L}}_i,T^{\texttt{R}}_i),
\end{equation}
with $j=2i$ for $1\leq i\leq l-1$.

Now we define an operator $\mc{S}_j$ on ${\bf T}$ for $j=2i$ for $1\leq i\leq l-1$ by 
\begin{equation} \label{eq:definition_S}
\mathcal{S}_j = 
\begin{cases}
\quad \quad \quad \mathcal{F}_j^{a_i}, & \textrm{if $T_{i+1}  \lt T_i$}, \\
\mathcal{E}_j\mathcal{E}_{j-1}\mathcal{F}_j^{a_i-1}\mathcal{F}_{j-1}, & \textrm{if $T_{i+1} \not\triangleleft T_i$},			
\end{cases}
\end{equation}
where $\mathcal{S}_j$ is understood as the identity operator when $a_i = 0$, and $\lt$ is given in Definition \ref{df:order}.

\begin{rem} \label{rem:generalization}
{\rm
The operator $\mc{S}_j$ in \eqref{eq:definition_S} agrees with the one in \cite{JK19-2}, which is defined only on the set of $\mf l$-highest weight elements.}	
\end{rem}

The following lemma is crucial in Section \ref{subsec:separation}.

\begin{lem} \label{lem:sliding}
Let ${\bf T} = (\dots, T_{i+1}, T_i, \dots) \in {\bf T}_{\lambda}$ be given.
\begin{itemize}
	\item[(1)] We have  
$
	\mathcal{S}_j {\bf T} = (\dots, U_{j+2}, \widetilde{U}_{j+1}, \widetilde{U}_j, U_{j-1}, \dots)
$ for some $\td{U}_{j+1}$ and $\td{U}_j$, and $($$U_{j+2}$, $\widetilde{U}_{j+1}$, $\widetilde{U}_j$, $U_{j-1}$$)$ is semistandard along $L$.
\vskip 1mm

	\item[(2)] 	Suppose that $\te_k {\bf T} \neq 0$ for some $k\in J$ and put
	${\bf S} = \te_k {\bf T} = (\dots, S_{i+1}, S_i,\dots)$.
	Then $T_{i+1}\lt T_i$ if and only if $S_{i+1} \lt S_i$.  
\end{itemize}
\end{lem}
\pf The proof is given in Section \ref{subsec:proof of sliding lemma}.
\qed



\subsection{Separation when $\la_n\geq 0$} \label{subsec:separation}
Let us assume $\lambda \in \mc{P}_n$ with $\lambda_n \ge 0$. 
The case when $\lambda_n < 0$ is considered in Section \ref{subsec:separation for odd}.

Let ${\bf T} = (T_{l}, \dots, T_1, T_0) \in {\bf T}_{\lambda}$ be given.
Since $T_i\in \mathbb{P}_L$ for $0\leq i\leq l$, we may consider 
the $(l+1)$-tuples 
\begin{equation}\label{eq:body and tail}
(T_{\ell}^{\texttt{body}}, \dots, T_{1}^{\texttt{body}}, \,T_0^{\texttt{body}}), \quad (T_{\ell}^{\texttt{tail}}, \dots, T_{1}^{\texttt{tail}}, \, T_0^{\texttt{tail}})
\end{equation} 
to form tableaux in $\mathbb{P}_L$.
But in general, they are not necessarily semistandard along $L$, and $(T_{l}^{\texttt{body}}, \dots, T_{1}^{\texttt{body}}, T_0^{\texttt{body}})$ may not be of a partition shape along $L$.
So instead of cutting ${\bf T}$ with respect to $L$ directly as in \eqref{eq:body and tail}, we will introduce an algorithm to separate ${\bf T}$ into two semistandard tableaux, which preserves ${\mf l}$-crystal equivalence. 

More precisely, we introduce an algorithm to get a semistandard tableau $\ov{\bf T}$ in $\mathbb{P}_L$ such that
\begin{itemize}
\item[(S1)] $\ov{\bf T}$ is Knuth equivalent to ${\bf T}$, that is, $\ov{\bf T}\equiv_{\mf l}{\bf T}$,

\item[(S2)] $\ov{\bf T}^{\texttt{tail}}\in SST_{[\ov{n}]}(\mu)$ and  $\ov{\bf T}^{\texttt{body}}\in SST_{[\ov{n}]}(\delta^\pi)$ for some $\delta\in \cP^{(1,1)}_n$,
where $\mu\in \cP_n$ is given by
\begin{equation}\label{eq:mu for positive case}
\mu' = (a_{\ell}, \dots, a_1)
\end{equation} with $a_i$ as in \eqref{eq:weight expression}.
\end{itemize}
We call this algorithm {\em separation} which can be viewed as a generalization of the one in \cite{JK19-2} (see \cite{K18-2} for type $B$ and $C$). Let us explain this with an example before we deal it in general.

\begin{ex} \label{ex:separation for positive case}
{\em 
Let ${\bf T}$ $=$ $($$T_4$, $T_3$, $T_2$, $T_1$$)$ $\in$ ${\bf T}_{(4,4,4,4,4,2)}$ be given in Example \ref{ex:spinor}.
\begin{equation*}
\begin{split}
\ytableausetup {mathmode, boxsize=1.0em} 
& \begin{ytableau}
\none & \none[\tl{\color{gray}{8}}] & \none[\tl{\color{gray}{7}}] & \none & \none[\tl{\color{gray}{6}}] & \none[\tl{\color{gray}{5}}] & \none & \none[\tl{\color{gray}{4}}] & \none[\tl{\color{gray}{3}}] & \none & \none[\tl{\color{gray}{2}}] & \none[\tl{\color{gray}{1}}]   \\
\none & \none & \none & \none & \none & \none & \none & \none & \none & \none & \none & \none   \\
\none &\none & \none & \none & \none & \none & \none & \none & \none & \none & \none & \tl{$\ov{6}$} & \none  \\
\none &\none & \none & \none & \none & \none & \none & \none & \none & \none & \none & \tl{$\ov{5}$} & \none  \\
\none &\none & \tl{$\ov{7}$} & \none & \none & \tl{$\ov{6}$} & \none & \none & \tl{$\ov{6}$} & \none & \tl{$\ov{5}$} & \tl{$\ov{3}$} & \none \\
\none[\!\!\!\!\mathrel{\raisebox{-0.7ex}{$\scalebox{0.45}{\dots\dots\dots\dots}$}}] &\none & \tl{$\ov{5}$} & \none[\mathrel{\raisebox{-0.7ex}{$\scalebox{0.45}{\dots\dots}$}}] & \none & \tl{$\ov{4}$} & \none[\mathrel{\raisebox{-0.7ex}{$\scalebox{0.45}{\dots\dots}$}}] & \none & \tl{$\ov{3}$} & \none[\mathrel{\raisebox{-0.7ex}{$\scalebox{0.45}{\dots\dots}$}}]& \tl{$\ov{4}$} & \tl{$\ov{2}$} & \none[\ \mathrel{\raisebox{-0.7ex}{\quad $\scalebox{0.45}{\dots\dots\dots\dots}$\ ${}_{\scalebox{0.75}{$L$}}$}}] \\
\none &\tl{$\ov{7}$} & \none & \none & \tl{$\ov{6}$} & \none & \none & \tl{$\ov{4}$} & \none & \none & \tl{$\ov{3}$} & \none & \none  \\
\none &\tl{$\ov{4}$} & \none & \none & \tl{$\ov{2}$} & \none & \none & \tl{$\ov{2}$} & \none & \none & \tl{$\ov{2}$} & \none & \none  \\
\none &\tl{$\ov{3}$} & \none & \none & \tl{$\ov{1}$} & \none & \none & \none & \none & \none & \none & \none & \none  \\
\none &\none & \none & \none & \none & \none & \none & \none & \none & \none & \none & \none & \none  \\
\none &\none[\tl{${U}_8$}] & \none[\,\,\,\,\tl{${U}_7$}] & \none & \none[\tl{${U}_6$}] & \none[\,\,\,\,\tl{${U}_5$}] & \none & \none[\tl{${U}_4$}] & \none[\,\,\,\,\tl{${U}_3$}] & \none & \none[\tl{${U}_2$}] & \none[\,\,\,\,\tl{${U}_1$}] & \none  \\
\end{ytableau} \quad\quad \\
\end{split}
\end{equation*} 
where $(U_1,\dots,U_8)$ denotes the image of $(T_4, T_3, T_2, T_1)$ under \eqref{eq:identification}. 

First we consider $\lt$ in Definition \ref{df:order} on $(T_{i+1}, T_i)$ for $1 \le i \le 3$. Then we can check that
$T_4 \not\triangleleft T_3$, $T_3 \not\triangleleft T_2$ and $T_2 \lt T_1$. 
By \eqref{eq:definition_S}, we have
\begin{equation*}
	\mathcal{S}_6 = \mathcal{E}_6\mathcal{E}_{5}\mathcal{F}_6^{2}\mathcal{F}_{5}, \quad 
	\mathcal{S}_4 = \mathcal{E}_4\mathcal{E}_{3}\mathcal{F}_4^{2}\mathcal{F}_{3}, \quad 
	\mathcal{S}_2 = \mathcal{F}_2^{2}.
\end{equation*}
Now, we apply these operators $\mc{S}_6$, $\mc{S}_4$, and then $\mc{S}_2$ to ${\bf T}$ to have
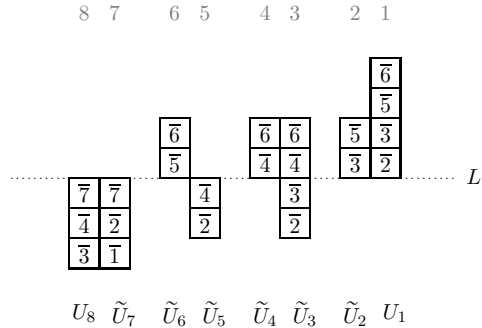
\begin{equation*}
\begin{split}
\ytableausetup {mathmode, boxsize=1.0em} 
& \begin{ytableau}
\none & \none[\tl{\color{gray}{8}}] & \none[\tl{\color{gray}{7}}] & \none & \none[\tl{\color{gray}{6}}] & \none[\tl{\color{gray}{5}}] & \none & \none[\tl{\color{gray}{4}}] & \none[\tl{\color{gray}{3}}] & \none & \none[\tl{\color{gray}{2}}] & \none[\tl{\color{gray}{1}}]   \\
\none & \none & \none & \none & \none & \none & \none & \none & \none & \none & \none & \none   \\
\none &\none & \none & \none & \none & \none & \none & \none & \none & \none & \none & \tl{$\ov{6}$} & \none  \\
\none &\none & \none & \none & \none & \none & \none & \none & \none & \none & \none & \tl{$\ov{5}$} & \none  \\
\none &\none & \none & \none & \tl{$\ov{6}$} & \none & \none & \tl{$\ov{6}$} & \tl{$\ov{6}$} & \none & \tl{$\ov{5}$} & \tl{$\ov{3}$} & \none \\
\none[\!\!\!\!\mathrel{\raisebox{-0.7ex}{$\scalebox{0.45}{\dots\dots\dots\dots}$}}] &\none & \none & \none[\mathrel{\raisebox{-0.7ex}{$\scalebox{0.45}{\dots\dots}$}}] & \tl{$\ov{5}$} & \none & \none[\mathrel{\raisebox{-0.7ex}{$\scalebox{0.45}{\dots\dots}$}}] & \tl{$\ov{4}$} & \tl{$\ov{4}$} & \none[\mathrel{\raisebox{-0.7ex}{$\scalebox{0.45}{\dots\dots}$}}]& \tl{$\ov{3}$} & \tl{$\ov{2}$} & \none[\ \mathrel{\raisebox{-0.7ex}{\quad $\scalebox{0.45}{\dots\dots\dots\dots}$\ ${}_{\scalebox{0.75}{$L$}}$}}] \\
\none &\tl{$\ov{7}$} & \tl{$\ov{7}$} & \none & \none & \tl{$\ov{4}$} & \none & \none & \tl{$\ov{3}$} & \none & \none & \none & \none  \\
\none &\tl{$\ov{4}$} & \tl{$\ov{2}$} & \none & \none & \tl{$\ov{2}$} & \none & \none & \tl{$\ov{2}$} & \none & \none & \none & \none  \\
\none &\tl{$\ov{3}$} & \tl{$\ov{1}$} & \none & \none & \none & \none & \none & \none & \none & \none & \none & \none  \\
\none &\none & \none & \none & \none & \none & \none & \none & \none & \none & \none & \none & \none  \\
\none &\none[\tl{${U}_8$}] & \none[\,\,\,\,\tl{$\widetilde{U}_7$}] & \none & \none[\tl{$\widetilde{U}_6$}] & \none[\,\,\,\,\tl{$\widetilde{U}_5$}] & \none & \none[\tl{$\widetilde{U}_4$}] & \none[\,\,\,\,\tl{$\widetilde{U}_3$}] & \none & \none[\tl{$\widetilde{U}_2$}] & \none[\,\,\,\,\tl{${U}_1$}] & \none  \\
\end{ytableau} \quad\quad \\
\end{split}
\end{equation*} 
We observe that (recall Definition \ref{def:admissibility})
\begin{equation}\label{eq:new spinor tableau}
	(\widetilde{U}_7, \widetilde{U}_6) \prec (\widetilde{U}_5, \widetilde{U}_4) \prec (\widetilde{U}_3, \widetilde{U}_2),
\end{equation}
(we will show in in Lemma \ref{lem:inductive step} that this holds in general).
So we can apply the above process to \eqref{eq:new spinor tableau}, and repeat it until there is no tail to move to the left horizontally.
Consequently we have
\begin{equation*}
\ytableausetup {mathmode, boxsize=1.0em} 
\begin{ytableau}
\none & \none & \none & \none & \none & \none & \none & \tl{$\ov{6}$} & \none  \\
\none & \none & \none & \none & \none & \none & \none & \tl{$\ov{5}$} & \none \\
\none & \none &\none & \tl{$\ov{6}$} & \tl{$\ov{6}$} & \tl{$\ov{6}$} & \tl{$\ov{5}$} & \tl{$\ov{3}$}  & \none\\
\none[\!\!\!\!\mathrel{\raisebox{-0.7ex}{$\scalebox{0.45}{\dots\dots\dots\dots}$}}] &\none & \none &\tl{$\ov{4}$} & \tl{$\ov{4}$} & \tl{$\ov{4}$} & \tl{$\ov{3}$} & \tl{$\ov{2}$}   & \none[\!\!\!\!\mathrel{\raisebox{-0.7ex}{$\scalebox{0.45}{\dots\dots\dots\dots}$\ ${}_{\scalebox{0.75}{$L$}}$}}] \\
\tl{$\ov{7}$} & \tl{$\ov{7}$} & \tl{$\ov{5}$} & \tl{$\ov{3}$} &\none & \none & \none & \none   & \none\\
\tl{$\ov{4}$} & \tl{$\ov{2}$} & \tl{$\ov{2}$} & \tl{$\ov{2}$} &\none & \none & \none & \none  & \none \\
\tl{$\ov{3}$} & \tl{$\ov{1}$} & \none & \none & \none & \none & \none & \none  & \none\\ 
\end{ytableau}
\end{equation*} 
\vskip 2mm
Hence we obtain two semistandard tableaux of shape $\delta^{\pi}$ and $\mu$, where $\delta = (4,4,1,1)$ and $\mu = (4,4,2)$.
}
\qed
\end{ex}

Now let ${\bf T} = (T_{l}, \dots, T_1, T_0) \in {\bf T}_{\lambda}$ be given and let ${\bf U}=(U_{2l},\dots,U_1, U_0)$ be its image under \eqref{eq:identification}.
We use the induction on the number of columns in ${\bf T}$ to define $\ov{\bf T}$. 
If $n\leq 3$, then let $\ov{\bf T}$ is given by putting together the columns in ${\bf U}$ horizontally along $L$. 

Suppose that $n\geq 4$. First, we consider
\begin{equation*} 
	\mathcal{S}_2 \dots \mathcal{S}_{2l-2}{\bf T} = \mathcal{S}_2 \dots \mathcal{S}_{2l-2} {\bf U} = (U_{2l}, \widetilde{U}_{2l-1}, \dots, \widetilde{U}_2, U_1, U_0) \in {\bf E}^{2l+1},
\end{equation*}
and let 
\begin{equation*}
\widetilde{\bf U} = (\widetilde{U}_{2l-1}, \dots, \widetilde{U}_2, U_1, U_0) \in {\bf E}^{2l}.
\end{equation*}
Note that applying $\mc{S}_{2i}$ to $\mathcal{S}_{2i+2} \dots \mathcal{S}_{2l-2}{\bf T}$ for $1\leq i\leq l-1$ is well-defined by Lemma \ref{lem:sliding}(1).

\begin{lem} \label{lem:inductive step}
Let $\td{\la}\in \mc{P}_n$ be such that $\omega_{\lambda} - \omega_{\widetilde{\lambda}} = \omega_{n-a}$ with $a = {\rm ht}(U_{2l}^{\texttt{\em tail}})$.
Then there exits a unique $\widetilde{\bf T} \in {\bf T}_{\widetilde{\lambda}}$ such that the image of $\widetilde{\bf T}$ under \eqref{eq:identification} is equal to $\widetilde{\bf U}$.
\end{lem}
\pf
There exists an ${\mf l}$-highest weight element ${\bf H} \in {\bf T}_{\lambda}$ such that
	${\bf H}	= \widetilde{e}_{i_1} \dots \widetilde{e}_{i_r} {\bf U}$
for some $i_1, \dots, i_r \in I$. 
Put ${\bf U}^{\sharp}=(U_{2l-1}, \dots, U_1)$ and let ${\bf H}^{\sharp}$ be obtained from ${\bf H}$ by removing its leftmost column, say $U'_{2l}$.
By tensor product rule of crystals, ${\bf H}^{\sharp}$ is also an ${\mf l}$-highest weight element.
Identifying ${\bf U}={\bf U}^\sharp\otimes U_{2l}$ We observe that
\vskip 2mm
\begin{equation*}
\begin{split}
	{\bf H}^{\sharp} \otimes U'_{2l}
	&  = \widetilde{e}_{i_1}\dots \widetilde{e}_{i_r} {\bf U} \\
	& = \widetilde{e}_{i_1}\dots \widetilde{e}_{i_r} \left( {\bf U}^{\sharp} \otimes U_{2l} \right) \\
	& =  \left( \widetilde{e}_{j_1} \dots \widetilde{e}_{j_s} {\bf U}^{\sharp} \right) \otimes \left( \widetilde{e}_{k_1} \dots \widetilde{e}_{k_t} U_{2l} \right)
\end{split}
\end{equation*}
where $\{ i_1, \dots, i_r \} = \{ j_1, \dots, j_s \} \cup \{ k_1, \dots, k_l \}$.
Hence
${\bf H}^{\sharp} = 	\widetilde{e}_{j_1} \dots \widetilde{e}_{j_s} {\bf U}^{\sharp}$,
and
\begin{equation*}
\widetilde{\bf U} = \mathcal{S}_2 \dots \mathcal{S}_{2l-2} {\bf U}^{\sharp}
= \mathcal{S}_2 \dots \mathcal{S}_{2l-2}\left( \widetilde{f}_{j_s} \dots \widetilde{f}_{j_1} {\bf H}^{\sharp} \right)
=  \widetilde{f}_{j_s} \dots \widetilde{f}_{j_1} \left( \mathcal{S}_2 \dots \mathcal{S}_{2l-2} {\bf H}^{\sharp} \right),
\end{equation*}
since ${\bf E}^{2l}$ is an $(\mf{l},\mf{sl}_{2l})$-bicrystal.
Note that the operator $\mathcal{S}_2 \dots \mathcal{S}_{2l-2}$ \eqref{eq:definition_S} is well-defined on ${\bf H}^{\sharp}$ by Remark \ref{rem:generalization} and Lemma  \ref{lem:sliding}(2).  Then $\mathcal{S}_2 \dots \mathcal{S}_{2l-2} {\bf H}^{\sharp}\in {\bf T}_{\td{\la}}$ by \cite[Lemma 3.10, 3.18]{JK19-2}, which implies that $\td{\bf U}\in {\bf T}_{\td{\la}}$.
\qed
\vskip 2mm

Put $\mathbb{T} = \widetilde{\bf T}$, which is given in Lemma \ref{lem:inductive step}. By induction hypothesis, there exists a tableau $\ov{\mathbb{T}}$ satisfying (S1) and (S2) associated to $\mathbb{T}$. Then we define $\ov{\bf T}$ to be the tableau in $\mathbb{P}_L$ obtained by putting together the leftmost column of ${\bf T}$, that is, $U_{2l}$, and $\ov{\mathbb T}$ along $L$. By definition, ${\rm sh}(\ov{\bf T})=\eta$ is of the following form:

\begin{equation}\label{eq:shape after separation}
\raisebox{.9ex}{$\eta=$} 
\resizebox{.38\hsize}{!}{
\def\lr#1{\multicolumn{1}{|@{\hspace{.75ex}}c@{\hspace{.75ex}}|}{\raisebox{-.04ex}{$#1$}}}
\def\l#1{\multicolumn{1}{|@{\hspace{.75ex}}c@{\hspace{.75ex}}}{\raisebox{-.04ex}{$#1$}}}
\def\r#1{\multicolumn{1}{@{\hspace{.75ex}}c@{\hspace{.75ex}}|}{\raisebox{-.04ex}{$#1$}}}
\raisebox{-.6ex}
{$\begin{array}{cccccccccc}
\cline{8-9}
& & & & & & &\l{\ \ } & \r{ } & \\
\cline{6-7}
 & & & & & \l{\ \ } & & & \r{ }\\
\cline{5-5} 
& & & & \l{\ \ } &  &  & & \r{\!\!\!\!\!\!\!\!\!\!\!\!\!\!\!\!\!\!\!\!\!\!\!\!\!\!{}_{\delta^\pi} } \\
\cline{3-4}
& & \l{\ \ } & & & & & & \r{ } &  \\
\cline{2-2}\cdashline{1-10}[0.5pt/1pt]\cline{6-9}
& \l{\ \ } &  & & \r{\ \ }  &  & &  \\ 
\cline{4-5}
& \l{\ \ } & \r{\!\!\!\!\!\!\!\!\! {}_{\mu}} & & &  & &  \\
& \l{\ \ } & \r{\ \ } & & & & &  \\  
\cline{3-3} 
& \lr{\ \ } & & & &  & &  \\  
\cline{2-2} 
\end{array}$}}\ \ \raisebox{-.7ex}{$L$}
\end{equation}
for some $\delta\in \cP_n^{(1,1)}$.

\begin{prop} \label{lem:separation}
Under the above hypothesis, $\ov{\bf T}$ satisfies (S1) and (S2).
\end{prop}
\pf By definition, it is clear that $\ov{\bf T}\equiv_{\mf l} \ov{T}$, which implies (S1). By Lemma \ref{lem:sliding},
$(U_{2l}^{\texttt{body}},\ov{\mathbb{T}}^{\texttt{body}})$ and $(U_{2l}^{\texttt{tail}},\ov{\mathbb{T}}^{\texttt{tail}})$ 
are semistandard along $L$, which implies (S2). 
\qed

\subsection{Separation when $\lambda_n < 0$} \label{subsec:separation for odd}
We consider the algorithm for separation when $\lambda_n < 0$ in a similar way as in \cite[Section 3.4]{JK19-2}.

Let us assume that $\lambda_n < 0$.
Recall that $-2\lambda_n = 2\ov{q}+\ov{r}$ with  $\ov{r} \in \{ \, 0, 1 \, \}$. 
For ${\bf T} \in {\bf T}_{\lambda}$, we may write 
\begin{equation*}
	{\bf T} = (T_l, \dots, T_{m+1}, T_m, \dots, T_1, T_0),
\end{equation*}
for some $m\ge 1$ such that $T_i \in {\bf T}(a_i)$ for some $a_i$ $(m+1 \le i \le l)$, $T_i \in \ov{\bf T}(0)$ $(1 \le i \le m)$ and $T_0 \in {\bf T}^{\rm sp-}$ (resp. $T_0 = \emptyset$) if $\ov{r} = 1$ (resp. $\ov{r} = 0$).
We identify ${\bf T}$ with
\begin{equation*}
	{\bf U}=(U_{2l}, \dots, U_{2m}, U_{2m-1}, \dots, U_1, U_0)
\end{equation*}
under \eqref{eq:identification}.
%
\begin{rem} \label{rem:sliding on spin column}
{\em 
For the spin column $U_{2m}$, we apply the sliding algorithm in Section \ref{subsec:generalization} as follows.
Let $U = H_{(1^n)}$ in \eqref{eq:H_(1^a)}.
Consider the pair
\begin{equation*}
	(U_{2m+2}, \,U_{2m+1}, \,U_{2m}, \,U),
\end{equation*}
where we regard $U = U \boxplus \emptyset\in \mathbb{P}_L$ and $(U_{2m}, U)\in {\bf T}(\varepsilon)$ as in Remark \ref{rem:admissibilty for spin column}(2).
By Lemma \ref{lem:sliding}, we have
\begin{equation*}
	\quad \quad (U_{2m+2}, \,U_{2m+1}, \,U_{2m}, \,U) \,
	\overset{\mathcal{S}_{2m}}{\xrightarrow{\hspace*{1.2cm}}} \,
	(U_{2m+2}, \,\widetilde{U}_{2m+1}, \,\widetilde{U}_{2m}, \,U),
\end{equation*}
for some $\widetilde{U}_{2m+1}$ and $\widetilde{U}_{2m}$, and by our choice of $U$, we have
\begin{equation*}
	(U_{2m+2}, \,U_{2m+1}, \,U_{2m}) \equiv_{\mf l}
	(U_{2m+2}, \,\widetilde{U}_{2m+1}, \,\widetilde{U}_{2m}).
\end{equation*}
}
\end{rem}

Now, we use the induction on the number of columns in ${\bf T}$ to define $\ov{\bf T}$ satisfying (S1) and (S2) where $\mu\in \cP_n$ in this is given by
\begin{equation} \label{eq:mu for negative case}
	\mu' = (a_{\ell}, \dots, a_1, \LaTeXunderbrace{1, \dots, 1}_{-2\lambda_n}).
\end{equation}
If $n\leq 3$, then let $\ov{\bf T}$ be given by putting together the columns in ${\bf U}$ along $L$. 

Suppose that $n\geq 4$. First, we consider $\mathcal{S}_{2l-2}\mathcal{S}_{2l-4} \dots \mathcal{S}_{2m} {\bf T}$, where $\mathcal{S}_{2m}$ is understood as in Remark \ref{rem:sliding on spin column}.. Then we have
\begin{equation} \label{eq:initial step of sliding for negative case}
	\mathcal{S}_{2l-2}\mathcal{S}_{2l-4} \dots \mathcal{S}_{2m} {\bf T} = (U_{2l}, \widetilde{U}_{2l-1}, \dots, \widetilde{U}_{2m+1}, \widetilde{U}_{2m}, U_{2m-1}, \dots, U_1, U_0),
\end{equation}
for some $\td{U}_i$ for $2m\le i\le 2l-1$.
Let 
$$\widetilde{{\bf U}} = (\widetilde{U}_{2l-1}, \dots, \widetilde{U}_{2m+1}, \widetilde{U}_{2m}, U_{2m-1}, \dots, U_1, U_0).$$
The following is an analogue of Lemma \ref{lem:inductive step} for $\lambda_n < 0$.

\begin{lem} \label{lem:inductive one step for negative case}
Let $\widetilde{\lambda}$ be such that $\omega_{\lambda} - \omega_{\widetilde{\lambda}} = \omega_{n-a}$ with $a = {\rm ht}(U_{2l}^{\texttt{\em tail}})$.
Then there exists a unique $\widetilde{\bf T} \in {\bf T}_{\widetilde{\lambda}}$ such that the image of $\widetilde{\bf T}$ is equal to $\widetilde{\bf U}$ under \eqref{eq:identification}.
\end{lem}
\pf 
Let ${\bf U}' = (U_{2l}, U_{2l-1}, \dots, U_{2m+1}, U_{2m}, \, U)$. 
Let $\nu = {\rm wt}({\bf U}')$ and let
$\widetilde{\nu}$ be given by $\nu - \widetilde{\nu} = \omega_{n-a}$ with $a={\rm ht}(U_{2l}^{\texttt{tail}})$.
By Lemma \ref{lem:inductive step} and Remark \ref{rem:sliding on spin column}, 
there exists a unique $\widetilde{\bf S} \in {\bf T}_{\widetilde{\nu}}$
whose image under  \eqref{eq:identification} is $\widetilde{\bf U'} = (\widetilde{U}_{2l-1}, \dots, \widetilde{U}_{2m}, U)$.
By semistandardness of $(U_{2m}, U_{2m-1})$ (cf. Remark \ref{rem:sliding on spin column}), we have
$$(\widetilde{U}_{2m+1}, \widetilde{U}_{2m}) \prec (U_{2m-1}, U_{2m-2}).
$$
Note that we may regard 
$(U_{2j-1}, U_{2j-2})\in \ov{{\bf T}}(0)$ $(2 \le j \le m)$, and $(U_1, U_0) \in \ov{\bf T}(0)$ if $U_0 \neq \emptyset$, $(U_1, U_0) \in {\bf T}^{\rm sp-}$ otherwise. Therefore there exists a unique
$\td{\bf T} \in {\bf T}_{\widetilde{\lambda}}$ which is equal to $\td{\bf U}$ under \eqref{eq:identification}.
\qed
\vskip 2mm

Put $\mathbb{T} = \widetilde{\bf T}$ given in Lemma \ref{lem:inductive one step for negative case}. By induction hypothesis, there exists a tableau $\ov{\mathbb{T}}$ satisfying (S1) and (S2) associated to $\mathbb{T}$. Then we define $\ov{\bf T}$ to be the tableau in $\mathbb{P}_L$ obtained by putting together the leftmost column of ${\bf T}$, that is, $U_{2l}$, and $\ov{\mathbb T}$ along $L$. By definition, ${\rm sh}(\ov{\bf T})=\eta$ is of the form \eqref{eq:shape after separation} with $\mu$ in \eqref{eq:mu for negative case}.

\begin{prop} \label{prop:separation for negative case}
Under the above hypothesis, $\ov{\bf T}$ satisfies (S1) and (S2).
\end{prop}
\pf It follows from the same argument as in the proof of Proposition \ref{lem:separation} with Lemma \ref{lem:inductive one step for negative case}.
\qed

\begin{ex} \label{ex:separation for negative case}
{\em 
Let $n=5$ and $\lambda = \left(\frac{5}{2}, \tfrac{3}{2}, \tfrac{3}{2}, \tfrac{1}{2}, -\tfrac{1}{2}\right)$. Then we have
\begin{equation*}
	\omega_{\lambda} = \omega_1+\omega_3+\omega_4, \quad 
	{\bf T}_{\lambda} \subset {\bf T}(4) \times {\bf T}(2) \times {\bf T}^{\rm sp-}.
\end{equation*}
Let us consider ${\bf T} = (T_2, T_1, T_0) \in {\bf T}_{\lambda}$ given by 
\begin{equation*}
\quad \quad \begin{split}
\ytableausetup {mathmode, boxsize=1.0em} 
& \begin{ytableau}
\none &\none & \none & \none & \none & \none & \none & \tl{$\ov{5}$} & \none & \none  \\
\none &\none & \tl{$\ov{5}$} & \none & \none & \tl{$\ov{3}$} & \none & \tl{$\ov{4}$} \\
\none[\!\!\!\!\mathrel{\raisebox{-0.7ex}{$\scalebox{0.45}{\dots\dots\dots\dots}$}}] &\none & \tl{$\ov{4}$} & \none[\mathrel{\raisebox{-0.7ex}{$\scalebox{0.45}{\dots\dots}$}}] & \none & \tl{$\ov{1}$} & \none[\mathrel{\raisebox{-0.7ex}{$\scalebox{0.45}{\dots\dots}$}}] & \tl{$\ov{1}$} &  \none[\ \mathrel{\raisebox{-0.7ex}{\quad $\scalebox{0.45}{\dots\dots\dots\dots}$\ ${}_{\scalebox{0.75}{$L$}}$}}] \\
\none &\tl{$\ov{5}$} & \none & \none & \tl{$\ov{2}$} & \none & \none & \none & \none & \none & \none  \\
\none &\tl{$\ov{3}$} & \none & \none & \tl{$\ov{1}$} & \none & \none & \none & \none & \none & \none  \\
\none &\tl{$\ov{2}$} & \none & \none & \none & \none & \none & \none & \none & \none & \none   \\
\none &\tl{$\ov{1}$} & \none & \none & \none & \none & \none & \none & \none & \none & \none   \\
\none & \none & \none & \none & \none & \none & \none & \none & \none & \none & \none   \\
\none & \none[\quad \, \tl{$T_2$}] & \none & \none & \none[\quad \, \tl{$T_1$}] & \none & \none & \none[\,\tl{$T_0$}] & \none & \none & \none   \\
\end{ytableau} \quad\quad \\
\end{split}
\end{equation*} 
We regard $T_0$ as in Remark \ref{rem:admissibilty for spin column}(1).
Then we have $T_2 \not\triangleleft T_1$ and $T_1 \not\triangleleft T_0$.
By applying \eqref{eq:initial step of sliding for negative case} (cf. Remark \ref{rem:sliding on spin column})
we obtain 
\vskip 1mm
\begin{equation*}
\quad \,\,\, \begin{split}
\ytableausetup {mathmode, boxsize=1.0em} 
& \begin{ytableau}
\none & \none & \none & \none & \none & \none & \none & \none & \none[\color{gray}{\tl{$\ov{5}$}}] & \none & \none & \none   \\
\none &\none & \none & \none & \none & \none & \none & \none & \none[\color{gray}{\tl{$\ov{4}$}}] & \none & \none   \\
\none &\none & \none & \none & \none & \none & \none & \none & \none[\color{gray}{\tl{$\ov{3}$}}] & \none  \\
\none &\none & \tl{$\ov{5}$} & \none & \none & \tl{$\ov{3}$} & \none & \tl{$\ov{5}$} & \none[\color{gray}{\tl{$\ov{2}$}}] \\
\none[\!\!\!\!\mathrel{\raisebox{-0.7ex}{$\scalebox{0.45}{\dots\dots\dots\dots}$}}] &\none & \tl{$\ov{4}$} & \none[\mathrel{\raisebox{-0.7ex}{$\scalebox{0.45}{\dots\dots}$}}] & \none & \tl{$\ov{1}$} & \none[\mathrel{\raisebox{-0.7ex}{$\scalebox{0.45}{\dots\dots}$}}] & \tl{$\ov{4}$} &  \none[\color{gray}{\tl{$\ov{1}$}} ]& \none[\ \mathrel{\raisebox{-0.7ex}{\quad $\scalebox{0.45}{\dots\dots\dots\dots}$}}] \\
\none &\tl{$\ov{5}$} & \none & \none & \tl{$\ov{2}$} & \none & \none & \tl{$\ov{1}$} & \none & \none & \none  \\
\none &\tl{$\ov{3}$} & \none & \none & \tl{$\ov{1}$} & \none & \none & \none & \none & \none & \none  \\
\none &\tl{$\ov{2}$} & \none & \none & \none & \none & \none & \none & \none & \none & \none   \\
\none &\tl{$\ov{1}$} & \none & \none & \none & \none & \none & \none & \none & \none & \none   \\
\none & \none & \none & \none & \none & \none & \none & \none & \none & \none & \none   \\
\none & \none[\quad \, \tl{$T_2$}] & \none & \none & \none[\quad \, \tl{$T_1$}] & \none & \none & \none[\quad \,\tl{$T_0$}] & \none & \none & \none   \\
\end{ytableau}
\end{split} \overset{\mathcal{S}_{2}\mathcal{S}_{0}}{\xrightarrow{\hspace*{1.2cm}}}
\begin{split}
\quad \quad \ytableausetup {mathmode, boxsize=1.0em} 
& \begin{ytableau}
\none & \none & \none & \none & \none & \none & \none & \none & \none[\color{gray}{\tl{$\ov{5}$}}] & \none & \none & \none   \\
\none &\none & \none & \none & \none & \none & \none & \tl{$\ov{5}$} & \none[\color{gray}{\tl{$\ov{4}$}}] & \none & \none   \\
\none &\none & \none & \none & \none & \none & \none & \tl{$\ov{4}$} & \none[\color{gray}{\tl{$\ov{3}$}}] & \none  \\
\none &\none & \none & \none & \tl{$\ov{5}$} & \none & \none & \tl{$\ov{3}$} & \none[\color{gray}{\tl{$\ov{2}$}}] \\
\none[\!\!\!\!\mathrel{\raisebox{-0.7ex}{$\scalebox{0.45}{\dots\dots\dots\dots}$}}] &\none & \none[\!\!\!\!\mathrel{\raisebox{-0.7ex}{$\scalebox{0.45}{\dots\dots\dots\dots}$}}] & \none[\mathrel{\raisebox{-0.7ex}{$\scalebox{0.45}{\dots\dots}$}}] & \tl{$\ov{2}$} & \none[\mathrel{\raisebox{-0.7ex}{$\scalebox{0.45}{\dots\dots}$}}] & \none[\mathrel{\raisebox{-0.7ex}{$\scalebox{0.45}{\dots\dots}$}}] & \tl{$\ov{1}$} &  \none[\color{gray}{\tl{$\ov{1}$}} ]& \none[\ \mathrel{\raisebox{-0.7ex}{\quad $\scalebox{0.45}{\dots\dots\dots\dots}$}}] \\
\none &\tl{$\ov{5}$} & \tl{$\ov{4}$} & \none & \none & \tl{$\ov{1}$} & \none & \none & \none & \none & \none  \\
\none &\tl{$\ov{3}$} & \tl{$\ov{1}$} & \none & \none & \none & \none & \none & \none & \none & \none  \\
\none &\tl{$\ov{2}$} & \none & \none & \none & \none & \none & \none & \none & \none & \none   \\
\none &\tl{$\ov{1}$} & \none & \none & \none & \none & \none & \none & \none & \none & \none   \\
\none & \none & \none & \none & \none & \none & \none & \none & \none & \none & \none   \\
\none & \none[\,\scalebox{0.7}{$U_4$}] & \none[\,\,\,\,\,\scalebox{0.7}{$\widetilde{U}_3$}] & \none & \none[\,\scalebox{0.7}{$\widetilde{U}_2$}] & \none[\,\,\,\,\,\scalebox{0.7}{$\widetilde{U}_1$}] & \none & \none[\,\,\scalebox{0.7}{$\widetilde{U}_0$}] & \none & \none & \none   \\
\end{ytableau}
\end{split}
\end{equation*}
where $U = H_{(1^5)}$ is the single-column tableau consisting of the numbers in gray. 

Finally, we apply $\mathcal{S}_1$ to $(U_4, \widetilde{U}_3, \widetilde{U}_2, \widetilde{U}_1, \widetilde{U}_0)$ and then we have $\ov{\bf T}$ given by
\begin{equation*}
\hspace{5cm} \ov{\bf T} = \hspace{-5cm}
\begin{split}
\quad \quad \ytableausetup {mathmode, boxsize=1.0em} 
& \begin{ytableau}
\none & \none &\none & \none & \none & \none & \none & \none & \none \\
\none & \none &\none & \none & \none & \none & \tl{$\ov{5}$} & \none & \none\\
\none & \none &\none & \none & \none & \none & \tl{$\ov{4}$} & \none & \none\\
\none & \none &\none & \none & \none & \tl{$\ov{5}$} & \tl{$\ov{3}$} & \none & \none\\
\none & \none[\mathrel{\raisebox{-0.7ex}{$\scalebox{0.45}{\dots\dots\dots\dots}$}}] &\none & \none & \none & \tl{$\ov{1}$} & \tl{$\ov{1}$} & \none[\mathrel{\raisebox{-0.7ex}{$\scalebox{0.45}{\dots\dots\dots\dots}$ \, ${}_{\scalebox{0.75}{$L$}}$}}] & \none\\
\none & \none &\tl{$\ov{5}$} & \tl{$\ov{4}$} & \tl{$\ov{2}$} & \none & \none & \none & \none\\
\none & \none &\tl{$\ov{3}$} & \tl{$\ov{1}$} & \none & \none & \none & \none & \none\\
\none & \none &\tl{$\ov{2}$} & \none & \none & \none & \none & \none & \none\\
\none & \none &\tl{$\ov{1}$} & \none & \none & \none & \none & \none & \none\\
\end{ytableau}
\end{split}
\end{equation*}
where $\ov{\bf T}^{\texttt{body}}$ (resp. $\ov{\bf T}^{\texttt{tail}}$) is the semistandard tableau located above $L$ (resp. below L) whose shape is $(2,2,1,1)^{\pi}$ (resp. $(3, 2, 1, 1)$). 
}	
\end{ex}

\subsection{Proof of Lemma \ref{lem:sliding}} \label{subsec:proof of sliding lemma}
Let ${\bf B}$ be one of ${\bf T}(a)$ $(0 \le a \le n-1)$ and ${\bf T}^{\rm sp}$ in \eqref{eq:def_T}.
When ${\bf B} = {\bf T}^{\rm sp}$,
we regard an element of ${\bf B}$ as in the sense of Remark \ref{rem:admissibilty for spin column}\,(2).
For $T \in {\bf B}$,
we define
\begin{equation} \label{eq:definition of s_T}
	\texttt{s}_T := \max_{1 \le s \le {\rm ht}(T^{\texttt{R}})}\,\left\{ s \,\, | \,\, T^{\texttt{L}}(s - \mf{r}_T - 1 + a) >
	T^{\texttt{R}}(s) \right\} \cup \{\,1\,\},
\end{equation}
where $T^{\texttt{L}}[k] := -\infty$ for $k \le 0$.
Note that $a=\mf{r}_T$ when ${\bf B} = {\bf T}^{\rm sp}$.

Let $(T, S) \in {\bf T}(a_2) \times {\bf B}$ be an admissible pair such that
$$T \in SST_{[\ov{n}]}(\lambda(a_2, b_2, c_2)),
\quad\,\, S \in SST_{[\ov{n}]}(\lambda(a_1, b_1, c_1))$$
for $a_i \in \mathbb{Z}_+$ and $b_i, c_i \in 2\mathbb{Z}_+$ $(i = 1, 2)$ with $a_1 \le a_2$.
If the pair $(T, S)$ satisfies 
${\rm ht}(T^{\texttt{R}}) > {\rm ht}(S^{\texttt{L}})-a_1\,$,
then put
$$S^{\texttt{L}}[0] := -\infty.$$
Note that above inequality occurs only for the case $\mf{r}_T \cdot \mf{r}_{S} = 1$.

\begin{lem} \label{lem1:appendix}
\mbox{}
\begin{itemize}
	\item[(1)] If $\mf{r}_T \cdot \mf{r}_S = 0$ \,$($resp. $\mf{r}_T \cdot \mf{r}_S = 1$$)$, then $T \lt S$ \,$\left(\textrm{resp.}\,\, T \lt (S^{\texttt{\em L}*}, \, S^{\texttt{\em R}*})\right)$.
	\vskip 1mm
	
	\item[(2)] If $\mf{r}_T \cdot \mf{r}_S = 1$, then $T \lt S$ is equivalent to the following condition:
\begin{equation} \label{eq:equivalent definition of the relation} 
	\quad T^{\texttt{\em R}}(k) \le S^{\texttt{\em L}}(k-1+a_1) \quad
	\textrm{for \,\,$\texttt{\em s}_T \le k \le {\rm ht}(T^{\texttt{\em R}})$.}
\end{equation}
\end{itemize}
\end{lem}
\pf
\vskip 1mm
(1) : If $\mf{r}_T \cdot \mf{r}_S = 0$ (resp. $\mf{r}_T \cdot \mf{r}_S = 1$), then
$T \lt S$ (resp. $T \lt (S^{\texttt{L}*}, \, S^{\texttt{R}*})$) follows immediately from Definition \ref{def:admissibility}\,(iii).
\vskip 2mm
 
(2) : Assume that $\mf{r}_T \cdot \mf{r}_S = 1$. 
	The relation $T \lt S$ implies
\begin{equation*}
	\quad \quad T^{\texttt{R}}(k) ={}^{\texttt{R}}T(k-1+a_2) \le S^{\texttt{L}}(k-1+a_1) \quad
	\, \textrm{for} \quad \texttt{s}_T \le k \le {\rm ht}(T^{\texttt{R}}).
\end{equation*}

Conversely, we assume that \eqref{eq:equivalent definition of the relation} holds.
Note that by definition of $S^{\texttt{L}*}$,
\begin{equation} \label{lem1_eq:Case 2}
	S^{\texttt{L}*}(k) = S^{\texttt{L}}(k),
\end{equation}
if there exists $k$ such that $1 \le k \le \texttt{s}_S-2+a_1$. 
Now we consider two cases.
\vskip 1mm

{\em Case 1}. $\texttt{s}_T > \texttt{s}_S$.
\,If $\texttt{s}_S \le k < \texttt{s}_T$, then by Definition \ref{def:admissibility}\,(ii),
\begin{equation} \label{lem1_eq:Case 1}
T^{\texttt{R}}(k) 
	= T^{\texttt{R}*}(k) 
	\le {}^{\texttt{L}}S(k) 
	= S^{\texttt{L}}(k-1+a_1).	
\end{equation}
Combining \eqref{eq:equivalent definition of the relation},
\eqref{lem1_eq:Case 2}, \eqref{lem1_eq:Case 1}
and Definition \ref{def:admissibility}\,(iii),
we have
\begin{equation*}
	{}^{\texttt{R}}T(k+a_2-a_1) \le
	S^{\texttt{L}}(k) \quad \textrm{for $1 \le k \le {\rm ht}(S^{\texttt{L}})$.}
\end{equation*}
\vskip 1mm

{\em Case 2}. $\texttt{s}_T \le \texttt{s}_S$.
\,In this case, the relation $T \lt S$ follows directly from \eqref{eq:equivalent definition of the relation}, \eqref{lem1_eq:Case 2} and Definition \ref{def:admissibility}\,(iii).
\qed
\vskip 2mm

\begin{lem} \label{lem2:appendix}
Assume that $\mf{r}_T \cdot \mf{r}_S=1$.
For $1 \le k \le {\rm ht}(T^{\texttt{R}})$,
\begin{itemize}
	\item[(1)] $
T^{\texttt{\em R}}(k) \leq S^{\texttt{\em L}*}(k).
$

	\item[(2)] If $T \lt S$, then
$
T^{\texttt{\em R}}(k) \leq S^{\texttt{\em L}}(k).
$
\end{itemize}
\end{lem}
\pf
By \eqref{lem1_eq:Case 1} and Definition \ref{def:admissibility}\,(iii),
in any case, we have
\begin{equation} \label{lem2_eq:Case 1}
	T^{\texttt{R}}(\texttt{s}_S) \le S^{\texttt{L}}(\texttt{s}_S-1+a_1) \le S^{\texttt{R}}(\texttt{s}_S)
	= S^{\texttt{L}*}(\texttt{s}_S-1+a_1).
\end{equation}
Then (1) follows from Definition \ref{def:admissibility}\,(ii)--(iii) and \eqref{lem2_eq:Case 1}.
By the same argument in the proof of Lemma \ref{lem1:appendix}\,(2), we obtain (2).
\qed
\vskip 3mm
Under the map \eqref{eq:identification-1}, put
\begin{equation} \label{eq:U_i} 
\begin{split}
	(T, S) &= (U_4,\, U_3,\, U_2,\, U_1), \\
	\mathcal{S}_2(T, S) &= (\widetilde{U}_4,\, \widetilde{U}_3,\, \widetilde{U}_2,\, \widetilde{U}_1).
\end{split}
\end{equation}
Here $\mathcal{S}_2$ is the operator given as in \eqref{eq:definition_S}.
For $1 \le i \le 4$, we regard a tableau $U_i$ as follow:
\begin{equation*}
\begin{split}
	U_1^{\texttt{body}} = U_1 \boxminus \emptyset, & \quad
	U_2^{\texttt{body}} = U_2 \boxminus S^{\texttt{tail}}, \\
	U_3^{\texttt{body}} = U_3 \boxminus \emptyset, & \quad 
	U_4^{\texttt{body}} = U_4 \boxminus T^{\texttt{tail}},
\end{split}
\end{equation*}
where $T^{\texttt{tail}} = (T^{\texttt{L}}(a_2)\,,\, \dots ,\, T^{\texttt{L}}(1))$ and $S^{\texttt{tail}} = (S^{\texttt{L}}(a_1)\,,\, \dots ,\, S^{\texttt{L}}(1))$.
Then we consider \eqref{eq:U_i} in $\mathbb{P}_L$.
For simplicity, put 
$$h = {\rm ht}(T^{\texttt{R}}), \quad g = {\rm ht}(S^{\texttt{R}}).$$
We consider some sequences given inductively as follows:
\begin{itemize}
	\item[(i)] Define a sequence $v_1 < \dots < v_h$ by 
	\begin{equation} \label{eq:seq v_k}
	\begin{split}
		v_1 & = \quad \,\,\,\min_{1 \le k \le a_2} \quad\,\,\,\,
		\left\{ k \,\left|\,\, T^{\texttt{L}}(k) \le T^{\texttt{R}}(1) \right. \right\}, \\
		v_s & = \min_{v_{s-1}+1\le k \le s+a_2}
		\left\{ k \,\left|\,\, T^{\texttt{L}}(k) \le T^{\texttt{R}}(s) \right. \right\}.
	\end{split}
	\end{equation}
	\vskip 2mm

	\item[(ii)] Define a sequence $w_1 < \dots < w_h$ by
	\begin{equation} \label{eq:seq w_k}
	\begin{split}
		w_h & =\, \max_{h \le k \le h+a_1} \,\,\,
		\left\{ k \, \left| \,\, T^{\texttt{R}}(h) \le S^{\texttt{L}}(k) \right. \right\}, \\
		w_t & = \max_{t \le k \le w_{t+1}-1}
		\left\{ k \, \left|\,\, T^{\texttt{R}}(t) \le S^{\texttt{L}}(k) \right. \right\}.
	\end{split}
	\end{equation}
	\vskip 2mm

	\item[(iii)] Define a sequence $x_1 < \dots < x_g$ by
	\begin{equation} \label{eq:seq x_k}
	\begin{split}
		x_1 & = \quad \,\,\,\min_{1 \le k \le a_1} \quad\,\,\,\,
		\left\{ k \,\left|\,\, S^{\texttt{L}}(k) \le S^{\texttt{R}}(1) \right. \right\}, \\
		x_u & = \min_{x_{u-1}+1\le k \le u+a_1}
		\left\{ k \,\left|\,\, S^{\texttt{L}}(k) \le S^{\texttt{R}}(u) \right. \right\}.
	\end{split}
	\end{equation}

\end{itemize}
\vskip 5mm

{\em Proof of Lemma \ref{lem:sliding}\,(1).}
By Lemma \ref{lem1:appendix}\,(1), the proof for the case $\mf{r}_T \cdot \mf{r}_S = 0$ is identical with the argument in \cite[Lemma 5.2]{K18-2}.
So we consider only the case $\mf{r}_T \cdot \mf{r}_S = 1$.
\vskip 2mm

{\em Case 1}. $T \lt S$.
In this case, $\mathcal{S}_2 = \mathcal{F}_2^{a_1}$.
This implies that $\widetilde{U}_1 = U_1$ and $\widetilde{U}_4=U_4$.
By Lemma \ref{lem2:appendix}\,(2),
\begin{equation*}
(T^{\texttt{R}},\, S^{\texttt{L}}) \in SST_{[\ov{n}]}(\lambda(0, b, c)),	
\end{equation*}
where $b = a_1 + c_1 - b_2 - c_2$ and $c = b_2+c_2$.
Since $T \prec S$ with $T \lt S$,\, $c_1 - b_2 - c_2 \ge 0$.
Therefore $\mathcal{F}_2^{a_1}(T, S)$ is well-defined.
Now we show that $(U_4, \widetilde{U}_3)$ and $(\widetilde{U}_2, U_1)$ are semistandard along $L$.
\vskip 3mm

{\em (1)} {\em $(U_4, \widetilde{U}_3)$ is semistandard along $L$.}
\,For $1 \le k \le {\rm ht}(T^{\texttt{R}})$ satisfying $v_k \ge a_2-a_1+1$, 
the relation $T \lt S$ implies
\begin{equation} \label{eq:relation v_k and w_k}
	v_k-a_2+a_1 \le w_k.
\end{equation}
\begin{itemize}
	\item[(i)] Assume that $\widetilde{U}_3(k) = T^{\texttt{R}}(k')$ for some $k'$. 
	By definition of $w_{k'}$, we have
	\begin{equation} \label{eq1:for case U_4 and U_3}
		k = w_{k'}.
	\end{equation}
	If $v_{k'} \ge a_2 - a_1 + 1$, then
	\begin{equation} \label{eq2:for case U_4 and U_3}
	\begin{split}
		U_4(k+a_2-a_1) &= T^{\texttt{L}}(w_{k'}+a_2-a_1) \quad \quad \textrm{by \eqref{eq:U_i} and \eqref{eq1:for case U_4 and U_3}} \\
		&\le T^{\texttt{L}}(v_{k'}) \quad \quad \quad \quad \quad \quad\,\,\,\, \textrm{by \eqref{eq:relation v_k and w_k}} \\
		&\le T^{\texttt{R}}(k') = \widetilde{U}_3(k) \quad \,\,\quad \,\,\,\,\,\, \textrm{by \eqref{eq:seq v_k}} \\
	\end{split}
	\end{equation}
	If $v_{k'} < a_2-a_1+1$, then
	\begin{equation} \label{eq3:for case U_4 and U_3}
		U_4(k+a_2-a_1) = T^{\texttt{L}}(k+a_2-a_1) < T^{\texttt{L}}(v_{k'}) \le T^{\texttt{R}}(k') = \widetilde{U}_3(k).
	\end{equation}
	\vskip 1mm
	
	\item[(ii)] Assume that $\widetilde{U}_3(k) \neq T^{\texttt{R}}(k')$ for any $k'$.
	In this case, we have
	\begin{equation*}
		\widetilde{U}_3(k) = S^{\texttt{L}}(k).
	\end{equation*}
	Then
	\begin{equation} \label{eq4:for case U_4 and U_3}
	\begin{split}
		U_4(k+a_2-a_1) &= T^{\texttt{L}}(k+a_2-a_1) \\
						&\le {}^{\texttt{R}}T(k+a_2-a_1) \quad \textrm{by definition of ${}^{\texttt{R}}T$} \\
						&\le S^{\texttt{L}}(k) = \widetilde{U}_3(k) \quad \,\,\,\, \textrm{by $T \lt S$}
	\end{split}
	\end{equation}
\end{itemize}
By \eqref{eq2:for case U_4 and U_3}, \eqref{eq3:for case U_4 and U_3} and \eqref{eq4:for case U_4 and U_3}, 
$(U_4, \widetilde{U}_3)$ is semistandard along $L$.
	\vskip 3mm

{\em (2)} {\em $(\widetilde{U}_2, U_1)$ is semistandard along $L$.}
By \eqref{eq:seq w_k}, \eqref{eq:seq x_k} and Definition \ref{def:admissibility}\,(ii),  we have
\begin{equation} \label{eq1:for case U_2 and U_1}
	x_k \le w_k \quad \textrm{for} \quad 1 \le k \le \texttt{s}_T-1.
\end{equation}
Also the relation $T \lt S$ implies
\begin{equation} \label{eq2:for case U_2 and U_1}
	k-1+a_1 \le w_k \quad \textrm{for} \quad k \ge \texttt{s}_T.
\end{equation}
\vskip 1mm

\begin{itemize}
	\item[(i)] For $1 \le k \le \texttt{s}_T-1$, 
	\begin{equation*}
	\begin{split}
		\widetilde{U}_2(k) &= S^{\texttt{L}}(w_k) \quad \quad \quad \quad \, \textrm{by \eqref{eq:seq w_k}}\\
						&\le S^{\texttt{L}}(x_k) \quad  \quad \quad \quad \,\, \textrm{by \eqref{eq1:for case U_2 and U_1}} \\
						&\le S^{\texttt{R}}(k) = U_1(k) \,\quad \textrm{by \eqref{eq:seq x_k}}
	\end{split}
	\end{equation*}	
	\vskip 2mm
	
	\item[(ii)] For $k \ge \texttt{s}_T$, \eqref{eq2:for case U_2 and U_1} implies
	\begin{equation*}
		\widetilde{U}_2(k) = S^{\texttt{L}}(w_k) \le S^{\texttt{L}}(k-1+a_1) \le S^{\texttt{R}}(k) = U_1(k).
	\end{equation*}
	Note that $S^{\texttt{L}}(k-1+a_1) \le S^{\texttt{R}}(k)$ holds since $\mf{r}_S = 1$.
\end{itemize}
\vskip 2mm
By (i)--(ii), $(\widetilde{U}_2, U_1)$ is semistandard along $L$.
\vskip 3mm

{\em Case 2}. $T \not\triangleleft S$.
In this case, $\mathcal{S}_2 = \mathcal{E}_2\mathcal{E}_1\mathcal{F}_2^{\,a_1-1}\mathcal{F}_1$.
Note that by definition
\begin{equation*}
	\mathcal{F}_1 (S^{\texttt{L}}, S^{\texttt{R}}) = (S^{\texttt{L}*}, S^{\texttt{R}*}).
\end{equation*}
By Lemma \ref{lem2:appendix}\,(1),
\begin{equation*}
(T^{\texttt{R}},\, S^{\texttt{L}*}) \in SST_{[\ov{n}]}(\lambda(0, b, c)),	
\end{equation*}
where $b = a_1+c_1-b_2-c_2+1$ and $c = b_2+c_2$.
Note that by Definition \ref{def:admissibility}\,(i),
\begin{equation*}
	b = a_1 - 1 + \left\{ c_1-(b_2+c_2-2) \right\} \ge a_1 - 1.
\end{equation*}
Therefore, $\mathcal{F}_2^{\,a_1-1}\mathcal{F}_1 (T, S)$ is well-defined.
Put
\begin{equation*}
	\mathcal{F}_2^{\,a_1-1}\mathcal{F}_1 (T, S) = (\dot{U}_4,\, \dot{U}_3,\, \dot{U}_2,\, \dot{U}_1).
\end{equation*}
Note that $\dot{U}_4 = U_4$ by definition of $\mathcal{F}_2^{\,a_1-1}\mathcal{F}_1$.
We use sequences $(w_k)$ and $(x_k)$ in \eqref{eq:seq w_k} and \eqref{eq:seq x_k} replacing 
$S^{\texttt{L}}$, $S^{\texttt{R}}$ and $a_1$ with $S^{\texttt{L}*}$, $S^{\texttt{R}*}$ and $a_1-1$, respectively.
\vskip 2mm

{\em (1)} {\em $(U_4, \dot{U}_3)$ is semistandard along $L$.}
By Lemma \ref{lem1:appendix}\,(1), we use the similar argument in the proof of {\em Case 1}\,{\em (1)}.
\vskip 3mm

{\em (2)} {\em $(\dot{U}_2, \dot{U}_1)$ is semistandard along $L$.}
We observe
\begin{equation*}
	\mathcal{E}^{a_1}(S^{\texttt{L}*}, S^{\texttt{R}*}) 
	= \mathcal{E}^{a_1}\,\mathcal{F}(S^{\texttt{L}}, S^{\texttt{R}})
	= \mathcal{E}^{a_1-1}(S^{\texttt{L}}, S^{\texttt{R}})
	= ({}^{\texttt{L}}S, {}^{\texttt{R}}S).
\end{equation*}
Therefore we use the similar argument in the proof of {\em Case 1\,(2)}.
\vskip 3mm

Note that \eqref{lem2_eq:Case 1}
implies that $S^{\texttt{R}}(\texttt{s}_S)$ is contained in $\dot{U}_2$.
Then the operator $\mathcal{E}_1$ on 
$(U_4,\, \dot{U}_3,\, \dot{U}_2,\, \dot{U}_1)$
moves $S^{\texttt{R}}(\texttt{s}_S)$ by one position to the right.
Therefore we have
\begin{equation} \label{eq1:for special case}
	\widetilde{U}_1 = U_1.
\end{equation}

On the other hand, \eqref{lem1_eq:Case 1} and Definition \ref{def:admissibility}\,(iii) implies that 
the operator $\mathcal{E}_2$ on $\mathcal{E}_1(U_4,\, \dot{U}_3,\, \dot{U}_2,\, \dot{U}_1)$ moves 
$$\quad \quad \dot{U}_3(w_k) = T^{\texttt{R}}(k) \quad \textrm{for some $k \ge \texttt{s}_T$.}$$ 
by one position to the right.
By the choice of $\texttt{s}_T$ and $\texttt{s}_S$ \eqref{eq:definition of s_T} with \eqref{eq1:for special case},
$(U_4, \widetilde{U}_3)$ and $(\widetilde{U}_2, U_1)$ are semistandard along $L$.

We complete the proof of Lemma \ref{lem:sliding}\,(1). 
\qed
\vskip 5mm

{\em Proof of Lemma \ref{lem:sliding}\,(2).}
Put 
$
	\widetilde{e}_k (T, \,S) = (T', S').
$
Then it is not difficult to check that
\begin{equation} \label{eq:Lemma 4.2(2)}
	\mf{r}_T = \mf{r}_{T'}, \quad \mf{r}_S = \mf{r}_{S'}.
\end{equation}

$(\Rightarrow)$\, 
Assume that
\begin{equation} \label{eq:Lemma 4.2(2) Case 1}
	\quad \quad  \widetilde{e}_k (T, \,S) = (T, \,\widetilde{e}_k\,S) \quad \textrm{and} \quad \widetilde{e}_k\,S = (\widetilde{e}_k\,S^{\texttt{L}}, S^{\texttt{R}}) \quad \,\, (k \in J).
\end{equation}
Otherwise it is clear that $T' \lt S'$.

If $\mf{r}_T \cdot \mf{r}_S = 0$, then Lemma \ref{lem1:appendix}\,(1) and \eqref{eq:Lemma 4.2(2)} implies that $T' \lt S'$ holds.

If $\mf{r}_T \cdot \mf{r}_S = 1$, then suppose $T' \not\triangleleft S'$. By Lemma \ref{lem1:appendix}\,(2),
there exists $s \ge \texttt{s}_T$ such that
\begin{equation*}
	T^{\texttt{R}}(s) = \ov{k}
\end{equation*}
which contradicts to \eqref{eq:Lemma 4.2(2) Case 1} by the tensor product rule.
Hence we have $T' \lt S'$.
\vskip 3mm

$(\Leftarrow)$\,
It follows from the similar argument of the previous proof.
\qed


\subsection{Embedding} \label{subsec:g-crystal equivalence}
For $\lambda \in \mc{P}_n$, let
\begin{equation*}
{\bf V}_{\lambda} :=
\Bigg( \raisebox{1mm}{$\underset{{\delta\in \cP_n^{(1,1)}}}{\scalebox{1.3}{$\bigsqcup$}} {SST}_{[\ov{n}]}\big(\delta^{\pi}\big)$} \Bigg) \raisebox{0.8mm}{$\times \, {SST}_{[\ov{n}]}(\mu)$}\,,
\end{equation*}
where 
$\mu\in\cP_n$ is given by \eqref{eq:mu for positive case} if $\lambda_n \ge 0$, and by \eqref{eq:mu for negative case} otherwise. 
Recall that
\begin{equation*}
{\bf V}:=\bigsqcup_{\delta\in \cP_n^{(1,1)}} {SST}_{[\ov{n}]}(\delta^{\pi}),
\end{equation*}
has a ${\mf g}$-crystal structure (see \cite[Section 5.2]{K13} for details).
On the other hand, we regard the ${\mf l}$-crystal 
\begin{equation*}
{\bf S}_\mu:={SST}_{[\ov{n}]}(\mu),
\end{equation*}
as a $\mf{g}$-crystal,  by defining $\widetilde{e}_n T = \widetilde{f}_n T = {\bf 0}$ with $\varphi_n(T) = \varepsilon_n(T) = -\infty$ for $T \in {SST}_{[\ov{n}]}(\mu)$.
Then we may regard ${\bf V}_\la$ as a $\mf g$-crystal by letting 
$${\bf V}_\la = {\bf V}\otimes {\bf S}_\mu,$$
which can be viewed as the crystal of a parabolic Verma module induced from a highest weight ${\mf l}$-module with highest weight $\la$.
The following is the main theorem in this section.
\begin{thm} \label{thm:separation}
The map   
\begin{equation} \label{eq:chi}
\xymatrixcolsep{3pc}\xymatrixrowsep{0pc}
	\xymatrix{
 {\bf T}_{\lambda} \ \ar@{->}[r] & {\bf V}_{\lambda} \otimes T_{r\omega_n} \\
{\bf T} \ar@{|->}[r] & \ov{\bf T}^{\texttt{\em body}}\otimes \ov{\bf T}^{\texttt{\em tail}} \otimes t_{r\omega_n}}
\end{equation}
is an embedding of crystals, where $r = ( \lambda, \omega_n)$.
\end{thm}

Before proving Theorem \ref{thm:separation}, let us recall the actions of $\te_n$ and $\tf_n$ on ${\bf T}_\la$ and ${\bf V}_\la$ in more details. We let
\begin{equation*}
{\texttt {vd}}
=\ytableausetup {mathmode, boxsize=1.5em} 
\begin{ytableau}
\tl{$\ov{n}$} \\ \tl{${}_{\ov{n-1}}$}
\end{ytableau}
\end{equation*}
be the vertical domino with entries $\ov{n}$ and $\ov{n-1}$.

Suppose that ${\bf T}=(T_l,\dots,T_1,T_0) \in {\bf T}_{\lambda}$ is given, where ${\bf T} = (U_{2\ell}, \dots, U_i, \dots, U_1, U_0)$ under \eqref{eq:identification}.
We define a sequence $\sigma = (\sigma_0,\sigma_1,\dots,\sigma_{2l})$ by
\begin{equation}\label{eq:signature}
\sigma_i = 
\left\{
		\begin{array}{ll}
 				+ \, &\textrm{if $U_i=\emptyset$ or $U_i[1]\geq \ov{n-2}$},\\
				- \, & \textrm{if $\texttt {vd}$ is located in the top of $U_i$}, \\
			\, \cdot \, \, \, & \textrm{otherwise},
			\end{array}
		\right.
\end{equation}
where we put $\sigma_0 = \cdot$ if $U_0 = \emptyset$.
Let $\sigma^{\rm red}$ be the (reduced) sequence obtained from $\sigma$ by replacing the pairs of neighboring signs $(+, -)$ (ignoring $\cdot$) with $(\,\cdot\,, \,\cdot\,)$ as far as possible.
If there exists a $-$ in $\sigma^{\rm red}$, then we define $\te_n{\bf T}$ to be the one obtained by removing $\texttt {vd}$
in $U_i$ corresponding to the rightmost $i$ such that $\sigma_i = -$ in $\sigma^{\rm red}$.
We define $\widetilde{e}_n {\bf T} = {\bf 0}$, otherwise.
Similarly, we define $\tf_n{\bf T}$ by adding $\texttt {vd}$ 
in $U_i$ corresponding to the leftmost $i$ such that $\sigma_i = +$ in $\sigma^{\rm red}$.

Next, suppose that $T\in {\bf V}_\la$ is given. 
We define a sequence $\tau = (\tau_0,\tau_1,\dots)$ by \eqref{eq:signature}. Note that $\tau$ is an infinite sequence where $\tau_i=+$ for all sufficiently large $i$. Then we define $\tau^{\rm red}$ and hence $\te_nT$ and $\tf_nT$ in the same way as in ${\bf T}_\la$.

\begin{lem} \label{lem:invariant nth signature}
Under the above hypothesis, let $\ov{\sigma}=(\ov{\sigma}_0,\dots, \ov{\sigma}_{2l})$ be the subsequence of $\tau^{\rm red}$ consisting of its first $2l+1$ components. If we ignore $\cdot$\,, then
$$
\sigma^{\rm red} = \ov{\sigma}^{\rm red}.
$$
\end{lem}
\pf 
We first consider a pair $(T_{i+1} ,T_i)$ in ${\bf T}$. Let $($$\sigma_{j-1}$, $\sigma_{j}$, $\sigma_{j+1}$, $\sigma_{j+2}$$)$ be the subsequence of $\sigma$ corresponding to $(T_{i+1}, T_i)$ with $j=2i$.
Let 
$$\mathcal{S}_j {\bf T} = (\dots, U_{j+2}, \widetilde{U}_{j+1}, \widetilde{U}_j, U_{j-1}, \dots)$$ 
(see Lemma \ref{lem:sliding}).
We denote by $(\sigma_{j-1}, \widetilde{\sigma}_{j}, \widetilde{\sigma}_{j+1}, \sigma_{j+2})$ the sequence defined by \eqref{eq:signature} corresponding to $($$U_{j+2}$, $\widetilde{U}_{j+1}$, $\widetilde{U}_j$, $U_{j-1}$$)$.

\vskip 2mm

\noindent {\em Case 1.} $\mf{r}_{i+1}\,\mf{r}_i = 0$.
Note that $\mathcal{S}_j = \mathcal{F}_j^{\, a_i}$ by Lemma \ref{lem1:appendix}\,(1).
Then we have 
$(\sigma_{j+1}, \sigma_j)=(\widetilde{\sigma}_{j+1}, \widetilde{\sigma}_j)$
by the similar argument in the proof of \cite[Theorem 5.7]{K18-2}.

%
\vskip 2mm

\noindent {\em Case 2.} $\mf{r}_{i+1} \, \mf{r}_i = 1$.
The relation between the two pairs $(\sigma_{j+1}, \sigma_j)$ and $(\widetilde{\sigma}_{j+1}, \widetilde{\sigma}_j)$ is given in Table \ref{table:2}. 
\begin{table}[h!]
\centering
\begin{tabular}{|c|c|} 
\hline
$(\sigma_{j}, \sigma_{j+1})$ & $(\widetilde{\sigma}_{j}, \widetilde{\sigma}_{j+1})$  \\ 
\hline
\hline
\quad \quad $(+, +)$ \quad \quad  & $(+, +)$  \\ 
\hline
\quad \quad $(+, -)$ \quad \quad & $(+, -)$ or $(\, \cdot \, , \, \cdot \, )$  \\ 
\hline
\quad \quad $(+, \, \cdot \,)$ \quad \quad & $(+, \, \cdot \,)$ or $(\, \cdot \, , +)$ \\ 
\hline
\quad \quad $(-,+)$ \quad \quad & $(-,+)$  \\ 
\hline
\quad \quad $(-,-)$ \quad \quad & $(-,-)$  \\ 
\hline
\quad \quad $(-,\, \cdot \,)$ \quad \quad & $(-,\, \cdot \,)$ or $(\, \cdot \, , -)$  \\ 
\hline
\quad \quad $(\, \cdot \, , +)$ \quad \quad & $(\, \cdot \, , +)$ or $( +, \, \cdot \,)$  \\ 
\hline
\quad \quad $(\, \cdot \, , -)$ \quad \quad & $(\, \cdot \, , -)$ or $(-, \, \cdot\,)$  \\ 
\hline
\quad \quad $(\, \cdot \, , \, \cdot \, )$ \quad \quad & $(\, \cdot \, , \, \cdot \, )$  \\
\hline
\end{tabular}
\vskip 3mm
\caption{The relation between $(\sigma_{j}, \sigma_{j+1})$ and $(\widetilde{\sigma}_{j}, \widetilde{\sigma}_{j+1})$ when $\mf{r}_{i+1} \mf{r}_i = 1$}
\label{table:2}
\end{table}
\vskip -7mm
Let ${\bf U} := \mathcal{S}_2 \dots \mathcal{S}_{2l-2}{\bf T} := (U_{2l}, \widetilde{U}_{2l-1}, \dots, \widetilde{U}_2, U_1, U_0)$ and let $\dot{\sigma}$ be the sequence given by \eqref{eq:signature} corresponding to ${\bf U}$.
Then we have
$\dot{\sigma} = (\sigma_{2l}, \widetilde{\sigma}_{2l-1}, \dots, \widetilde{\sigma}_2, \sigma_1, \sigma_0).$
It is straightforward to check that $\sigma^{\rm red} = \dot{\sigma}^{\rm red}$.
Note that $\sigma_{2l} = \ov{\sigma}_{2l}$ by definition. 
By Lemma \ref{lem:inductive step} and \ref{lem:inductive one step for negative case}, we may use an inductive argument for $(\widetilde{U}_{2\ell-1}, \dots, \widetilde{U}_2, U_1, U_0)$ to have
$\sigma^{\rm red} = \dot{\sigma}^{\rm red} = \ov{\sigma}^{\rm red}$.
This completes the proof.
%
%
\qed
\vskip 2mm

\noindent {\em Proof of Theorem \ref{thm:separation}}.
We will use the following notations under the identification \eqref{eq:identification} in this proof.
\vskip 2mm
\begin{itemize}
	\item[(1)] $\widetilde{\bf T}= (\widetilde{U}_{2l}, \widetilde{U}_{2l-1}, \dots)$ : the sequence of tableaux with a column shape obtained from ${\bf T}$
by shifting the tails one position to the left as in Sections \ref{subsec:separation} and \ref{subsec:separation for odd}. 
\vskip 1mm 

	\item[(2)] $\mathbb{T} = (\widetilde{U}_{2l-1}, \widetilde{U}_{2l-2}, \dots)$ : the sequence of tableaux obtained from $\widetilde{\bf T}$ by removing the column $\widetilde{U}_{2l}$.
\vskip 1mm

	\item[(3)] $\texttt{T} = (\texttt{U}_{2l}, \texttt{U}_{2l-1}, \dots )$ : the sequence obtained from ${\bf T}$ by applying $\widetilde{f}_n$, that is, 
	$$\texttt{T} = \widetilde{f}_n\,{\bf T}.$$

	\item[(4)] $\td{\texttt{T}}$ : the sequence obtained from $\texttt{T}$ by shifting the tails one position to the left as in Sections \ref{subsec:separation} and \ref{subsec:separation for odd} and then removing the left-most column of $\texttt{T}$. 
\vskip 1mm

	\item[(5)] $\widetilde{\sigma} = (\widetilde{\sigma}_0, \widetilde{\sigma}_1, \dots,\widetilde{\sigma}_{2l})$ : the sequence of \eqref{eq:signature} associated to $\widetilde{\bf T}$.
\end{itemize}
\vskip 1mm
Recall that $\mathbb{T} \in {\bf T}_{\widetilde{\lambda}}$ by Lemmas \ref{lem:inductive step} and \ref{lem:inductive one step for negative case}, where $\widetilde{\lambda}$ satisfies $\omega_{\lambda} - \omega_{\widetilde{\lambda}} = \omega_{n-a}$ with $a = {\rm ht}(U_{2l}^{\texttt{tail}})$.
\vskip 2mm

Let us denote by $\chi_\la$ the map in \eqref{eq:chi}.
First we show that $\chi_\la$ is injective.
Suppose that $\chi_{\lambda}({\bf T}) = \chi_{\lambda}({\bf T}')$ for some ${\bf T}, {\bf T}' \in {\bf T}_{\lambda}$.
There exists $i_1, \dots, i_r \in J$ such that $\widetilde{e}_{i_1} \dots \widetilde{e}_{i_r} \chi_{\lambda}({\bf T}) = \widetilde{e}_{i_1} \dots \widetilde{e}_{i_r} \chi_{\lambda}({\bf T}')$ is an ${\mf l}$-highest weight element. 
Note that $\chi_{\lambda}$ is a morphism of ${\mf l}$-crystals by Proposition \ref{lem:separation} (cf. Section \ref{subsec:g-crystal equivalence}).
Therefore, we have
	$\chi_{\lambda}(\widetilde{e}_{i_1} \dots \widetilde{e}_{i_r}{\bf T}) = \chi_{\lambda}(\widetilde{e}_{i_1} \dots \widetilde{e}_{i_r} {\bf T}')$.
By \cite[Proposition 3.4, Lemma 6.5]{JK19-2}, we have
$$\widetilde{e}_{i_1} \dots \widetilde{e}_{i_r}{\bf T} = \widetilde{e}_{i_1} \dots \widetilde{e}_{i_r} {\bf T}'.$$ 
Then ${\bf T} = {\bf T}'$.
Hence $\chi_{\lambda}$ is injective.
\vskip 2mm

Now it remains to show that
\begin{equation} \label{appendix B claim}
	\textrm{$\widetilde{f}_n {\bf T} \neq {\bf 0}$ and $\chi_{\lambda}({\bf T}) \neq {\bf 0}$ for ${\bf T} \in {\bf T}_{\lambda}$ \quad $\Longrightarrow$ \quad $\chi_{\lambda}(\widetilde{f}_n{\bf T}) = \widetilde{f}_n\chi_{\lambda}({\bf T})$}.
\end{equation}
If $\widetilde{f}_n$ acts on $U_{2l}$ or $U_{k}$\,, where\,
$k = 0$
if $U_0 \neq \emptyset$, and 
$k= 1$
if $U_0 = \emptyset$,
then the claim \eqref{appendix B claim} follows from Lemmas \ref{lem:sliding} and \ref{lem:invariant nth signature}. 

Suppose that $\widetilde{f}_n$ acts on $U_k$ for $k < 2l$.
%
%
%
%
%
To prove \eqref{appendix B claim}, it is enough to show that 
\begin{equation} \label{appendix B claim 2}
	\widetilde{\texttt{T}} = 
\widetilde{f}_n \mathbb{T}.
\end{equation}
Indeed, if \eqref{appendix B claim 2} holds, then by induction on the number of columns in ${\bf T}$, we have
\begin{equation*} \label{appendix B induction argument}
\begin{split}
	\ov{\widetilde{f}_n {\bf T}} 
	& =	\ov{\left(U_{2l},\,\, \widetilde{\texttt{T}} \right)} \\ 
	& = \ov{\left(U_{2l},\,\,\widetilde{f}_n \mathbb{T}\right)} \quad \quad\, \textrm{by \eqref{appendix B claim 2}} \\
	& = \left(U_{2l},\,\,\ov{\widetilde{f}_n \mathbb{T}}\right) \\
	& = \left(U_{2l},\,\,\widetilde{f}_n\ov{\mathbb{T}}\right) \quad \quad \textrm{by induction hypothesis}\\
	& = \widetilde{f}_n\left(U_{2l},\,\, \ov{\mathbb{T}}\right) \quad \quad \, \textrm{by Lemma \ref{lem:invariant nth signature}}, \\
\end{split}
\end{equation*}
Hence we have $\ov{\widetilde{f}_n {\bf T}} = \widetilde{f}_n \ov{\bf T}$ which implies \eqref{appendix B claim}.
%

Now we verify \eqref{appendix B claim 2}. 
Let us recall Table \ref{table:2} and then we have
\begin{equation} \label{eq:one step red. sig.}
	\sigma^{\rm red} = \widetilde{\sigma}^{\rm red}.
\end{equation}

We will prove \eqref{appendix B claim 2} for the case $U_k^{\texttt{tail}} \neq \emptyset$ with non-trivial sub-cases. 
The proof of other cases or the case $U_k^{\texttt{tail}} = \emptyset$ is almost identical, so we leave it to the reader to verify \eqref{appendix B claim 2} for these cases.

Let us recall the action of $\widetilde{f}_n$ in Section \ref{subsec:g-crystal equivalence}.
Then for the case $U_k^{\texttt{tail}} \neq \emptyset$, 
the signature $\sigma_{k+1}$ \eqref{eq:signature} associated to $U_{k+1}$ must be $\scalebox{1.2}{$\cdot$}$ or $+$\,.
\vskip 2mm

\begin{itemize}
	\item[{\em Case 1.}] If ${\rm ht}(U_{k+1}^{\texttt{body}}) \le {\rm ht}(U_k^{\texttt{body}}) < {\rm ht}(U_{k-1}^{\texttt{body}})={\rm ht}(U_{k}^{\texttt{body}})+2$,
then by definition of $\mathcal{S}_k$ the top entry of $U_{k+1}$ can not be moved to the right in $\widetilde{\bf T}$.
Thus we obtain
$(\sigma_{k}, \sigma_{k+1}) = (\widetilde{\sigma}_{k}, \widetilde{\sigma}_{k+1})$, which implies \eqref{appendix B claim 2}.
	\vskip 2mm

	\item[{\em Case 2.}]	If ${\rm ht}(U_{k+1}^{\texttt{body}}) = {\rm ht}(U_k^{\texttt{body}})+2 = {\rm ht}(U_{k-1}^{\texttt{body}})$, then the signature $\sigma_{k+1}$ \eqref{eq:signature} can not be $+$. Otherwise
	$(\texttt{U}_{k+2}, \texttt{U}_{k+1}) \not\prec (\texttt{U}_k, \texttt{U}_{k-1})$ for Definition \ref{def:admissibility}\,(ii), which is a contradiction. 
	Thus $\sigma_{k+1} = \scalebox{1.2}{$\cdot$}$\,.
		We observe that 
		\begin{equation} \label{eq:sign change in Case 1 (ii)}
		\begin{split}
			& (\sigma_{k}, \sigma_{k+1}) = (+,\,\scalebox{1.2}{$\cdot$}\,) \longrightarrow (\widetilde{\sigma}_{k}, \widetilde{\sigma}_{k+1}) = (\,\scalebox{1.2}{$\cdot$}\,, +), \\
			& U_{k+1}[1] =  
				\raisebox{-.6ex}{{\tiny ${\def\lr#1{\multicolumn{1}{|@{\hspace{.8ex}}c@{\hspace{.8ex}}|}{\raisebox{-.3ex}{$#1$}}}\raisebox{-.0ex}
{$\begin{array}[b]{c}
\cline{1-1}
\lr{ \!\,\ov{n}\,\! }\\
\cline{1-1}
\end{array}$}}$}}
				\,\, \textrm{or} \,\,
				\raisebox{-.6ex}{{\tiny ${\def\lr#1{\multicolumn{1}{|@{\hspace{.8ex}}c@{\hspace{.8ex}}|}{\raisebox{-.4ex}{$#1$}}}\raisebox{-.0ex}
{$\begin{array}[b]{c}
\cline{1-1}
\lr{ \!\,\ov{n-1}\,\! }\\
\cline{1-1}
\end{array}$}}$}}\,,
				\quad
		U_{k+1}[2] \ge \ov{n-2}.
		\end{split}
		\end{equation}
	It is straightforward to check from the definition of $\mathcal{S}_k$ that when we apply $\mathcal{S}_k$ to $\texttt{T}$, the domino
	{\texttt {vd}}
	in $\texttt{U}_k$ is changed as follows:
	\begin{equation} \label{eq:domino change in Case 1 (ii)}
		\begin{cases}
			\textrm{
				\raisebox{-.6ex}{{\tiny ${\def\lr#1{\multicolumn{1}{|@{\hspace{.8ex}}c@{\hspace{.8ex}}|}{\raisebox{-.3ex}{$#1$}}}\raisebox{-.0ex}
{$\begin{array}[b]{c}
\cline{1-1}
\lr{ \!\,\ov{n-1}\,\! }\\
\cline{1-1}
\end{array}$}}$}}
			in $\texttt{U}_k$ is moved to
				\raisebox{-.6ex}{{\tiny ${\def\lr#1{\multicolumn{1}{|@{\hspace{.8ex}}c@{\hspace{.8ex}}|}{\raisebox{-.3ex}{$#1$}}}\raisebox{-.0ex}
{$\begin{array}[b]{c}
\cline{1-1}
\lr{ \!\,\ov{n}\,\! }\\
\cline{1-1}
\end{array}$}}$}} in $\texttt{U}_{k+1}$ below
				} \, & \textrm{if} \,\, \texttt{U}_{k+1}[1] = 
				\raisebox{-.6ex}{{\tiny ${\def\lr#1{\multicolumn{1}{|@{\hspace{.8ex}}c@{\hspace{.8ex}}|}{\raisebox{-.3ex}{$#1$}}}\raisebox{-.0ex}
{$\begin{array}[b]{c}
\cline{1-1}
\lr{ \!\,\ov{n}\,\! }\\
\cline{1-1}
\end{array}$}}$}} \,,\\
			\textrm{
				\raisebox{-.6ex}{{\tiny ${\def\lr#1{\multicolumn{1}{|@{\hspace{.8ex}}c@{\hspace{.8ex}}|}{\raisebox{-.3ex}{$#1$}}}\raisebox{-.0ex}
{$\begin{array}[b]{c}
\cline{1-1}
\lr{ \!\,\ov{n}\,\! }\\
\cline{1-1}
\end{array}$}}$}}
			in $\texttt{U}_k$ is moved to
				\raisebox{-.6ex}{{\tiny ${\def\lr#1{\multicolumn{1}{|@{\hspace{.8ex}}c@{\hspace{.8ex}}|}{\raisebox{-.3ex}{$#1$}}}\raisebox{-.0ex}
{$\begin{array}[b]{c}
\cline{1-1}
\lr{ \!\,\ov{n-1}\,\! }\\
\cline{1-1}
\end{array}$}}$}} in $\texttt{U}_{k+1}$ above
				} \, & \textrm{if} \,\, \texttt{U}_{k+1}[1] = 
				\raisebox{-.6ex}{{\tiny ${\def\lr#1{\multicolumn{1}{|@{\hspace{.8ex}}c@{\hspace{.8ex}}|}{\raisebox{-.3ex}{$#1$}}}\raisebox{-.0ex}
{$\begin{array}[b]{c}
\cline{1-1}
\lr{ \!\,\ov{n-1}\,\! }\\
\cline{1-1}
\end{array}$}}$}}\,.
		\end{cases}
	\end{equation}
	Combining \eqref{eq:one step red. sig.}, \eqref{eq:sign change in Case 1 (ii)}
 and \eqref{eq:domino change in Case 1 (ii)}, we conclude that \eqref{appendix B claim 2} holds in this case.
	\vskip 2mm
	
	\item[{\em Case 3.}] If ${\rm ht}(U_{k+1}^{\texttt{body}}) \le {\rm ht}(U_k^{\texttt{body}})$ and ${\rm ht}(U_k^{\texttt{body}}) + 2 < {\rm ht}(U_{k-1}^{\texttt{body}})$, then
	we have
	\begin{equation*}
		(U_{k+2}, U_{k+1}) \lt (U_k, U_{k-1}) \longrightarrow 
		(\texttt{U}_{k+2}, \texttt{U}_{k+1}) \lt (\texttt{U}_k, \texttt{U}_{k-1}).
	\end{equation*}
	By \eqref{eq:one step red. sig.}, this implies that \eqref{appendix B claim 2} holds in this case.
	\vskip 2mm
	
	\item[{\em Case 4.}] Suppose ${\rm ht}(U_{k+1}^{\texttt{body}}) = {\rm ht}(U_k^{\texttt{body}})+2$ and ${\rm ht}(U_k^{\texttt{body}}) + 2 < {\rm ht}(U_{k-1}^{\texttt{body}})$.
	First we consider the case $\sigma_{k+1} = \scalebox{1.2}{$\cdot$}$\,.
	If there exists $i$ such that $U_{k+1}(i) > U_k(i+a_k-1)$, where $a_k = {\rm ht}(U_k^{\texttt{ail}})$, then by definition
	\begin{equation*}
		(U_{k+2}, U_{k+1}) \not\triangleleft (U_k, U_{k-1}) \longrightarrow 
		(\texttt{U}_{k+2}, \texttt{U}_{k+1}) \not\triangleleft (\texttt{U}_k, \texttt{U}_{k-1}).
	\end{equation*}
	Form this we observe that the domino
	{\texttt {vd}} in $\texttt{U}_k$
	is not moved when we apply the {\em sliding} to $\texttt{T}$.
	Also we note that
	\begin{equation*}
		(\sigma_{k}, \sigma_{k+1}) = (+,\,\scalebox{1.2}{$\cdot$}\,) \longrightarrow (\widetilde{\sigma}_{k}, \widetilde{\sigma}_{k+1}) = (+, \,\scalebox{1.2}{$\cdot$}\,).
	\end{equation*}
	Combining these observations with \eqref{eq:one step red. sig.}, we obtain \eqref{appendix B claim 2} in this case.
	If there is no such $i$, we have
	\begin{equation*}
	\begin{split}
	& (U_{k+2}, U_{k+1}) \not\triangleleft (U_k, U_{k-1}) \longrightarrow 
		(\texttt{U}_{k+2}, \texttt{U}_{k+1}) \lt (\texttt{U}_k, \texttt{U}_{k-1}), \\
	& \quad \quad \quad (\sigma_{k}, \sigma_{k+1}) = (+,\,\scalebox{1.2}{$\cdot$}\,) \longrightarrow (\widetilde{\sigma}_{k}, \widetilde{\sigma}_{k+1}) = (\,\scalebox{1.2}{$\cdot$}\,, +).	
	\end{split}
	\end{equation*}
	Note that $\widetilde{U}_k[1] = U_{k+1}[1] \le \ov{n-1}$
	and we observe that \eqref{eq:domino change in Case 1 (ii)} also holds in this case.
	Hence we have \eqref{appendix B claim 2}.
	
	Second we consider the case $\sigma_{k+1}=+$.
	In this case, we obtain 
	\begin{equation*}
	(U_{k+2}, U_{k+1}) \not\triangleleft (U_k, U_{k-1}) \longrightarrow 
		(\texttt{U}_{k+2}, \texttt{U}_{k+1}) \not\triangleleft (\texttt{U}_k, \texttt{U}_{k-1}),
	\end{equation*}
	since $U_{k+1}[1] > \ov{n-1}$.
	Also it is clear that 
	$(\widetilde{\sigma}_{k}, \widetilde{\sigma}_{k+1}) = (+, +)$ by definition of $\mathcal{S}_k$.
	On the other hand,
	the domino
	{\texttt {vd}}
	in $\texttt{U}_k$ can not be moved to the left when we apply the {\em sliding} to $\texttt{T}$.
	Consequently, we have \eqref{appendix B claim 2}.
\end{itemize}
\vskip 2mm
We complete the proof of Theorem \ref{thm:separation}.
\qed


\section{Embedding into the crystal of Lusztig data}\label{sec:parabolic into PBW}
In this section, we describe the embedding of ${\bf V}_\la$ into the crystal of Lusztig data. 

\subsection{PBW crystal}
Let $N=n^2-n$, which is the length of the longest element $w_0\in W$ of type $D_n$. 
Let $\bi_0=(i_1,\dots,i_N)$ be the reduced expression of $w_0\in W$ corresponding to the following convex ordering on $\Phi^+$:
\begin{equation} \label{eq:convex order i_0}
\begin{split}
&\ep_i+\ep_j \prec \ep_k-\ep_l,\\
\ep_i+\ep_j \prec \ep_k+& \ep_l \quad \Longleftrightarrow \quad 
\text{$(j>l)$ or $(j=l$, $i>k)$}, \\
\ep_i-\ep_j \prec \ep_k-& \ep_l \quad \Longleftrightarrow \quad
\text{$(i<k)$ or $(i=k$, $j<l)$},
\end{split}
\end{equation}
for $1\leq i<j\leq n$ and $1\leq k<l\leq n$ (see \cite[Section 3.1]{JK19-1}). 
We have 
\begin{equation}\label{eq:beta_k}
\Phi^+=\{\beta_1:=\alpha_{i_1}\prec \beta_2:=s_{i_1}(\alpha_{i_2})\prec  \ldots \prec\beta_N:=s_{i_1}\cdots s_{i_{N-1}}(\alpha_{i_N})\},
\end{equation}
where $\beta_1 = \alpha_n$.
Let $\Phi^+_J=\{\,\ep_i-\ep_j\,|\,1\leq i<j\leq n\,\}$ be the set of positive roots of $\mf l$ and $\Phi^+(J)=\{\,\ep_i+\ep_j\,|\,1\leq i<j\leq n\,\}$ be the set of roots of the nilradical $\mf u$ of the parabolic subalgebra of $\mf g$ associated to ${\mf l}$. 
Then we have
\begin{equation*}
\Phi^+(J)=\{\,\beta_i\,|\,1\leq i\leq M\,\},\quad \Phi_{J}=\{\,\beta_i\,|\,M+1\leq i\leq N\,\},
\end{equation*}
where $M=N/2$.
Let 
\begin{equation} \label{eq:def of B}
\B=\{\,{\bf c}=(c_{\beta_1},\dots,c_{\beta_N})\,|\,c_{\beta_i}\in\Z_+ \ (1\leq i\leq N)\,\}=\Z_+^N.
\end{equation} 
Then $\B$ becomes a crystal isomorphic to that the negative part of $U_q(\mf g)$, which is called the crystal of Lusztig data associated to $\bi_0$. Let
\begin{equation} \label{eq:def of B^J and B_J}
\begin{split}
\B^J &=\left\{\,{\bf c}=(c_{\beta})\in \B\,\big\vert \,c_{\beta}=0 \text{ unless $\beta\in \Phi^{+}(J)$}\,\right\},\\
\B_J &=\left\{\,{\bf c}=(c_{\beta})\in \B\,\big\vert \,c_{\beta}=0 \text{ unless $\beta\in \Phi_J$}\,\right\},
\end{split}
\end{equation}
be the subcrystals of $\B$, where we assume that 
$\te_n{\bf c}=\tf_n{\bf c}={\bf 0}$ with  
$\varepsilon_n({\bf c})=\varphi_n({\bf c})=-\infty$ for ${\bf c}\in \B_J$.
The crystal $\B^J$ can be viewed as a crystal of a quantum nilpotent subalgebra, while $\B_J$ is isomorphic to the crystal of $U^-_q(\mf l)$ as an ${\mf l}$-crystal.

The reduced expression $\bi_0$ enables us to describe the crystal structure on $\B$ very explicitly \cite[Proposition 3.2]{JK19-1}. In particular, we have

\begin{prop}\label{thm:decomposition}\cite[Corollary 3.5]{JK19-1}
The map
\begin{equation} \label{eq:Phi}
\xymatrixcolsep{3pc}\xymatrixrowsep{0pc}\xymatrix{
\B^J \, \otimes \B_J \ar@{->}[r]  & \ \B \\
{\bf c}^J\otimes{\bf c}_{J} \ar@{|->}[r] &  {\bf c} }
\end{equation}
is an  isomorphism of crystals, where ${\bf c}$ is the concatenation of ${\bf c}^J$ and ${\bf c}_J$ given by 
\begin{equation*}
	{\bf c} = ( \LaTeXunderbrace{c_{\beta_1}, \dots, c_{\beta_M}}_{{\bf c}^J}, \LaTeXunderbrace{c_{\beta_{M+1}}, \dots, c_{\beta_N}}_{{\bf c}_J} ).
\end{equation*}
\end{prop}

\begin{rem} \label{rem:uniqueness of i^J}
{\em 
One may find an explicit description of the reduced expression ${\bf i}_0=(i_1,\dots,i_N)$ of $w_0$ corresponding to \eqref{eq:convex order i_0} \cite[Section 3.1]{JK19-1}. 
Recall that there is a one-to-one correspondence between
the reduced expressions of $w_0$ and the convex orderings on $\Phi^+$ \cite{Pap94}.
Then the subexpression $(i_1, \dots, i_M)$ corresponding to the roots of $\mf u$ always appears as the first $M$ entries (up to 2-term braid moves) in any reduced expression of $w_0$  such that the positive roots of ${\mf u}$ precede those of ${\mf l}$ with respect to the corresponding convex ordering.
Here a 2-term braid move means $ij=ji$ for $i,j\in I$ such that $|i-j|>1$.

}	
\end{rem}

\begin{rem}
\mbox{}
{\em 
In \cite[Proposition 3.2]{JK19-1}, we give an explicit description of the operator $\widetilde{f}_i$ on $\B$, which is obtained by applying \cite[Theorem 4.5]{SST1} to the reduced expression ${\bf i}_0$ (see also \cite[Remark 3.3]{JK19-1}).
}
\end{rem}

\begin{rem}
{\em 
There is another reduced expression ${\bf i}$ of $w_0$, which gives a nice combinatorial description of the crystal ${\bf B}_{\bi}$ of Lusztig data associated to ${\bf i}$ \cite{SST1}. Furthermore, there is an explicit construction of {\em isomorphism} from the crystal of {\em marginally large tableaux} of type $D$ to ${\bf B}_{\bi}$ \cite{SST2}, where a marginally large tableau is a combinatorial model for the crystal of the negative part of the quantized enveloping algebra of classical type \cite{HL} ({\em not} the crystal of a highest weight module).

We remark that the reduced expression $\bi$ in \cite{SST2} is not equal to $\bi_0$ here (up to $2$-term braid moves).
It would be interesting to see whether we have similar results as in \cite{JK19-1} with respect to $\bi$.
}
\end{rem}

\subsection{RSK of type $D_n$} \label{subsec:RSK of type D}
Let $\Omega$ be the set of biwords $(\ba,\bb)$ such that
\begin{itemize}
\item[(1)] $\ba=a_1\cdots a_r$ and $\bb=b_1\cdots b_r$ are finite words in $[\ov{n}]$ for some $r\geq 0$,

\item[(2)] $a_i < b_i$ for $1\leq i\leq r$,

\item[(3)] $(a_1,b_1)\leq \cdots \leq (a_r,b_r)$, 
\end{itemize}
where $(a,b)< (c,d)$ if and only if $(a<c)$ or ($a=c$ and $b>d$) for $(a,b)$, $(c,d)$. 
One may naturally identify $\Omega$ with $\B^J$ by regarding $\left|\{\,k\,|\,(a_k,b_k)=(a,b) \,\}\right|$ as the multiplicity of $\epsilon_i+\epsilon_j$ when $(a,b)=(\ov{j},\ov{i})$ with $i<j$.

There is an analogue of RSK due to \cite{B}, which gives a bijection
\begin{equation}\label{eq:kappa}
\xymatrixcolsep{3pc}\xymatrixrowsep{0pc}\xymatrix{
\V  \ \ar@{->}[r] & \ \B^J  \\
\quad T \ar@{|->}[r] & {\bf c}^J(T) }.
\end{equation}

For the reader's convenience, let us briefly recall the bijection \eqref{eq:kappa}. For $\ov{i}\in [\ov{n}]$ and $T\in SST_{[\ov{n}]}(\la^\pi)$ with $\la\in \cP_n$, we denote by $T \leftarrow \ov{i}$ the tableau obtained by applying the Schensted's column insertion of $\ov i$ into $T$ in a reverse way starting from the rightmost column of $T$ so that 
$(T\leftarrow \ov i)\in SST_{[\ov{n}]}(\mu^\pi)$ for some $\mu$ obtained by adding a box in a corner of $\la$.  

Now, let $T \in {SST}(\delta^{\pi})\subset \V$ be given. Then ${\bf c}^J(T)$ is given by the following steps:
\begin{itemize}
	\item[(1)] Let $\ov{x}_1$ be the smallest entry in $T$ such that it is located at the leftmost among such entries and let $T'$ be the tableau obtained from $T$ by removing $\ov{x}_1$.
	
	\item[(2)] Take $\ov{y}_1$ which is the entry in $T$ below $\ov{x}_1$. Let $T''$ be the tableau obtained from $T'$ by applying the inverse of the Schensted's column insertion to $\ov{y}_1$.
	Then we obtain an entry $\ov{z}_1$ such that $(T'' \leftarrow {\ov{z}_1}) = T'$.
	
	\item[(3)] We apply the above steps to $T_1 := T''$ instead of $T$. Then we denote by $\ov{x}_2$ and $\ov{z}_2$ the entries obtained from $T_1$ and let $T_2 := T_1''$. 
	
	\item[(4)] In general, repeat this process for $T_{i+1} := T_i''$ $(i=1, 2, \dots)$ until there is no entries in $T_{i+1}$, and let $\ov{x}_{i+1}$ and $\ov{z}_{i+1}$ be the entries from $T_{i+1}$ by this process.
	\item[(5)] Finally we obtain a biword $({\bf a}, {\bf b})\in \Omega$ given by 
		\begin{equation*}
			\left( 
				\begin{array}{c}
				{\bf a} \\
				{\bf b}
			\end{array} \right) =
			\left( 
				\begin{array}{cccc}
				\ov{x}_1 & \ov{x}_2 & \dots & \ov{x}_{\ell} \\
				\ov{z}_1 & \ov{z}_2 & \dots & \ov{z}_{\ell}
			\end{array} \right)\,,
		\end{equation*}
and define ${\bf c}^J(T)\in \B^J$ to be the one corresponding to $(\ba,\bb)$. 

\end{itemize}

\begin{thm}\label{thm:isomorphism theorem} \cite[Theorem 4.6]{JK19-1}
The bijection in \eqref{eq:kappa} is an isomorphism of crystals.
\end{thm}

\subsection{Embedding of ${\bf KN}_\la$ into $\B$}
Let $\mu\in\cP_n$ be given and let $\epsilon_\mu =\sum_{i=1}^n\mu_i\epsilon_{n-i}$.  Consider 
the map
	\begin{equation} \label{eq:SST to B_J}
		\xymatrixcolsep{3pc}\xymatrixrowsep{0pc}
		\xymatrix{
			 {\bf S}_\mu \ \ar@{->}[r] & {\bf B}_J \otimes T_{\epsilon_\mu} \\
 				T \ar@{|->}[r] & {\bf c}_J(T) \otimes t_{\epsilon_\mu}
		}\,,
	\end{equation}
where ${\bf c}_J(T)$ is the one such that the multiplicity $c_{\epsilon_i-\epsilon_j}$ is equal to the number of $\ov{i}$'s appearing in the $(n-j+1)$th row of $T$ for $1 \le i < j \le n$.

\begin{prop} \label{prop:SST to B_J}	
The map \eqref{eq:SST to B_J} is an embedding of crystals.
\end{prop}
\pf It is a well-known fact that the map is an embedding of ${\mf l}$-crystals (cf.~\cite{K18-1}).
By definition of $\te_n$ and $\tf_n$ on ${\bf S}_\mu$ and $\B_J$, it becomes a morphism of $\mf g$-crystals. 
\qed

\begin{cor} \label{prop:T to B^J otimes B_J}
For $\lambda \in \mc{P}_n$, the map  
\begin{equation} \label{eq:phi}
\xymatrixcolsep{3pc}\xymatrixrowsep{0pc}
	\xymatrix{
{\bf V}_{\lambda} \otimes T_{r\omega_n}  \ \ar@{->}[r] & {\bf B}^J\otimes {\bf B}_J \otimes T_{\omega_{\lambda}}\\
 \ S \otimes T \otimes t_{r\omega_n}  \ar@{|->}[r] & {\bf c}^J(S) \otimes {\bf c}_J(T) \otimes t_{\omega_{\lambda}}
}\,,
\end{equation}
is an embedding of crystals, where $r = ( \lambda, \omega_n )$.
\qed
\end{cor}

 We are now in a position to state the main result in this paper. 

\begin{thm} \label{thm:main}
For $\lambda \in \mc{P}_n$, 
we have an embedding of crystals given by
\begin{center} 
\begin{tikzpicture}
  \matrix (m) [matrix of math nodes,row sep=3em,column sep=4em,minimum width=2em]
  {
     {\bf KN}_{\lambda} & & {\bf B} \otimes T_{\omega_{\lambda}} \\
     {\bf T}_{\lambda} & {\bf V}_{\lambda} \otimes T_{r\omega_n}  & {\bf B}^J \otimes {\bf B}_J \otimes T_{\omega_{\lambda}}\\};
  \path[-stealth]
    (m-1-1) edge node [left] {\eqref{eq:Psi}} (m-2-1)
            edge node [above] {$\Xi_{\lambda}$} (m-1-3)
    (m-2-1) edge node [below] {\eqref{eq:chi}} (m-2-2)
    (m-2-2) edge node [below] {\eqref{eq:phi}} (m-2-3)
    (m-2-3) edge node [right] {\eqref{eq:Phi}} (m-1-3);
\end{tikzpicture}
\end{center}
\end{thm}
\pf
It follows from Theorems \ref{thm:KN to spinor}, \ref{thm:separation}, Proposition \ref{thm:decomposition} and Corollary \ref{prop:T to B^J otimes B_J}.
\qed

\begin{ex}
{\em 
Set $n=5$ and $\lambda = \left(\frac{5}{2}, \frac{3}{2}, \frac{3}{2}, \frac{1}{2}, -\frac{1}{2}\right)$.
Let $T \in {\bf KN}_{\lambda}$ be given by
\begin{equation*}
\hspace{6.5cm} T = \hspace{-6.5cm} 
\begin{split}
\ytableausetup {mathmode, boxsize=0.88em} 
\begin{ytableau}
	\none & \none & \none \\
	\none & \none & \tl{$2$} \\
	\none & \none & \tl{$3$} \\
	\none & \tl{$4$} & \tl{$\ov{5}$} \\
	\none & \tl{$5$} & \tl{$\ov{4}$} \\
	\tl{$\ov{5}$} & \tl{$\ov{1}$} & \tl{$\ov{1}$} \\
\end{ytableau}
\end{split}\,.
\end{equation*}

By \eqref{eq:Psi} and \eqref{eq:chi} (see Example \ref{ex:separation for negative case}), we have 
\begin{equation*}
\hspace{1.5cm} \Psi_{\lambda}(T) = \hspace{-2.3cm}
\begin{split}
\ytableausetup {mathmode, boxsize=1.0em} 
& \begin{ytableau}
\none & \none & \none & \none & \none & \none & \none & \none & \none & \none & \none & \none \\
\none &\none & \none & \none & \none & \none & \none & \none & \none & \none & \none & \none  \\
\none &\none & \none & \none & \none & \none & \none & \tl{$\ov{5}$} & \none & \none & \none \\
\none &\none & \tl{$\ov{5}$} & \none & \none & \tl{$\ov{3}$} & \none & \tl{$\ov{4}$} & \none\\
\none &\none & \tl{$\ov{4}$} & \none[\mathrel{\raisebox{-0.7ex}{$\scalebox{0.45}{\dots\dots}$}}] & \none & \tl{$\ov{1}$} & \none[\mathrel{\raisebox{-0.7ex}{$\scalebox{0.45}{\dots\dots}$}}] & \tl{$\ov{1}$} &  \none& \none \\
\none &\tl{$\ov{5}$} & \none & \none & \tl{$\ov{2}$} & \none & \none & \none & \none & \none & \none  \\
\none &\tl{$\ov{3}$} & \none & \none & \tl{$\ov{1}$} & \none & \none & \none & \none & \none & \none  \\
\none &\tl{$\ov{2}$} & \none & \none & \none & \none & \none & \none & \none & \none & \none & \none  \\
\none &\tl{$\ov{1}$} & \none & \none & \none & \none & \none & \none & \none & \none & \none & \none  \\
\none & \none & \none & \none & \none & \none & \none & \none & \none & \none & \none & \none  \\
\end{ytableau}
\end{split}
\hspace{-1cm} \overset{\eqref{eq:chi}}{\xrightarrow{\hspace*{1.2cm}}} \quad
\begin{split}
\ytableausetup {mathmode, boxsize=1.0em} 
& \begin{ytableau}
\none & \none &\none & \none & \none & \none & \none & \none & \none \\
\none & \none &\none & \none & \none & \none & \tl{$\ov{5}$} & \none & \none\\
\none & \none &\none & \none & \none & \none & \tl{$\ov{4}$} & \none & \none\\
\none & \none &\none & \none & \none & \tl{$\ov{5}$} & \tl{$\ov{3}$} & \none & \none\\
\none & \none[\mathrel{\raisebox{-0.7ex}{$\scalebox{0.45}{\dots\dots\dots\dots}$}}] &\none & \none & \none & \tl{$\ov{1}$} & \tl{$\ov{1}$} & \none[\mathrel{\raisebox{-0.7ex}{$\scalebox{0.45}{\dots\dots\dots\dots}$ \, }}] & \none\\
\none & \none &\tl{$\ov{5}$} & \tl{$\ov{4}$} & \tl{$\ov{2}$} & \none & \none & \none & \none\\
\none & \none &\tl{$\ov{3}$} & \tl{$\ov{1}$} & \none & \none & \none & \none & \none\\
\none & \none &\tl{$\ov{2}$} & \none & \none & \none & \none & \none & \none\\
\none & \none &\tl{$\ov{1}$} & \none & \none & \none & \none & \none & \none\\
\none
\end{ytableau}
\end{split}\, \quad
\end{equation*} 
Note that $\sigma = (\, -,\, +,\, +,\, -,\,\, \cdot\,\,)$ and $\ov{\sigma} = (\,-,\,\,\cdot\,\,,\,+,\,\,\cdot\,\,,\,\cdot\,\,)$.
Also
\begin{equation*}
	\sigma^{\rm red} = (\, -,\, +,\,\, \cdot\,\,,\, \cdot\,\,,\,\, \cdot\,\,), \quad \ov{\sigma}^{\rm red} = (\,-,\,\,\cdot\,\,,\,+,\,\,\cdot\,\,,\,\cdot\,\,)
\end{equation*}
(cf. Lemma \ref{lem:invariant nth signature}).
%

Put ${\bf T} = \Psi_{\lambda}(T)$.
Then
\begin{equation*}
\hspace{3.5cm} {\bf T}^{\texttt{body}} = \hspace{-4.3cm}
\begin{split} 
\ytableausetup {mathmode, boxsize=1.0em} 
\begin{ytableau}
\none & \tl{$\ov{5}$} \\
\none & \tl{$\ov{4}$} \\
\tl{$\ov{5}$} & \tl{$\ov{3}$} \\
\tl{$\ov{1}$} & \tl{$\ov{1}$} \\
\none & \none \\
\none & \none \\
\none & \none \\
\none & \none 
\end{ytableau}
\end{split}\,\,, \quad \quad \quad
{\bf T}^{\texttt{tail}} = \,\,
\begin{split} 
\ytableausetup {mathmode, boxsize=1.0em} 
\begin{ytableau}
\none & \none & \none \\
\none & \none & \none \\
\none & \none & \none \\
\none & \none & \none \\
\tl{$\ov{5}$} & \tl{$\ov{4}$} & \tl{$\ov{2}$} \\
\tl{$\ov{3}$} & \tl{$\ov{1}$} & \none \\
\tl{$\ov{2}$} & \none & \none \\
\tl{$\ov{1}$} & \none & \none 
\end{ytableau}
\end{split}\,.
\end{equation*}

Let us recall \eqref{eq:def of B} and \eqref{eq:def of B^J and B_J}.
For simplicity, we use the notation in \cite[Section 3.2]{JK19-1} for ${\bf c}^J({\bf T}^{\texttt{body}})$ associated to $\Phi^+(J)$ with the convex order \eqref{eq:convex order i_0}, that is, we identify $(c_{\beta_1}, \dots, c_{\beta_{10}}) \in {\bf B}^J$ with  
\begin{equation*}
\begin{split}
\scalebox{1}{\xymatrixcolsep{-0.2pc}\xymatrixrowsep{0.3pc}\xymatrix{
& & & & c_{\beta_4} & & & & \\
& & & c_{\beta_3} & & c_{\beta_7} & & & \\
& & c_{\beta_2} & & c_{\beta_6} & & c_{\beta_9} & \\  
& c_{\beta_1} & & c_{\beta_5} & & c_{\beta_8} & & c_{\beta_{10}}}}
\end{split}\,.
\end{equation*}
Here $\beta_1 = \alpha_5$.
Similarly, we use the above notation for ${\bf c}_J({\bf T}^{\texttt{tail}})$ with respect to $\Phi_J$ and the convex order \eqref{eq:convex order i_0}, that is, 
we identify $(\beta_{11}, \dots, \beta_{20}) \in {\bf B}_J$ with 
\begin{equation*}
\begin{split}
\scalebox{1}{\xymatrixcolsep{-0.3pc}\xymatrixrowsep{0.3pc}\xymatrix{
& & & & c_{\beta_{14}} & & & & \\
& & & c_{\beta_{13}} & & c_{\beta_{17}} & & & \\
& & c_{\beta_{12}} & & c_{\beta_{16}} & & c_{\beta_{19}} & \\  
& c_{\beta_{11}} & & c_{\beta_{15}} & & c_{\beta_{18}} & & c_{\beta_{20}}}}
\end{split}\,.
\end{equation*}
Here $\beta_{11} = \alpha_1$, $\beta_{15} = \alpha_2$, $\beta_{18} = \alpha_3$ and $\beta_{20} = \alpha_4$
(cf. \cite[Example 3.1]{JK19-1}).

Now we find ${\bf c}^J\,({\bf T}^{\texttt{body}})$ by the steps (1)--(5) in Section \ref{subsec:RSK of type D} as follows.

\begin{equation*}
\begin{split}
\ytableausetup {mathmode, boxsize=1.0em} 
\begin{ytableau}
\none & \tl{$\ov{5}$} \\
\none & \tl{$\ov{4}$} \\
\tl{$\ov{\bf 5}$} & \tl{$\ov{3}$} \\
\tl{$\ov{\bf 1}$} & \tl{$\ov{1}$} \\
\none & \none
\end{ytableau}
\end{split} \quad \overset{}{\xrightarrow{\hspace*{0.5cm}}} \,
\begin{split}
\ytableausetup {mathmode, boxsize=1.0em} 
\begin{ytableau}
\none & \tl{$\ov{\bf 5}$} \\
\none & \tl{$\ov{\bf 4}$} \\
\none & \tl{$\ov{3}$} \\
\none & \tl{$\ov{1}$} \\
\none & \none
\end{ytableau}
\end{split}\,\,,\quad
\begin{split}
	\scalebox{0.8}{$\left( \begin{array}{c}
		 \ov{5} \\
		 \ov{1}
	\end{array}
\right)$}
\end{split} \quad \overset{}{\xrightarrow{\hspace*{0.5cm}}} \,
\begin{split}
\ytableausetup {mathmode, boxsize=1.0em} 
\begin{ytableau}
\none & \none \\
\none & \none \\
\none & \tl{$\ov{\bf 3}$} \\
\none & \tl{$\ov{\bf 1}$} \\
\none & \none
\end{ytableau}
\end{split}\,\,,\quad
\begin{split}
	\scalebox{0.8}{$\left( \begin{array}{cc}
		 \ov{5} & \ov{5} \\
		 \ov{4} & \ov{1}
	\end{array}
\right)$}
\end{split} \quad \overset{}{\xrightarrow{\hspace*{0.5cm}}} \,\,\,\,\,
\begin{split}
\emptyset
\end{split}\,\,,\quad
\begin{split}
	\scalebox{0.8}{$\left( \begin{array}{ccc}
		 \ov{5} & \ov{5} & \ov{3} \\
		 \ov{4} & \ov{1} & \ov{1}
	\end{array}
\right)\,.
$}
\end{split}
\end{equation*}
Thus we have
\begin{equation*}
\begin{split}
\ytableausetup {mathmode, boxsize=1em} 
\begin{ytableau}
	\none & \tl{$\ov{5}$} \\
	\none & \tl{$\ov{4}$} \\
	\tl{$\ov{5}$} & \tl{$\ov{3}$} \\
	\tl{$\ov{1}$} & \tl{$\ov{1}$} \\
\end{ytableau}
\end{split} \quad  \overset{\eqref{eq:kappa}}{\xrightarrow{\hspace*{1.2cm}}} \,
\begin{split}
\scalebox{0.7}{\xymatrixcolsep{-0.2pc}\xymatrixrowsep{0.3pc}\xymatrix{
& & & & 1 & & & & \\
& & & 0 & & 0 & & & \, \in \B^J \\
& & 0 & & 0 & & 1 & \\  
& 1 & & 0 & & 0 & & 0}}
\end{split}
\end{equation*}
(cf. Example \cite[Example 4.5]{JK19-1}).
Next by definition \eqref{eq:SST to B_J}, we have ${\bf c}_J({\bf T}^{\texttt{tail}})$ as follows.


\begin{equation*}
\begin{split}
\quad \ytableausetup {mathmode, boxsize=1em} 
\begin{ytableau}
	\tl{$\ov{5}$} & \tl{$\ov{4}$} & \tl{$\ov{2}$} \\
	\tl{$\ov{3}$} & \tl{$\ov{1}$} & \none \\
	\tl{$\ov{2}$} & \none & \none \\
	\tl{$\ov{1}$} & \none & \none 
\end{ytableau}
\end{split} \quad  \overset{\eqref{eq:SST to B_J}}{\xrightarrow{\hspace*{1.2cm}}} \,
\begin{split}
\scalebox{0.7}{\xymatrixcolsep{-0.2pc}\xymatrixrowsep{0.3pc}\xymatrix{
& & & & 0 & & & & \\
& & & 1 & & 1 & & & \, \in \B_J\\
& & 0 & & 0 & & 0 & \\  
& 1 & & 1 & & 1 & & 1}}
\end{split}\,. \quad \quad 
\end{equation*}

Hence we obtain the Lusztig data for the KN tableau $T$ associated to ${\bf i}_0$ (cf. Proposition \ref{thm:decomposition}), that is,
\begin{equation*}
	\Xi_{\lambda}(T) = (1,\, 0,\, 0,\, 1,\, 0,\, 0,\, 0,\, 0,\, 1,\, 0,\, 1,\, 0,\, 1,\, 0,\, 1,\, 0,\, 1,\, 1,\, 0,\, 1) \otimes t_{\omega_{\lambda}}.
\end{equation*}
}	
\end{ex}

\end{document}